\documentclass[fleqn,12pt,a4paper]{article}

\usepackage{ddwmacs,latexsym,amssymb}
\usepackage[dvips,xdvi]{graphics,color}
\usepackage{epsfig}

\title{%
  Automatic Evaluation of the Links--Gould Invariant for all
  Prime Knots of up to $10$ Crossings
}

\author{%
  David~~De Wit%
  \footnote{
    RIMS, Kyoto University 606-8502, Japan.
    \texttt{ddw@kurims.kyoto-u.ac.jp}
  }
}

\begin{document}

\maketitle

\begin{center}
  \begin{picture}(0,0)%
\epsfig{file=graphics/BraidArt.pstex}%
\end{picture}%
\setlength{\unitlength}{0.00083300in}%
\begingroup\makeatletter\ifx\SetFigFont\undefined
\def\x#1#2#3#4#5#6#7\relax{\def\x{#1#2#3#4#5#6}}%
\expandafter\x\fmtname xxxxxx\relax \def\y{splain}%
\ifx\x\y   
\gdef\SetFigFont#1#2#3{%
  \ifnum #1<17\tiny\else \ifnum #1<20\small\else
  \ifnum #1<24\normalsize\else \ifnum #1<29\large\else
  \ifnum #1<34\Large\else \ifnum #1<41\LARGE\else
     \huge\fi\fi\fi\fi\fi\fi
  \csname #3\endcsname}%
\else
\gdef\SetFigFont#1#2#3{\begingroup
  \count@#1\relax \ifnum 25<\count@\count@25\fi
  \def\x{\endgroup\@setsize\SetFigFont{#2pt}}%
  \expandafter\x
    \csname \romannumeral\the\count@ pt\expandafter\endcsname
    \csname @\romannumeral\the\count@ pt\endcsname
  \csname #3\endcsname}%
\fi
\fi\endgroup
\begin{picture}(1159,2438)(1147,-3077)
\end{picture}

\end{center}

\begin{abstract}
  \noindent
  This paper describes a method for the automatic evaluation of the
  Links--Gould two-variable polynomial link invariant ($LG$) for any
  link, given only a braid presentation. This method is currently
  feasible for the evaluation of $LG$ for links for which we have a
  braid presentation of string index at most $5$.  Data are presented
  for the invariant, for all prime knots of up to $10$ crossings and
  various other links.  $LG$ distinguishes between these links, and
  also detects the chirality of those that are chiral. In this sense,
  it is more sensitive than the well-known two-variable HOMFLY and
  Kauffman polynomials.  When applied to examples which defeat the
  HOMFLY invariant, interestingly, $LG$ `almost' fails.
  The automatic method is in fact applicable to the
  evaluation of any such state sum invariant for which an appropriate R
  matrix and cap and cup matrices have been determined.
\end{abstract}


\vfill

\pagebreak

\section{Background}

The Links--Gould two-variable polynomial invariant of oriented links
($LG$) is based on the $(0,0|\alpha)$ family of $4$ dimensional highest
weight representations of the quantum superalgebra $U_q[gl(2|1)]$.  It
was originally described by Jon Links and Mark Gould in
\cite{LinksGould:92b}, although that work did not pursue its evaluation
for want of an efficient method.  In \cite{DeWitKauffmanLinks:99b},
Louis Kauffman, Jon Links and I
described a state model for its evaluation.%
\footnote{
  In conjunction with this work, the reader should peruse
  \cite{DeWitKauffmanLinks:99b}, of which a more thorough treatment is
  found in \cite{DeWit:98} (available from the author).
}
A very precise description of state models is contained within
\cite{AkutsuDeguchiOhtsuki:92}.

In \cite{DeWitKauffmanLinks:99b}, we demonstrated that $LG$ does not
detect the noninvertibility or mutation of links. Furthermore, we
provided evaluations of it for a small set of sundry interesting links,
in particular for the (prime) knots $9_{42}$ and $10_{48}$, these two
being infamous as the first chiral knots whose chirality is undetected
by the HOMFLY invariant.  Whilst it was discovered that $LG$ detected
the chirality of these two examples, the method used to evaluate it was
labour-intensive, and it was impracticable to evaluate $LG$ for
\emph{all} the knots up to $9_{42}$. Thus, we did not know whether $LG$
was in fact more sensitive than the HOMFLY invariant as $LG$ might have
failed to detect the chirality of some chiral knot before $9_{42}$.

Here, we improve upon that situation by describing a more sophisticated
method for evaluating $LG$.  This method requires as its input a braid
presentation $\beta$ for the (closed) link $\hat{\beta}$ in question.
Using it, we evaluate $LG$ for \emph{all} the (prime) knots of up to
$10$ crossings.


\section{Changes of Variable and Notation}

Importantly, we make a change of variable from
\cite{DeWitKauffmanLinks:99b}. In that work, the underlying
$\alpha$-parametric family of representations of $U_q[gl(2|1)]$ led us
to an invariant which manifested as an integer-coefficient `polynomial'
in variable $q$ with exponents being integer-linear functions of
$\alpha$.  Defining $p\defeq q^{\alpha}$ ensured that our invariant was
a two-variable (Laurent) polynomial in variables $q$ and $p$ (actually
$q^2$ and $p^2$).  We showed that $LG$ would be unchanged under the
transformation $\alpha\mapsto -\alpha-1$ (and hence unable to detect
noninvertibility). This meant that it satisfied the symmetry
$LG_K(q,p)=LG_K(q,q^{-1}p^{-1})$.

Here, we instead define $p\defeq q^{\alpha+1/2}$, which changes the
symmetry to the simpler $LG_K(q,p)=LG_K(q,p^{-1})$.  In these new
variables $q$ and $p$, the representation of the braid generator and
the caps and cups take on more symmetric forms than they do using our
previous variables.  These new variables are thus a more natural choice
for our invariant.

\vfill

\pagebreak

The tensor $\psi(\overline{\sigma})$ representing the
braid generator in the state model is now:
\def\scs{\scriptstyle}
\begin{eqnarray*}
  \hspace{-20pt}
  \m{[}{@{}*{3}{*{3}{c@{}@{}}c|}*{3}{c@{}@{}}c@{}}
    {\scs p^{-2}q} &.&.&.&.&.&.&.&.&.&.&.&.&.&.&.\\
    .&.&.&.&{\scs p^{-1} q^{1/2}} &.&.&.&.&.&.&.&.&.&.&.\\
    .&.&.&.&.&.&.&.&{\scs p^{-1} q^{1/2}} &.&.&.&.&.&.&.\\
    .&.&.&.&.&.&.&.&.&.&.&.&{\scs 1} &.&.&.\\
    \hline
    .&{\scs p^{-1} q^{1/2}} &.&.&{\scs p^{-2}q-1} &.&.&.&.&.&.&.&.&.&.&.\\
    .&.&.&.&.&{\scs -1} &.&.&.&.&.&.&.&.&.&.\\
    .&.&.&.&.&.&.&.&.&{\scs -q} &.&.&{\scs -q^{1/2}Y} &.&.&.\\
    .&.&.&.&.&.&.&.&.&.&.&.&.&{\scs p q^{1/2}} &.&.\\
    \hline
    .&.&{\scs p^{-1} q^{1/2}}&.&.&.&.&.&{\scs p^{-2}q-1}&.&.&.&.&.&.&.\\
    .&.&.&.&.&.&{\scs -q}&.&.&{\scs q^2-1}&.&.&{\scs q^{3/2}Y}&.&.&.\\
    .&.&.&.&.&.&.&.&.&.&{\scs -1}&.&.&.&.&.\\
    .&.&.&.&.&.&.&.&.&.&.&.&.&.&{\scs p q^{1/2}}&.\\
    \hline
    .&.&.&{\scs 1}&.&.&{\scs -q^{1/2}Y}&.&.&{\scs q^{3/2}Y}&.&.&{\scs q Y^2}&.&.&.\\
    .&.&.&.&.&.&.&{\scs p q^{1/2}} &.&.&.&.&.&{\scs p^2q-1} &.&.\\
    .&.&.&.&.&.&.&.&.&.&.&{\scs p q^{1/2}} &.&.&{\scs p^2 q -1} &.\\
    .&.&.&.&.&.&.&.&.&.&.&.&.&.&.&{\scs p^2 q}
  \me,
\end{eqnarray*}
where we have used the shorthand:
\begin{eqnarray*}
  Y
  \defeq
  \sqrt{
    ( p q^{1/2} - p^{-1} q^{-1/2} )
    ( p q^{-1/2} - p^{-1} q^{1/2} )
  }
  =
  \sqrt{
    p^{+2} + p^{-2} - q^{+1} - q^{-1}
  };
\end{eqnarray*}
note that this definition of $Y$ differs from that in our previous
work. This $Y$ is more symmetric than our previous one -- it is
invariant under $q\mapsto q^{-1}$ (and hence $p\mapsto p^{-1}$).
In our current variables,
$\psi(\overline{\sigma}^{-1})$
is obtainable from $\psi(\overline{\sigma})$ by a
twisting transformation combined with the mapping $q\mapsto q^{-1}$.
Explicitly:
\begin{eqnarray*}
  {\[
    \psi(\overline{\sigma}^{-1})
  \]}^{a~c}_{b~d}
  =
  {\[
    \psi(\overline{\sigma})
  \]}^{c~a}_{d~b}
  \quad
  \textrm{after~mapping}
  \quad
  q \mapsto q^{-1}.
\end{eqnarray*}

This transformation was opaque using our previous variables. It
effectively illustrates the `uniqueness' of our R matrix by
demonstrating explicitly a unitary transformation relating two
candidates. It is in fact also visible in the R matrix used for the
bracket polynomial \cite{Kauffman:87b,Kauffman:97a}, i.e.  that
R matrix originating in the $2$ dimensional vector representation of
$U_q[sl(2)]$. In that case, the same twist is used, and the mapping is
$A\mapsto A^{-1}$.

\pagebreak

To describe the caps ($\Omega^{\pm}$) and cups ($\mho^{\pm}$), we
introduce another slight change of notation from
\cite{DeWitKauffmanLinks:99b}; to whit, we exchange $\Omega^{+}$ with
$\Omega^{-}$. This notation ensures that the superscript ``$+$''
(respectively ``$-$'') indicates a sector of an anticlockwise
(respectively clockwise) arc, and is consistent with the introduction
of an additional notation: we shall use $C^{\pm}$ to denote vertically
oriented arcs defined by the composition of horizontally oriented arcs.
We shall call these new diagram components \emph{left handles}.%
\footnote{For our current purposes, we have no need for right handles.}
Our current collection of arcs is shown in \figref{CapsCupsHandles}.

\begin{figure}[htbp]
  \begin{center}
    \input{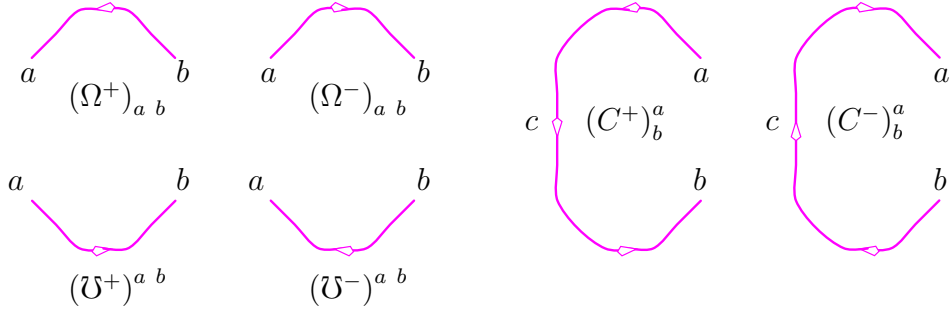}
    \caption{
      Caps, cups and (left) handles.
    }
    \figlabel{CapsCupsHandles}
  \end{center}
\end{figure}

The matrices corresponding to the left handles may be described in
terms of the caps and cups:
\begin{eqnarray*}
  {\(C^{\pm}\)}^{a}_{b}
  =
  {\(\Omega^{\pm}\)}_{c~a}
  {\(\mho^{\pm}\)}^{c~b}.
\end{eqnarray*}
These definitions are generally applicable to any state model.  In
\cite{DeWitKauffmanLinks:99b}, we determined appropriate caps and cups
for our particular state model.  Using our new diagram conventions, in
terms of our new variables that information becomes:
\begin{eqnarray*}
  \Omega^{-}
  \eq
  I_4
  \\
  \Omega^{+}
  \eq
  p^{-2} \; \textrm{diag} \{ q^{+1}, - q^{+1}, - q^{-1}, q^{-1} \}
  \\
  \mho^{+}
  \eq
  {\(\Omega^{-}\)}^{-1}
  =
  I_4
  \\
  \mho^{-}
  \eq
  {\(\Omega^{+}\)}^{-1}
  =
  p^{+2} \; \textrm{diag} \{ q^{-1}, - q^{-1}, - q^{+1}, q^{+1} \}
  \\
  C^{+}
  \eq
  p^{-2} \; \textrm{diag} \{ q^{+1}, - q^{+1}, - q^{-1}, q^{-1} \}
  \\
  C^{-}
  \eq
  p^{+2} \; \textrm{diag} \{ q^{-1}, - q^{-1}, - q^{+1}, q^{+1} \}.
\end{eqnarray*}
Again, note that these matrices possess greater symmetries than
they do using our previous variables.

Lastly, by inspection, it turns out that our new evaluations of $LG$
are actually Laurent polynomials in $q$ and $p^2$, so we can introduce
the further notational artifice $P\defeq p^2$.

\vfill


\section{The Algorithm}

The present algorithm used for evaluating state model invariants
abandons the prior drawing of a quasi-Morse diagram for the link (as
used in \cite{DeWitKauffmanLinks:99b}), instead recognising that a
closed braid corresponding to the link is already in Morse form.  If
the braid form is used to represent a link, then the previous method
may be simplified to eliminate all the auxiliary abstract tensors, as
we need only consider crossings with upwardly-pointing arrows,
represented by the tensors $\sigma$ and $\sigma^{-1}$.

We describe the process for a link $K=\hat{\beta}$, corresponding to an
$n$ string braid word $\beta$ (i.e. $\beta\in B_n$). We regard $\beta$
as an oriented $(n,n)$ tangle, and $\hat{\beta}$ as being obtained by
the (vertical) closures of all the strings of $\beta$.  Recall that for
our purposes, we actually require a $(1,1)$ tangle obtained by closing
all but one of the strings of $\beta$.
Firstly, we construct $Z$, a rank $2n$ tensor corresponding to
$\beta$.  Secondly, we construct $T$, a rank $2$ tensor obtained by
contracting all but $2$ of the indices of $Z$. In contrast to our
previous method, at each stage of the construction this method uses the
data contained in $\sigma$ -- the tensors $Z$ and $T$ are concrete, not
abstract.  The indices of all our tensors run from $1$ to $M$, where
$M$ is the dimension of the underlying representation (in our case
$M=4$, but the process is more generally applicable to any state sum
invariant).  Thus, $Z$ has $M^{2n}$ entries, and $T$ is an
$M\times M$ matrix.



We shall build $Z$ by an $m$ stage accretionary process corresponding
to the definition of $\beta$ as a string of $m$ braid generators.  At
each stage, the intermediate tensor will also be called $Z$.
Initially, $Z$ corresponds to the trivial $(n,n)$ tangle, so we
initialise it to be the identity rank $2n$ tensor $I_{M}^{\otimes n}$,
where $I_{M}$ is the identity rank $2$ tensor (i.e. the usual
$M\times M$ identity matrix).
Say that our braid $\beta\in B_n$ may be written as:
\begin{eqnarray*}
  \beta
  =
  \sigma^{i_1}_{j_1}
  \sigma^{i_2}_{j_2}
  \cdots
  \sigma^{i_m}_{j_m},
  \qquad
  \qquad
  i_k = \pm 1,
  \quad
  1 \leqslant j_k < n,
  \qquad
  k = 1, \dots, m,
\end{eqnarray*}
where $\sigma^{\pm 1}_{j_k}$ indicates that letter $k$ in $\beta$ is
a positive (respectively negative) crossing of strings $j_k$ and
$j_k+1$.  Where $Z$ corresponds to an $(n,n)$ tangle representing a
part $\overline{\beta}$ of $\beta$, the accretion of a letter $\sigma_j$
to $\overline{\beta}$ corresponds to the product of the tensors $Z$ and
$\sigma_j$ according to \figref{AccreteXtoZ}.%
\footnote{%
  These and the following comments of course apply equally to
  $\sigma^{-1}$.
}
Algebraically (Einstein summation convention used throughout):
\begin{eqnarray*}
  {(Z\sigma_j)}^{a_1~\cdots~a_n}_{b_1~\cdots~b_n}
  =
  {(Z)}^{
    a_1~\cdots~a_{j-1}~c_j~c_{j+1}~a_{j+2}~\cdots~a_n
  }_{
    b_1~\cdots~b_{j-1}~b_j~b_{j+1}~b_{j+2}~\cdots~b_n
  }
  \cdot
  {(\sigma_j)}^{a_j~a_{j+1}}_{c_j~c_{j+1}}.
\end{eqnarray*}
More generally, we may accrete any appropriate object in this fashion,
and thus, we may cheat somewhat in the accretionary process: say
$\sigma_j$ actually repeats $N$ times consecutively within $\beta$.
Having computed the auxiliary rank $4$ tensor $\sigma_j^N$ (described
in \cite{DeWitKauffmanLinks:99b}), we may then accrete $\sigma_j^N$ to
$Z$ in one, rather than $N$ stages, which improves efficiency in a small
but practical manner.

\begin{figure}[htbp]
  \begin{center}
    \input{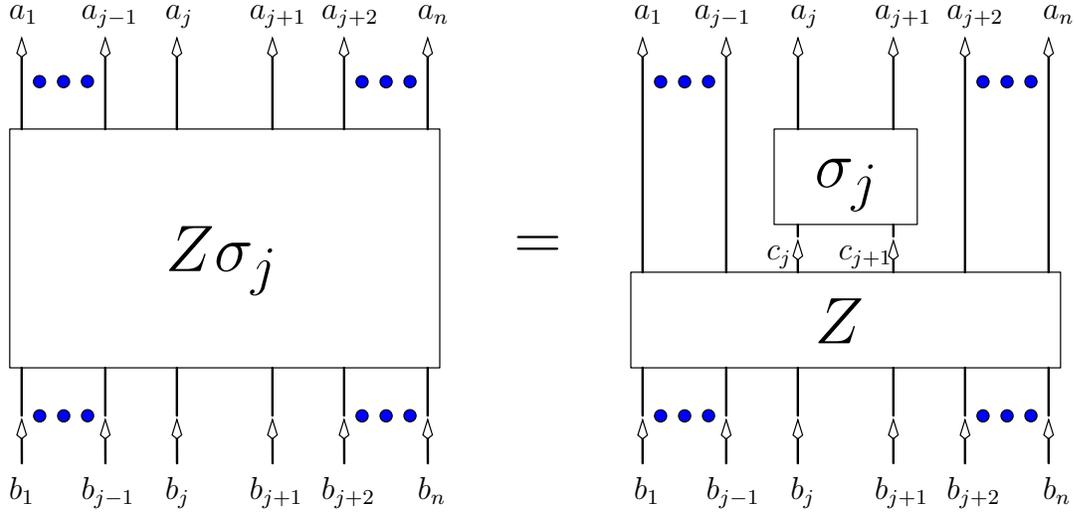}
    \caption{
      Accreting a braid generator to position $j$ of an $(n,n)$ tangle.
    }
    \figlabel{AccreteXtoZ}
  \end{center}
\end{figure}



The tensor $T$ corresponds to the $(1,1)$ tangle obtained by the
closure of all but one of the strings of $\beta$ -- we shall choose the
rightmost one for convenience. This process is depicted in
\figref{ContractZ}.

\begin{figure}[htbp]
  \begin{center}
    \input{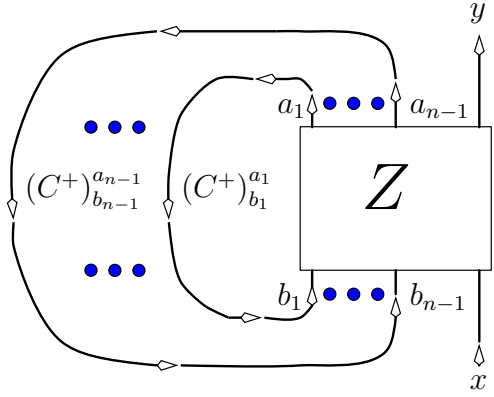}
    \caption{
      Closure of all but one pair of strings of an $(n,n)$ tangle.
    }
    \figlabel{ContractZ}
  \end{center}
\end{figure}

Thus, $T$ is constructed by contraction of all but two of the indices
of $Z$ (preferably done `string-at-a-time'), viz:
\begin{eqnarray*}
  {(T)}^y_x
  \defeq
  {(Z)}^{a_1~\cdots~a_{n-1}~y}_{b_1~\cdots~b_{n-1}~x}
  \cdot
  \prod_{j=1}^{n-1} {(C^{+})}^{a_j}_{b_j}.
\end{eqnarray*}


\section{Feasibility and Implementation}

Where our state model is based on a underlying representation of
dimension $M$, and we are evaluating an invariant based on a braid
presentation of $n$ strings, the tensor $Z$ contains $M^{2n}$
components (see \tabref{ComponentsofZ}), each of which may become a
possibly large polynomial expression (to be fair, it may remain zero
throughout).  Thus, the algorithm requires storage at least
proportional to $M^{2n}$, which is a serious limitation.

\begin{table}[ht]
  \centering
  \begin{tabular}{||c||*{4}{c|}|}
    \hline\hline
                   &                     &                     &
                                         &                     \\[-4mm]
    Algebra        & $U_q[sl(2)]$        & $U_q[gl(2|1)]$      &
                     $U_q[gl(3|1)]$      & $U_q[gl(4|1)]$      \\[1mm]
    \hline
                   &                     &                     &
                                         &                     \\[-4mm]
    Representation & $(1,0)$             & $(0,0|\alpha)$      &
                     $(0,0,0|\alpha)$    & $(0,0,0,0|\alpha)$  \\[1mm]
    \hline
                   &                     &                     &
                                         &                     \\[-4mm]
                   & $M=2$               & $4$                 &
                     $8$                 & $16$                \\[1mm]
    \hline\hline
                   &                     &                     &
                                         &                     \\[-4mm]
    $n=2$  & $16$                & $256$               &
                     $4096$              & $6.6 \times 10^{4}$ \\[1mm]
    \hline
                   &                     &                     &
                                         &                     \\[-4mm]
    $3$            & $64$                & $4096$               &
                     $2.6 \times 10^{5}$ & $1.7 \times 10^{7}$ \\[1mm]
    \hline
                   &                     &                     &
                                         &                     \\[-4mm]
    $4$            & $256$               & $6.6 \times 10^{4}$ &
                     $1.7 \times 10^{7}$ & $\cdot$             \\[1mm]
    \hline
                   &                     &                     &
                                         &                     \\[-4mm]
    $5$            & $1024$              & $1.0 \times 10^{6}$ &
                     $\cdot$             & $\cdot$             \\[1mm]
    \hline
                   &                     &                     &
                                         &                     \\[-4mm]
    $6$            & $4096$              & $1.7 \times 10^{7}$ &
                     $\cdot$             & $\cdot$             \\[1mm]
    \hline
                   &                     &                     &
                                         &                     \\[-4mm]
    $\vdots$       & $\vdots$            & $\cdot$             &
                     $\cdot$             & $\cdot$ \\[1mm]
    \hline
                   &                     &                     &
                                         &                     \\[-4mm]
    $10$           & $1.0 \times 10^{6}$ & $\cdot$             &
                     $\cdot$             & $\cdot$             \\[1mm]
    \hline\hline
  \end{tabular}
  \caption{%
    $M^{2n}$ for various state models and braid indices.
  }
  \tablabel{ComponentsofZ}
\end{table}

Experiments demonstrate that for $LG$, a current practical limit for
$M^{2n}$ is $10^{6}$, thus evaluation of $LG$ is viable for all links
where the braid presentation to hand has at most $5$ strings. All but
five of the prime knots to $10_{166}$ have braid index (see
\secref{BraidPresentations}) at most $5$; the five outstanding ones
have braid index $6$. Thus, our algorithm is applicable to almost all
of the prime knots to $10_{166}$, and $LG$ may be evaluated for the
outstanding cases using our previous method (see
\secref{FiveOutstandingKnots}).

Our algorithm is extremely simple to encode, and is readily adaptable
to other state models; the price paid for this generality is that it is
grossly inefficient.  In contrast, the use of representation-specific
properties can improve the efficiency of such computations at the cost
of less generality. For example, \cite{Deguchi:94} and
\cite{SmithiesButler:96} utilise Clebsch--Gordan coefficients and
vector coupling coefficients respectively to facilitate computations.
Our algorithm also compares poorly in performance with
that used by Morton and Short \cite{MortonShort:90} to evaluate the
HOMFLY polynomial -- their (non-state model) method requires much less
storage and time.

The algorithm has again been implemented within the interpreted
environment of \textsc{Mathematica}.  Conversion to a compiled language
would lead to an enormous speedup, but would not enhance its
applicability to links with braid presentations of more than $5$
strings.

\pagebreak


\section{Braid Presentations}
\seclabel{BraidPresentations}

To evaluate $LG$ for a particular link $K$, we require a braid
presentation, i.e. a braid word $\beta$, such that $K$ is equal to the
(vertical) closure $\hat{\beta}$ of $\beta$. For computational
efficiency, it is preferable that $\beta$ be of as few strings as
possible, and as short as possible. Of these two considerations, the
former is critical.

Recall that the \emph{braid index} of a link $K$ is the minimal
number of strings of all possible braid presentations $\beta$
such that $\hat{\beta}=K$.  Unfortunately, we do
not (as yet!) have an algorithm to determine the braid index of an
arbitrary link, much less an algorithm to construct minimal braid
presentations \cite{Birman:74,Murasugi:82}. Indeed, such an algorithm
would answer Artin's \emph{algebraic link problem},%
\footnote{
  Artin stated that the question of the \emph{topological}
  equivalence (or otherwise) of two arbitrary links could be answered
  by the solution of the \emph{algebraic} problem of the equivalence of
  (arbitrary examples of) closed braids associated with the links. Its
  solution depends on the solution of the \emph{Markov equivalence
  problem} \cite{Birman:74} for closed braids. Open!
}
and the braid word so constructed would itself be a complete invariant
for links.  Recall that because no such algorithm is known, it is
hoped that a polynomial invariant may eventually provide a complete
invariant.  Indeed, it is entirely possible that no such algorithm
\emph{exists}, as \cite{MagnusKarrassSolitar:76} describes a
finitely-presented group for which it is known that there can be no
algorithm to determine whether a word is equal to the unit element.

The construction of braid presentations from links is itself of
algorithmic interest. The classical algorithm of Alexander
\cite{Alexander:23} has been used to great effect, but is
labour-intensive, and doesn't yield particularly minimal
presentations.  In 1987, Yamada \cite{Yamada:87} showed that the braid
index of a link is equal to the minimum number of Seifert circles of
all braid presentations of that link. This was proven via the
demonstration of a braid-generating algorithm, improving the method of
Alexander. In 1990, Yamada's algorithm was further improved by Vogel
\cite{Vogel:90}.

In 1987, Jones \cite[p64]{Jones:87} combined the use of the
`Franks--Williams--Morton inequality' together with some explicit braid
presentations to demonstrate the braid indices of all ($249$) prime
knots of up to $10$ crossings%
\footnote{%
  Recall that $10_{161}=10_{162}$. Here, we maintain the artifice that
  these are distinct.
}
(i.e. up to $10_{166}$).
Jones went further \cite[pp66-71]{Jones:87}, \emph{manually}%
\footnote{%
Personal communication. This must have been a
Herculean task!
}
applying Alexander's algorithm to generate a set of braid presentations
for all prime knots to $10_{166}$. These presentations all have string
index equal to the braid index.  We have used this data for our
evaluations of $LG$, although we have no way of independently
confirming that the data does not contain errors originating in the
application of Alexander's algorithm.

Manual continuation of the application of Alexander's (or even Vogel's)
algorithm to yield braid presentations for the prime knots beyond
$10_{166}$ is an unattractive proposition.  Unfortunately, Vogel's
algorithm has not yet been implemented on a computer, hence we do not
yet have at our disposal a substantial list of braid presentations of
prime knots beyond $10_{166}$.

\vfill

\pagebreak


\section{The Knots $10_1,10_3,10_{13},10_{35}$ and $10_{58}$}
\seclabel{FiveOutstandingKnots}

The five prime knots up to $10_{166}$ of braid index $6$ are
$10_1,10_3,10_{13},10_{35}$ and $10_{58}$.  Here, we describe
quasi-Morse diagrams for these knots, from which we are able to
evaluate $LG$ using our previous method, by writing down abstract
tensors associated with them.  The diagrams (sometimes the
reflections!) are taken from \cite[pp257-259]{Kawauchi:96}, and the
conventions are the same as those of \cite{DeWitKauffmanLinks:99b}.

\begin{description}
\item[$\mathbf{10_{1}}$ and $\mathbf{10_{3}}$:]
  These two knots, presented in \figref{Ten1and3}, are structurally
  similar. Their abstract tensors are:
  \be
    {\( T_{10_{1}} \)}^y_x
    & \defeq &
    {\( \sigma^{-1}_{rl} \)}^{a~y}_{b~d}
    \cdot
    {\( \sigma^{8}_{dudu} \)}^{b~d}_{c~x}
    \cdot
    {\( C^{-} \)}^{a}_{c}
    \\
    {\( T_{10_{3}} \)}^y_x
    & \defeq &
    {\( \sigma^{-4}_{rlrl} \)}^{a~y}_{b~d}
    \cdot
    {\( \sigma^{6}_{dudu} \)}^{b~d}_{c~x}
    \cdot
    {\( C^{-} \)}^{a}_{c}.
  \end{eqnarray*}
  \begin{figure}[htbp]
    \begin{center}
      \input{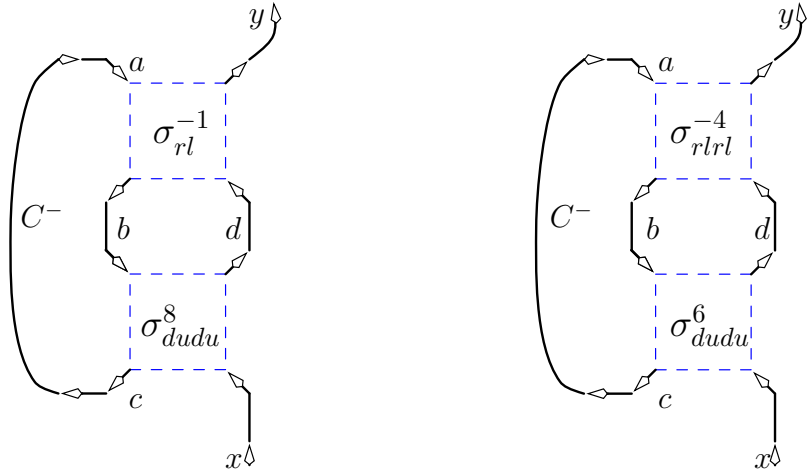}
      \caption{Tangle forms of $10_{1}$ and $10_{3}$.}
      \figlabel{Ten1and3}
    \end{center}
  \end{figure}
\item[$\mathbf{10_{13}}$ and $\mathbf{10_{35}}$:]
  These two knots, presented in \figref{Ten13and35}, are also
  structurally similar. Their abstract tensors are:
  \be
    {\( T_{10_{13}} \)}^y_x
    & \defeq &
    {\( \sigma^{-1}_{rl} \)}^{a~y}_{b~j}
    \cdot
    {\( \sigma^{-1}_{rl} \)}^{b~e}_{c~f}
    \cdot
    {\( \sigma^{4}_{dudu} \)}^{c~f}_{d~g}
    \cdot
    {\( \sigma_{du} \)}^{h~j}_{i~x}
    \cdot
    \\
    & &
    \qquad \qquad
    {\( \Omega^+ \)}_{e~h}
    \cdot
    {\( \mho^- \)}^{g~i}
    \cdot
    {\( C^{-} \)}^{a}_{d}
    \\
    {\( T_{10_{35}} \)}^y_x
    & \defeq &
    {\( \sigma_{rl} \)}^{a~y}_{b~j}
    \cdot
    {\( \sigma^4_{rlrl} \)}^{b~e}_{c~f}
    \cdot
    {\( \sigma^{-1}_{du} \)}^{c~f}_{d~g}
    \cdot
    {\( \sigma^{-1}_{du} \)}^{h~j}_{i~x}
    \cdot
    \\
    & &
    \qquad \qquad
    {\( \Omega^+ \)}_{e~h}
    \cdot
    {\( \mho^- \)}^{g~i}
    \cdot
    {\( C^{-} \)}^{a}_{d}.
  \end{eqnarray*}
  \begin{figure}[htbp]
    \begin{center}
      \input{graphics/Ten13and35.pstex_t}
      \caption{Tangle forms of $10_{13}$ and $10_{35}$.}
      \figlabel{Ten13and35}
    \end{center}
  \end{figure}
\item[$\mathbf{10_{58}}$:]
  This knot, depicted in \figref{Ten58}, has an elegant symmetry. Its
  abstract tensor is:
  \be
    {\( T_{10_{58}} \)}^y_x
    & \defeq &
    {\( \sigma_{du} \)}^{a~y}_{b~m}
    \cdot
    {\( \sigma^{-1}_{rl} \)}^{b~e}_{c~f}
    \cdot
    {\( \Omega^+ \)}_{e~i}
    \cdot
    {\( \mho^- \)}^{f~j}
    \cdot
    {\( \sigma_{du} \)}^{i~m}_{j~n}
    \cdot
    \\
    & &
    \hspace{5pt}
    {\( \sigma^{-1}_{rl} \)}^{c~g}_{d~h}
    \cdot
    {\( \Omega^+ \)}_{g~k}
    \cdot
    {\( \mho^- \)}^{h~l}
    \cdot
    {\( \sigma_{du} \)}^{k~n}_{l~x}
    \cdot
    {\( C^{-} \)}^{a}_{d}.
  \end{eqnarray*}
  \vfill
  \begin{figure}[ht]
    \begin{center}
      \input{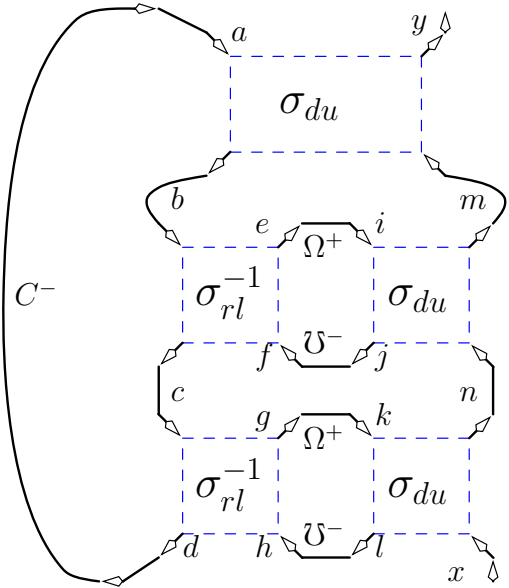}
      \caption{Tangle form of $10_{58}$.}
      \figlabel{Ten58}
    \end{center}
  \end{figure}
\end{description}

\vfill


\section{Behaviour of $LG$}
\seclabel{BehaviourofLG}

Evaluations of $LG$ for the prime knots to $10_{166}$ are presented in
the Appendix.  To reduce the volume of text, we have utilised the
symmetry $LG_K(q,p)=LG_K(q,p^{-1})$ (recall that this indicates that
$LG$ is unable to detect inversion), and only printed the
($q$-polynomial) coefficients of non-negative powers of $P\defeq p^2$.
To illustrate the encoding, consider the entry for $8_{17}$:
\begin{eqnarray*}
  & &
  4 q^{-4} +  68 q^{-2} + 139       + 68 q^{2} + 4 q^{4},
  \\
  & &
  - 22 q^{-3} - 102 q^{-1} - 102 q^{1} - 22 q^{3},
  \\
  & &
  3 q^{-4} +  40 q^{-2} +  82       + 40 q^{2} + 3 q^{4},
  \\
  & &
  -  7 q^{-3} -  36 q^{-1} -  36 q^{1} -  7 q^{3},
  \\
  & &
  7 q^{-2} +  18        +   7 q^{2},
  \\
  & &
  -  4 q^{-1} -   4 q^{1},
  \\
  & &
  1.
\end{eqnarray*}
This is intended to be read:
\begin{eqnarray*}
  LG_{8_{17}}
  \eq
  \hspace{23.6mm}
  (     4 q^{-4} +  68 q^{-2} + 139       + 68 q^{2} + 4 q^{4} )
  \\
  & &
  +
  (P^{-1} + P^{1})
  (  - 22 q^{-3} - 102 q^{-1} - 102 q^{1} - 22 q^{3})
  \\
  & &
  +
  (P^{-2} + P^{2})
  (  3 q^{-4} +  40 q^{-2} +  82       + 40 q^{2} + 3 q^{4})
  \\
  & &
  +
  (P^{-3} + P^{3})
  (  -  7 q^{-3} -  36 q^{-1} -  36 q^{1} -  7 q^{3})
  \\
  & &
  +
  (P^{-4} + P^{4})
  (  7 q^{-2} +  18        +   7 q^{2})
  \\
  & &
  +
  (P^{-5} + P^{5})
  (  -  4 q^{-1} -   4 q^{1})
  \\
  & &
  +
  (P^{-6} + P^{6}).
\end{eqnarray*}
The polynomial is invariant under $P\mapsto P^{-1}$, and in this case,
is also palindromic in $q$ (i.e. invariant under $q\mapsto q^{-1}$), as
$8_{17}$ is amphichiral, i.e. equal to its own reflection.%
\footnote{%
Strictly speaking, we mean that $8_{17}$ is $(+)$amphichiral
\cite[p~4]{Kawauchi:96}.
In fact, it is also $(-)$amphichiral.
}
To illustrate the economy and utility of the notation, replacing $P$
with $p^2$ yields the polynomial in a traditional form. Replacing
$p\mapsto q^{1/2} p$ in this yields the form presented in the original
work \cite{DeWitKauffmanLinks:99b}, viz, the opaque form:
\begin{eqnarray*}
  & &
  \hspace{-17mm}
    139
  - 7   p^{-6}
  + 40  p^{-4}
  - 102 p^{-2}
  - 102 p^{2}
  + 40  p^{4}
  - 7   p^{6}
  +     p^{-12} q^{-6}
  - 4   p^{-10} q^{-6}
  \\
  & &
  \hspace{-17mm}
  + 7   p^{-8}  q^{-6}
  - 7   p^{-6}  q^{-6}
  + 3   p^{-4}  q^{-6}
  + 4           q^{-4}
  - 4   p^{-10} q^{-4}
  + 18  p^{-8}  q^{-4}
  - 36  p^{-6}  q^{-4}
  \\
  & &
  \hspace{-17mm}
  + 40  p^{-4}  q^{-4}
  - 22  p^{-2}  q^{-4}
  + 68          q^{-2}
  + 7   p^{-8}  q^{-2}
  - 36  p^{-6}  q^{-2}
  + 82  p^{-4}  q^{-2}
  - 102 p^{-2}  q^{-2}
  \\
  & &
  \hspace{-17mm}
  - 22  p^{2}   q^{-2}
  + 3   p^{4}   q^{-2}
  + 68          q^{2}
  + 3   p^{-4}  q^{2}
  - 22  p^{-2}  q^{2}
  - 102 p^{2}   q^{2}
  + 82  p^{4}   q^{2}
  - 36  p^{6}   q^{2}
  \\
  & &
  \hspace{-17mm}
  + 7   p^{8}   q^{2}
  + 4           q^{4}
  - 22  p^{2}   q^{4}
  + 40  p^{4}   q^{4}
  - 36  p^{6}   q^{4}
  + 18  p^{8}   q^{4}
  - 4   p^{10}  q^{4}
  + 3   p^{4}   q^{6}
  - 7   p^{6}   q^{6}
  \\
  & &
  \hspace{-17mm}
  + 7   p^{8}   q^{6}
  - 4   p^{10}  q^{6}
  +     p^{12}  q^{6}.
\end{eqnarray*}

\vfill



\subsection{Observations for the Prime Knots to $10_{166}$}

To evaluate $LG$ for the prime knots to $10_{166}$, we have
utilised the braid presentations supplied by Jones in \cite{Jones:87}.
That table of braids also contains evaluations of the
Jones polynomial $V(t)$,%
\footnote{
  Actually the data in \cite{Jones:87} is that for a related polynomial
  $W(t)\defeq(1-V(t))/(1-t)(1-t^3)$, which has integral coefficients.
}
and we confirm that we have made no transcription errors in obtaining
our list of braids by checking that our method yields
Jones polynomials
in agreement with those of \cite{Jones:87}.
These comments of course apply to all but the five knots of braid
index $6$.  Recall that the diagrams for those knots
were used to generate abstract tensors
manually. Again, Jones polynomials that agree with those presented in
\cite{Jones:87} were obtained,
indicating the veracity of our abstract tensors.  Confirmation that our
algorithm and code are correct is provided by the agreement of our data
for $LG$ with that of the subset published in
\cite{DeWitKauffmanLinks:99b}.  Thus, our evaluations of $LG$ are
correct modulo errors in the original construction of the braids.

Generally, $LG$ for any of the prime knots has the following properties.
\begin{itemize}

\item
  The $q$-polynomials are nonzero, the only exceptions being those
  corresponding to $(P^{-5}+P^{5})$ for $LG$ of both
  $10_{154}$ and $10_{161}$.
  (This fact does not appear to be significant.)

\item
  The $q$-polynomials corresponding to even (respectively odd) powers
  of $P$ contain only even (respectively odd) powers of $q$.

\item
  The coefficients (of the powers of $q$) within the $q$-polynomials
  corresponding to odd (respectively even) powers of $P$ are all
  negative (respectively positive).  This is true for all the knots
  except the following ten:  $8_{19}$, $10_{124}$, $10_{128}$,
  $10_{132}$, $10_{139}$, $10_{145}$, $10_{152}$--$10_{154}$, and
  $10_{161}$.  Of these knots, the coefficients of the $q$-polynomials
  of $8_{19}$ and $10_{124}$ are all of the same sign within each
  $q$-polynomial, whilst those of the other eight
  sometimes change sign within a $q$-polynomial, so these latter
  knots are in a sense, even more unusual. Furthermore:

  \begin{itemize}
  \item
    These unusual knots all have braid index either $3$ or $4$.

  \item
    These unusual knots are all chiral.  For the record, the $20$
    amphichiral knots to $10_{166}$ are $4_1$, $6_3$, $8_3$, $8_9$,
    $8_{12}$, $8_{17}$, $8_{18}$, $10_{17}$, $10_{33}$, $10_{37}$,
    $10_{43}$, $10_{45}$, $10_{79}$, $10_{81}$, $10_{88}$, $10_{99}$,
    $10_{109}$, $10_{115}$, $10_{118}$ and $10_{123}$.

  \item
    These unusual knots are all invertible except for $10_{153}$.  For
    the record, the $36$ noninvertible knots to $10_{166}$ are
    $8_{17}$, $9_{32}$, $9_{33}$, $10_{67}$, $10_{79}$--$10_{88}$,
    $10_{90}$--$10_{95}$, $10_{98}$, $10_{102}$, $10_{106}$,
    $10_{107}$, $10_{109}$, $10_{110}$, $10_{115}$,
    $10_{117}$--$10_{119}$, $10_{147}$--$10_{151}$ and $10_{153}$.

  \item
    These unusual knots are all non-alternating. For the record,
    of the knots to $10$ crossings, only $8_{19}$--$8_{21}$,
    $9_{42}$--$9_{49}$ and $10_{124}$--$10_{166}$ are non-alternating.
  \item
    Furthermore, the knots $8_{19}$ and $10_{124}$ are independently of
    interest in that they are non-hyperbolic knots. In fact, of all the
    prime knots of up to $10$ crossings, apart from these two, only the
    $2$-string knots $3_1$, $5_1$, $7_1$ and $9_1$ are non-hyperbolic.

  \end{itemize}

\end{itemize}

\vfill

\pagebreak

Considering the data as a whole, we observe:
\begin{itemize}
\item
  $LG$ assigns \emph{distinct} polynomials to each of the (distinct)
  knots tested. (We also observe that there are no hidden equalities in
  the data by checking if there are equivalences between any pairs
  if the interchange $q\mapsto q^{-1}$ is made in \emph{one} of the
  pair.) This is a significant result, as it is \emph{not} a feature of
  the HOMFLY polynomial $P$; indeed \cite{Jones:87,Kawauchi:96} to
  $10_{166}$ we have four HOMFLY equivalences:  $P_{5_1}=P_{10_{132}}$,
  $P_{8_8}=P_{10_{129}}$, $P_{8_{16}}=P_{10_{156}}$ and
  $P_{10_{25}}=P_{10_{56}}$.  Investigating each of these
  `HOMFLY-pairs', we find that whilst $LG$ of $5_1$ and $10_{132}$ are
  quite distinct (recall that $LG$ of $10_{132}$ had unusual
  coefficients), the other three pairs have Links--Gould invariants
  that are remarkably similar, in that the coefficients of their
  $q$-polynomials are very close (in particular, compare $LG$ of
  $10_{25}$ with that of $10_{56}$).  In this sense, $LG$ \emph{only
  just} distinguishes these pairs.

  (There are also $9$ other coincidences of the $V(t)$ amongst the $10$
  crossing knots \cite{Jones:87}. Within each pair, the
  Alexander--Conway polynomials are distinct, hence they are not HOMFLY
  pairs.)

\item
  $LG$ correctly detects the chirality of all chiral prime knots
  tested. This is also a significant result, as it is not true for the
  HOMFLY polynomial.  Indeed, of the prime knots to
  $10_{166}$, the HOMFLY invariant fails to distinguish the chirality
  of the six knots $9_{42},10_{48},10_{71},10_{91},10_{104}$ and
  $10_{125}$ \cite{Jones:87,Kawauchi:96}.  (Recall that in
  \cite{DeWitKauffmanLinks:99b}, we discovered that $LG$ detects the
  chirality of $9_{42}$ and $10_{48}$.) The $LG$ polynomials of these
  six knots are not palindromic in $q$, but for four of them
  ($10_{48},10_{71},10_{91}$ and $10_{104}$) are curiously `almost
  palindromic', in that the $q$-polynomials are almost palindromic. In
  this sense, $LG$ \emph{only just} succeeds in detecting the chirality
  of these knots.

\end{itemize}


\subsection{$LG$ of a Set of Divers Links}

It is of interest to evaluate $LG$ for multicomponent links, to see if
the resultant polynomials have the same properties as those of the
prime knots. Here, we present evaluations of $LG$ for a small set of
divers links, mostly taken from \cite[pp108-110]{Jones:85} (same
notation).  The list includes several links of $2$ and $3$ components.

$4^2_{1a}$ and $4^2_{1b}$ are different orientations of the same
two-component link $4^2_1$; the latter may be obtained from the former
by reversing the orientation of only \emph{one} of its components --
they \emph{are} distinct oriented links \cite[p110]{Jones:85}.  The
notation ``JE'' means ``Jones' Example''.  ``KT'' is the
Kinoshita--Terasaka Knot \cite{KinoshitaTerasaka:57} and ``C'' is the
Conway Knot. These two are a mutant pair of $11$ crossing prime knots.

We have used our method to evaluate $V(t)$ for these links. That our
$V(t)$ are in agreement with \cite{Jones:85} gives us confidence that
our evaluations of $LG$ are correct modulo errors in the source
braids.

\pagebreak

\vspace{-10mm}
\scriptsize
\begin{eqnarray*}
  & & \hspace{-17mm}
  LG_{2^2_1}
  =
  - 1 - q^{2},
  q^{1}
  \\
  & & \hspace{-17mm}
  LG_{4^2_{1a}}
  =
  - 1 - 2 q^{2} - q^{4},
  q^{1} + 2 q^{3},
  - q^{2} - q^{4},
  q^{3}
  \\
  & & \hspace{-17mm}
  LG_{4^2_{1b}}
  =
  - 3 - 4 q^{2} - q^{4},
  3 q^{1} + q^{3}
  \\
  & & \hspace{-17mm}
  LG_{5^2_1}
  =
  - 8 q^{-2} - 10 - 2 q^{2},
  2 q^{-3} + 9 q^{-1} + 4 q^{1},
  - 3 q^{-2} - 3,
  q^{-1}
  \\
  & & \hspace{-17mm}
  LG_{6^2_1}
  =
  - 1 - 2 q^{2} - 2 q^{4} - q^{6},
  q^{1} + 2 q^{3} + 2 q^{5},
  - q^{2} - 2 q^{4} - q^{6},
  q^{3} + 2 q^{5},
  - q^{4} - q^{6},
  q^{5}
  \\
  & & \hspace{-17mm}
  LG_{6^2_2}
  =
  - 3 - 12 q^{2} - 10 q^{4} - q^{6},
  5 q^{1} + 12 q^{3} + 4 q^{5},
  - 5 q^{2} - 6 q^{4} - q^{6},
  3 q^{3} + q^{5}
  \\
  & & \hspace{-17mm}
  LG_{6^2_3}
  =
  - 4 q^{-2} - 23 - 22 q^{2} - 3 q^{4},
  8 q^{-1} + 20 q^{1} + 7 q^{3},
  - 5 - 5 q^{2},
  q^{1}
  \\
  & & \hspace{-17mm}
  LG_{6^3_1}
  =
  2 q^{-2} + 13 + 18 q^{2} + 3 q^{4},
  - 5 q^{-1} - 15 q^{1} - 10 q^{3},
  5 + 10 q^{2} + 2 q^{4},
  - 3 q^{1} - 3 q^{3},
  q^{2}
  \\
  & & \hspace{-17mm}
  LG_{6^3_2}
  =
  16 q^{-2} + 38 + 16 q^{2},
  - 3 q^{-3} - 25 q^{-1} - 25 q^{1} - 3 q^{3},
  6 q^{-2} + 16 + 6 q^{2},
  - 4 q^{-1} - 4 q^{1},
  1
  \\
  & & \hspace{-17mm}
  LG_{6^3_3}
  =
  1 + 2 q^{2} + q^{4},
  - q^{1} - q^{3},
  q^{4},
  - q^{3} - q^{5},
  q^{4}
  \\
  & & \hspace{-17mm}
  LG_{JE_{1,1}}
  =
  - 3 - 4 q^{2} - q^{4},
  3 q^{1} + q^{3}
  \\
  & & \hspace{-17mm}
  LG_{JE_{1,2}}
  =
  - 1 - 2 q^{2} - q^{4},
  q^{-1} + q^{1} + q^{3},
  - q^{-2} - 1,
  q^{-1}
  \\
  & & \hspace{-17mm}
  LG_{JE_{2,1}}
  =
  LG_{JE_{2,2}}
  =
  \\
  & & \hspace{-2mm}
  - 1 - 11 q^{2} - 32 q^{4} - 24 q^{6} - 2 q^{8},
  3 q^{1} + 19 q^{3} + 30 q^{5} + 9 q^{7},
  - 5 q^{2} - 19 q^{4} - 15 q^{6} - q^{8},
  5 q^{3} + 11 q^{5} + 3 q^{7},
  \\
  & & \hspace{-2mm}
  - 3 q^{4} - 3 q^{6},
  q^{5}
  \\
  & & \hspace{-17mm}
  LG_{KT}
  =
  LG_{C}
  =
  \\
  & & \hspace{-2mm}
  10 q^{-4} + 6 q^{-2} - 23 - 8 q^{2} + 14 q^{4} + 2 q^{6},
  - 2 q^{-5} - 11 q^{-3} + 9 q^{-1} + 18 q^{1} - 7 q^{3} - 7 q^{5},
  \\
  & & \hspace{-2mm}
  3 q^{-4} + 2 q^{-2} - 9 - 3 q^{2} + 6 q^{4} + q^{6},
  - q^{-3} + q^{-1} + 2 q^{1} - q^{3} - q^{5}
\end{eqnarray*}
\normalsize

\begin{itemize}
\item
  As for the prime knots, the $q$-polynomials corresponding to even
  (respectively odd) powers of $P$ contain only even (respectively odd)
  powers of $q$. Also, the coefficients (of the powers of $q$) within
  each $q$-polynomial are generally of consistent sign.

\item
  Only $LG_{6^3_2}$ is palindromic, and this is so as only $6^3_2$ is
  amphichiral.

\item
  $LG_{KT} = LG_{C}$, as shown in \cite{DeWitKauffmanLinks:99b}.  This
  is the first evidence that $LG$ doesn't always distinguish prime
  knots.

\item
  The two links $JE_{1,1}$ and $JE_{1,2}$ have homeomorphic
  complements, but are distinguished by $LG$ (and for that matter, also
  by $V(t)$).

\item
  Jones provides braid presentations
  $\sigma_1^2 \sigma_2^3 \sigma_3^3$ for $JE_{2,1}$ and
  $\sigma_1^3 \sigma_2^2 \sigma_3^3$ for $JE_{2,2}$,
  so these links are actually conjugate (the conjugacy of a similar
  pair is demonstrated in \cite{Birman:74}). Thus, any invariant should
  take on the same value for both of them, and this is true for both
  $V(t)$ and $LG$.

\end{itemize}


\subsection{$LG$ of the $(3,5,7)$ Pretzel Knot}

For good measure, we include the evaluation of the $(3,5,7)$ pretzel
(for comparison with our previously published result
\cite{DeWitKauffmanLinks:99b}):

\vspace{-10mm}
\scriptsize
\begin{eqnarray*}
  & &
  \hspace{-17mm}
  LG_{(3,5,7)}
  =
  57 + 308 q^{2} + 468 q^{4} + 464 q^{6} + 338 q^{8} + 174 q^{10}
  + 56 q^{12} + 8 q^{14},
  \\
  & & \hspace{-2mm}
  - 113 q^{1} - 279 q^{3} - 329 q^{5} - 279 q^{7} - 172 q^{9}
  - 71 q^{11} - 16 q^{13} - q^{15},
  \\
  & & \hspace{-2mm}
  57 q^{2} + 84 q^{4} + 83 q^{6} + 60 q^{8} + 30 q^{10} + 9 q^{12}
  + q^{14}.
\end{eqnarray*}
\normalsize


\section{Summary}

Significantly:

\begin{itemize}
\item
  $LG$ allocates \emph{distinct} polynomials to each of the ($249$)
  prime knots to $10_{166}$.
\item
  Furthermore, it detects the chirality of those that are chiral.
\end{itemize}

In both these features, $LG$ is more sensitive than both the well-known
two-variable HOMFLY and Kauffman invariants.  We do not yet have a
counterexample to show that it does not always detect the chirality of
a chiral knot.  The evidence that it almost fails to detect chirality
in some cases where the HOMFLY polynomial fails
(\secref{BehaviourofLG}) suggests that such counterexamples exist.
Investigation of the prime knots beyond $10_{166}$ might expose such
examples, but this would require braid presentations of these knots
(\secref{BraidPresentations}), and such are not currently to hand.

It is still an open question whether $LG$ always distinguishes between
nonmutant prime knots.  Also note that we do not have an example of a
nontrivial prime knot $K$ such that $LG_K=LG_{0_1}$.

The method used to evaluate $LG$ is generally applicable to any state
sum invariant once appropriate R matrices and caps and cups have been
determined. Stay tuned for future experiments with other such
invariants -- specifically ones based on $(\dot{0}_m|\dot{\alpha}_n)$
representations of $U_q[gl(m|n)]$.


\section*{Acknowledgements}

I am grateful to Jon Links of The University of Queensland, Australia,
who provided valuable advice with respect to changing the variables
used to describe $LG$.

My research at Kyoto University is funded by a Postdoctoral
Fellowship for Foreign Researchers (\# P99703), provided by the Japan
Society for the Promotion of Science.
D\={o}mo arigat\={o} gozaimashita.

This paper was first placed in the LANL archives as
\texttt{math.GT/9906059} on 11 June 1999.  It is also RIMS preprint
RIMS-1235.

\pagebreak


\appendix

\section{$LG$ of the Prime Knots to $10_{166}$}

\vspace{-10mm}
\scriptsize
\begin{eqnarray*}
  & & \hspace{-17mm}
  LG_{3_{1}}
  =
  1 + 2 q^{2},
  - q^{1} - q^{3},
  q^{2}
  \\
  & & \hspace{-17mm}
  LG_{4_{1}}
  =
  2 q^{-2} + 7 + 2 q^{2},
  - 3 q^{-1} - 3 q^{1},
  1
  \\
  & & \hspace{-17mm}
  LG_{5_{1}}
  =
  1 + 2 q^{2} + 2 q^{4},
  - q^{1} - 2 q^{3} - q^{5},
  q^{2} + 2 q^{4},
  - q^{3} - q^{5},
  q^{4}
  \\
  & & \hspace{-17mm}
  LG_{5_{2}}
  =
  3 + 10 q^{2} + 4 q^{4},
  - 5 q^{1} - 6 q^{3} - q^{5},
  3 q^{2} + q^{4}
  \\
  & & \hspace{-17mm}
  LG_{6_{1}}
  =
  4 q^{-2} + 17 + 10 q^{2} + 2 q^{4},
  - 7 q^{-1} - 10 q^{1} - 3 q^{3},
  3 + q^{2}
  \\
  & & \hspace{-17mm}
  LG_{6_{2}}
  =
  2 q^{-2} + 11 + 14 q^{2} + 2 q^{4},
  - 4 q^{-1} - 12 q^{1} - 8 q^{3},
  4 + 9 q^{2} + 2 q^{4},
  - 3 q^{1} - 3 q^{3},
  q^{2}
  \\
  & & \hspace{-17mm}
  LG_{6_{3}}
  =
  10 q^{-2} + 25 + 10 q^{2},
  - 2 q^{-3} - 16 q^{-1} - 16 q^{1} - 2 q^{3},
  4 q^{-2} + 11 + 4 q^{2},
  - 3 q^{-1} - 3 q^{1},
  1
  \\
  & & \hspace{-17mm}
  LG_{7_{1}}
  =
  1 + 2 q^{2} + 2 q^{4} + 2 q^{6},
  - q^{1} - 2 q^{3} - 2 q^{5} - q^{7},
  q^{2} + 2 q^{4} + 2 q^{6},
  - q^{3} - 2 q^{5} - q^{7},
  q^{4} + 2 q^{6},
  - q^{5} - q^{7},
  q^{6}
  \\
  & & \hspace{-17mm}
  LG_{7_{2}}
  =
  5 + 20 q^{2} + 14 q^{4} + 4 q^{6},
  - 9 q^{1} - 14 q^{3} - 6 q^{5} - q^{7},
  5 q^{2} + 3 q^{4} + q^{6}
  \\
  & & \hspace{-17mm}
  LG_{7_{3}}
  =
  3 + 12 q^{2} + 16 q^{4} + 4 q^{6},
  - 5 q^{1} - 14 q^{3} - 10 q^{5} - q^{7},
  5 q^{2} + 12 q^{4} + 4 q^{6},
  - 5 q^{3} - 6 q^{5} - q^{7},
  3 q^{4} + q^{6}
  \\
  & & \hspace{-17mm}
  LG_{7_{4}}
  =
  9 + 38 q^{2} + 28 q^{4} + 6 q^{6},
  - 17 q^{1} - 27 q^{3} - 11 q^{5} - q^{7},
  9 q^{2} + 6 q^{4} + q^{6}
  \\
  & & \hspace{-17mm}
  LG_{7_{5}}
  =
  3 + 20 q^{2} + 32 q^{4} + 10 q^{6},
  - 7 q^{1} - 26 q^{3} - 21 q^{5} - 2 q^{7},
  9 q^{2} + 20 q^{4} + 7 q^{6},
  - 7 q^{3} - 8 q^{5} - q^{7},
  3 q^{4} + q^{6}
  \\
  & & \hspace{-17mm}
  LG_{7_{6}}
  =
  4 q^{-2} + 31 + 50 q^{2} + 16 q^{4},
  - 10 q^{-1} - 37 q^{1} - 30 q^{3} - 3 q^{5},
  10 + 22 q^{2} + 7 q^{4},
  - 5 q^{1} - 5 q^{3},
  q^{2}
  \\
  & & \hspace{-17mm}
  LG_{7_{7}}
  =
  24 q^{-2} + 67 + 38 q^{2} + 4 q^{4},
  - 4 q^{-3} - 39 q^{-1} - 46 q^{1} - 11 q^{3},
  8 q^{-2} + 24 + 11 q^{2},
  - 5 q^{-1} - 5 q^{1},
  1
  \\
  & & \hspace{-17mm}
  LG_{8_{1}}
  =
  6 q^{-2} + 27 + 22 q^{2} + 10 q^{4} + 2 q^{6},
  - 11 q^{-1} - 18 q^{1} - 10 q^{3} - 3 q^{5},
  5 + 3 q^{2} + q^{4}
  \\
  & & \hspace{-17mm}
  LG_{8_{2}}
  =
  2 q^{-2} + 11 + 18 q^{2} + 14 q^{4} + 2 q^{6},
  - 4 q^{-1} - 13 q^{1} - 17 q^{3} - 8 q^{5},
  4 + 13 q^{2} + 14 q^{4} + 2 q^{6},
  - 4 q^{1} - 12 q^{3} - 8 q^{5},
  \\
  & & \hspace{-2mm}
  4 q^{2} + 9 q^{4} + 2 q^{6},
  - 3 q^{3} - 3 q^{5},
  q^{4}
  \\
  & & \hspace{-17mm}
  LG_{8_{3}}
  =
  4 q^{-4} + 26 q^{-2} + 53 + 26 q^{2} + 4 q^{4},
  - 7 q^{-3} - 29 q^{-1} - 29 q^{1} - 7 q^{3},
  3 q^{-2} + 10 + 3 q^{2}
  \\
  & & \hspace{-17mm}
  LG_{8_{4}}
  =
  2 q^{-4} + 14 q^{-2} + 33 + 30 q^{2} + 4 q^{4},
  - 4 q^{-3} - 19 q^{-1} - 31 q^{1} - 16 q^{3},
  4 q^{-2} + 16 + 21 q^{2} + 4 q^{4},
  \\
  & & \hspace{-2mm}
  - 3 q^{-1} - 10 q^{1} - 7 q^{3},
  1 + 3 q^{2}
  \\
  & & \hspace{-17mm}
  LG_{8_{5}}
  =
  2 q^{-2} + 15 + 30 q^{2} + 26 q^{4} + 4 q^{6},
  - 5 q^{-1} - 21 q^{1} - 30 q^{3} - 14 q^{5},
  7 + 22 q^{2} + 22 q^{4} + 3 q^{6},
  - 7 q^{1} - 17 q^{3} - 10 q^{5},
  \\
  & & \hspace{-2mm}
  5 q^{2} + 10 q^{4} + 2 q^{6},
  - 3 q^{3} - 3 q^{5},
  q^{4}
  \\
  & & \hspace{-17mm}
  LG_{8_{6}}
  =
  4 q^{-2} + 33 + 62 q^{2} + 28 q^{4} + 2 q^{6},
  - 10 q^{-1} - 44 q^{1} - 44 q^{3} - 10 q^{5},
  12 + 33 q^{2} + 17 q^{4} + 2 q^{6},
  \\
  & & \hspace{-2mm}
  - 9 q^{1} - 12 q^{3} - 3 q^{5},
  3 q^{2} + q^{4}
  \\
  & & \hspace{-17mm}
  LG_{8_{7}}
  =
  10 q^{-2} + 35 + 40 q^{2} + 10 q^{4},
  - 2 q^{-3} - 18 q^{-1} - 39 q^{1} - 25 q^{3} - 2 q^{5},
  4 q^{-2} + 20 + 31 q^{2} + 10 q^{4},
  \\
  & & \hspace{-2mm}
  - 4 q^{-1} - 18 q^{1} - 16 q^{3} - 2 q^{5},
  4 + 11 q^{2} + 4 q^{4},
  - 3 q^{1} - 3 q^{3},
  q^{2}
  \\
  & & \hspace{-17mm}
  LG_{8_{8}}
  =
  22 q^{-2} + 73 + 54 q^{2} + 12 q^{4},
  - 4 q^{-3} - 42 q^{-1} - 60 q^{1} - 24 q^{3} - 2 q^{5},
  10 q^{-2} + 35 + 23 q^{2} + 4 q^{4},
  \\
  & & \hspace{-2mm}
  - 9 q^{-1} - 12 q^{1} - 3 q^{3},
  3 + q^{2}
  \\
  & & \hspace{-17mm}
  LG_{8_{9}}
  =
  2 q^{-4} + 28 q^{-2} + 59 + 28 q^{2} + 2 q^{4},
  - 10 q^{-3} - 43 q^{-1} - 43 q^{1} - 10 q^{3},
  2 q^{-4} + 18 q^{-2} + 37 + 18 q^{2} + 2 q^{4},
  \\
  & & \hspace{-2mm}
  - 4 q^{-3} - 18 q^{-1} - 18 q^{1} - 4 q^{3},
  4 q^{-2} + 11 + 4 q^{2},
  - 3 q^{-1} - 3 q^{1},
  1
  \\
  & & \hspace{-17mm}
  LG_{8_{10}}
  =
  12 q^{-2} + 51 + 62 q^{2} + 16 q^{4},
  - 2 q^{-3} - 25 q^{-1} - 58 q^{1} - 38 q^{3} - 3 q^{5},
  5 q^{-2} + 29 + 43 q^{2} + 13 q^{4},
  \\
  & & \hspace{-2mm}
  - 6 q^{-1} - 23 q^{1} - 19 q^{3} - 2 q^{5},
  5 + 12 q^{2} + 4 q^{4},
  - 3 q^{1} - 3 q^{3},
  q^{2}
  \\
  & & \hspace{-17mm}
  LG_{8_{11}}
  =
  6 q^{-2} + 51 + 90 q^{2} + 38 q^{4} + 2 q^{6},
  - 16 q^{-1} - 65 q^{1} - 61 q^{3} - 12 q^{5},
  18 + 44 q^{2} + 21 q^{4} + 2 q^{6},
  \\
  & & \hspace{-2mm}
  - 11 q^{1} - 14 q^{3} - 3 q^{5},
  3 q^{2} + q^{4}
  \\
  & & \hspace{-17mm}
  LG_{8_{12}}
  =
  6 q^{-4} + 64 q^{-2} + 129 + 64 q^{2} + 6 q^{4},
  - 17 q^{-3} - 81 q^{-1} - 81 q^{1} - 17 q^{3},
  17 q^{-2} + 41 + 17 q^{2},
  - 7 q^{-1} - 7 q^{1},
  1
  \\
  & & \hspace{-17mm}
  LG_{8_{13}}
  =
  34 q^{-2} + 105 + 74 q^{2} + 14 q^{4},
  - 6 q^{-3} - 61 q^{-1} - 83 q^{1} - 30 q^{3} - 2 q^{5},
  14 q^{-2} + 46 + 29 q^{2} + 4 q^{4},
  \\
  & & \hspace{-2mm}
  - 11 q^{-1} - 14 q^{1} - 3 q^{3},
  3 + q^{2}
  \\
  & & \hspace{-17mm}
  LG_{8_{14}}
  =
  6 q^{-2} + 63 + 124 q^{2} + 60 q^{4} + 4 q^{6},
  - 18 q^{-1} - 85 q^{1} - 86 q^{3} - 19 q^{5},
  22 + 56 q^{2} + 28 q^{4} + 2 q^{6},
  \\
  & & \hspace{-2mm}
  - 13 q^{1} - 16 q^{3} - 3 q^{5},
  3 q^{2} + q^{4}
  \\
  & & \hspace{-17mm}
  LG_{8_{15}}
  =
  7 + 68 q^{2} + 128 q^{4} + 60 q^{6} + 4 q^{8},
  - 21 q^{1} - 93 q^{3} - 91 q^{5} - 19 q^{7},
  28 q^{2} + 67 q^{4} + 33 q^{6} + 2 q^{8},
  \\
  & & \hspace{-2mm}
  - 19 q^{3} - 24 q^{5} - 5 q^{7},
  6 q^{4} + 3 q^{6}
  \\
  & & \hspace{-17mm}
  LG_{8_{16}}
  =
  20 q^{-2} + 89 + 106 q^{2} + 28 q^{4},
  - 3 q^{-3} - 43 q^{-1} - 100 q^{1} - 65 q^{3} - 5 q^{5},
  8 q^{-2} + 50 + 72 q^{2} + 22 q^{4},
  \\
  & & \hspace{-2mm}
  - 10 q^{-1} - 38 q^{1} - 31 q^{3} - 3 q^{5},
  8 + 18 q^{2} + 6 q^{4},
  - 4 q^{1} - 4 q^{3},
  q^{2}
  \\
  & & \hspace{-17mm}
  LG_{8_{17}}
  =
  4 q^{-4} + 68 q^{-2} + 139 + 68 q^{2} + 4 q^{4},
  - 22 q^{-3} - 102 q^{-1} - 102 q^{1} - 22 q^{3},
  3 q^{-4} + 40 q^{-2} + 82 + 40 q^{2} + 3 q^{4},
  \\
  & & \hspace{-2mm}
  - 7 q^{-3} - 36 q^{-1} - 36 q^{1} - 7 q^{3},
  7 q^{-2} + 18 + 7 q^{2},
  - 4 q^{-1} - 4 q^{1},
  1
  \\
  & & \hspace{-17mm}
  LG_{8_{18}}
  =
  6 q^{-4} + 102 q^{-2} + 205 + 102 q^{2} + 6 q^{4},
  - 33 q^{-3} - 152 q^{-1} - 152 q^{1} - 33 q^{3},
  \\
  & & \hspace{-2mm}
  4 q^{-4} + 60 q^{-2} + 122 + 60 q^{2} + 4 q^{4},
  - 10 q^{-3} - 53 q^{-1} - 53 q^{1} - 10 q^{3},
  10 q^{-2} + 25 + 10 q^{2},
  - 5 q^{-1} - 5 q^{1},
  1
  \\
  & & \hspace{-17mm}
  LG_{8_{19}}
  =
  1 + 2 q^{2} + 2 q^{4},
  - q^{1} - q^{3},
  - q^{4} - q^{6},
  q^{3} + q^{5},
  q^{6},
  - q^{5} - q^{7},
  q^{6}
  \\
  & & \hspace{-17mm}
  LG_{8_{20}}
  =
  7 + 8 q^{2} + 4 q^{4},
  - 3 q^{-1} - 7 q^{1} - 5 q^{3} - q^{5},
  q^{-2} + 5 + 3 q^{2} + q^{4},
  - 2 q^{-1} - 2 q^{1},
  1
  \\
  & & \hspace{-17mm}
  LG_{8_{21}}
  =
  9 + 30 q^{2} + 18 q^{4} + 2 q^{6},
  - q^{-1} - 17 q^{1} - 23 q^{3} - 7 q^{5},
  3 + 14 q^{2} + 8 q^{4} + q^{6},
  - 3 q^{1} - 4 q^{3} - q^{5},
  q^{2}
\end{eqnarray*}

\begin{eqnarray*}
  & & \hspace{-17mm}
  LG_{9_{1}}
  =
  1 + 2 q^{2} + 2 q^{4} + 2 q^{6} + 2 q^{8},
  - q^{1} - 2 q^{3} - 2 q^{5} - 2 q^{7} - q^{9},
  q^{2} + 2 q^{4} + 2 q^{6} + 2 q^{8},
  - q^{3} - 2 q^{5} - 2 q^{7} - q^{9},
  \\
  & & \hspace{-2mm}
  q^{4} + 2 q^{6} + 2 q^{8},
  - q^{5} - 2 q^{7} - q^{9},
  q^{6} + 2 q^{8},
  - q^{7} - q^{9},
  q^{8}
  \\
  & & \hspace{-17mm}
  LG_{9_{2}}
  =
  7 + 30 q^{2} + 26 q^{4} + 14 q^{6} + 4 q^{8},
  - 13 q^{1} - 22 q^{3} - 14 q^{5} - 6 q^{7} - q^{9},
  7 q^{2} + 5 q^{4} + 3 q^{6} + q^{8}
  \\
  & & \hspace{-17mm}
  LG_{9_{3}}
  =
  3 + 12 q^{2} + 18 q^{4} + 16 q^{6} + 4 q^{8},
  - 5 q^{1} - 14 q^{3} - 18 q^{5} - 10 q^{7} - q^{9},
  5 q^{2} + 14 q^{4} + 16 q^{6} + 4 q^{8},
  \\
  & & \hspace{-2mm}
  - 5 q^{3} - 14 q^{5} - 10 q^{7} - q^{9},
  5 q^{4} + 12 q^{6} + 4 q^{8},
  - 5 q^{5} - 6 q^{7} - q^{9},
  3 q^{6} + q^{8}
  \\
  & & \hspace{-17mm}
  LG_{9_{4}}
  =
  5 + 24 q^{2} + 40 q^{4} + 20 q^{6} + 4 q^{8},
  - 9 q^{1} - 30 q^{3} - 30 q^{5} - 10 q^{7} - q^{9},
  9 q^{2} + 26 q^{4} + 16 q^{6} + 4 q^{8},
  \\
  & & \hspace{-2mm}
  - 9 q^{3} - 14 q^{5} - 6 q^{7} - q^{9},
  5 q^{4} + 3 q^{6} + q^{8}
  \\
  & & \hspace{-17mm}
  LG_{9_{5}}
  =
  15 + 70 q^{2} + 70 q^{4} + 32 q^{6} + 6 q^{8},
  - 29 q^{1} - 55 q^{3} - 36 q^{5} - 11 q^{7} - q^{9},
  15 q^{2} + 14 q^{4} + 6 q^{6} + q^{8}
  \\
  & & \hspace{-17mm}
  LG_{9_{6}}
  =
  3 + 20 q^{2} + 42 q^{4} + 40 q^{6} + 10 q^{8},
  - 7 q^{1} - 28 q^{3} - 44 q^{5} - 25 q^{7} - 2 q^{9},
  9 q^{2} + 30 q^{4} + 36 q^{6} + 10 q^{8},
  \\
  & & \hspace{-2mm}
  - 9 q^{3} - 28 q^{5} - 21 q^{7} - 2 q^{9},
  9 q^{4} + 20 q^{6} + 7 q^{8},
  - 7 q^{5} - 8 q^{7} - q^{9},
  3 q^{6} + q^{8}
  \\
  & & \hspace{-17mm}
  LG_{9_{7}}
  =
  5 + 42 q^{2} + 88 q^{4} + 52 q^{6} + 10 q^{8},
  - 13 q^{1} - 60 q^{3} - 69 q^{5} - 24 q^{7} - 2 q^{9},
  17 q^{2} + 48 q^{4} + 31 q^{6} + 7 q^{8},
  \\
  & & \hspace{-2mm}
  - 13 q^{3} - 20 q^{5} - 8 q^{7} - q^{9},
  5 q^{4} + 3 q^{6} + q^{8}
  \\
  & & \hspace{-17mm}
  LG_{9_{8}}
  =
  4 q^{-4} + 34 q^{-2} + 91 + 100 q^{2} + 28 q^{4},
  - 10 q^{-3} - 48 q^{-1} - 89 q^{1} - 56 q^{3} - 5 q^{5},
  10 q^{-2} + 37 + 48 q^{2} + 13 q^{4},
  \\
  & & \hspace{-2mm}
  - 5 q^{-1} - 16 q^{1} - 11 q^{3},
  1 + 3 q^{2}
  \\
  & & \hspace{-17mm}
  LG_{9_{9}}
  =
  3 + 22 q^{2} + 58 q^{4} + 62 q^{6} + 16 q^{8},
  - 7 q^{1} - 36 q^{3} - 64 q^{5} - 38 q^{7} - 3 q^{9},
  11 q^{2} + 42 q^{4} + 50 q^{6} + 13 q^{8},
  \\
  & & \hspace{-2mm}
  - 13 q^{3} - 36 q^{5} - 25 q^{7} - 2 q^{9},
  11 q^{4} + 22 q^{6} + 7 q^{8},
  - 7 q^{5} - 8 q^{7} - q^{9},
  3 q^{6} + q^{8}
  \\
  & & \hspace{-17mm}
  LG_{9_{10}}
  =
  9 + 62 q^{2} + 110 q^{4} + 54 q^{6} + 6 q^{8},
  - 21 q^{1} - 82 q^{3} - 82 q^{5} - 22 q^{7} - q^{9},
  25 q^{2} + 66 q^{4} + 39 q^{6} + 6 q^{8},
  \\
  & & \hspace{-2mm}
  - 21 q^{3} - 31 q^{5} - 11 q^{7} - q^{9},
  9 q^{4} + 6 q^{6} + q^{8}
  \\
  & & \hspace{-17mm}
  LG_{9_{11}}
  =
  4 q^{-2} + 33 + 76 q^{2} + 70 q^{4} + 16 q^{6},
  - 10 q^{-1} - 47 q^{1} - 76 q^{3} - 42 q^{5} - 3 q^{7},
  12 + 47 q^{2} + 58 q^{4} + 16 q^{6},
  \\
  & & \hspace{-2mm}
  - 12 q^{1} - 39 q^{3} - 30 q^{5} - 3 q^{7},
  10 q^{2} + 22 q^{4} + 7 q^{6},
  - 5 q^{3} - 5 q^{5},
  q^{4}
  \\
  & & \hspace{-17mm}
  LG_{9_{12}}
  =
  8 q^{-2} + 75 + 152 q^{2} + 88 q^{4} + 16 q^{6},
  - 22 q^{-1} - 99 q^{1} - 110 q^{3} - 36 q^{5} - 3 q^{7},
  24 + 64 q^{2} + 38 q^{4} + 7 q^{6},
  \\
  & & \hspace{-2mm}
  - 13 q^{1} - 18 q^{3} - 5 q^{5},
  3 q^{2} + q^{4}
  \\
  & & \hspace{-17mm}
  LG_{9_{13}}
  =
  9 + 72 q^{2} + 142 q^{4} + 80 q^{6} + 12 q^{8},
  - 23 q^{1} - 100 q^{3} - 110 q^{5} - 35 q^{7} - 2 q^{9},
  29 q^{2} + 80 q^{4} + 51 q^{6} + 9 q^{8},
  \\
  & & \hspace{-2mm}
  - 23 q^{3} - 35 q^{5} - 13 q^{7} - q^{9},
  9 q^{4} + 6 q^{6} + q^{8}
  \\
  & & \hspace{-17mm}
  LG_{9_{14}}
  =
  48 q^{-2} + 161 + 138 q^{2} + 44 q^{4} + 4 q^{6},
  - 8 q^{-3} - 87 q^{-1} - 134 q^{1} - 66 q^{3} - 11 q^{5},
  18 q^{-2} + 64 + 48 q^{2} + 11 q^{4},
  \\
  & & \hspace{-2mm}
  - 13 q^{-1} - 18 q^{1} - 5 q^{3},
  3 + q^{2}
  \\
  & & \hspace{-17mm}
  LG_{9_{15}}
  =
  8 q^{-2} + 89 + 196 q^{2} + 120 q^{4} + 20 q^{6},
  - 24 q^{-1} - 123 q^{1} - 143 q^{3} - 47 q^{5} - 3 q^{7},
  28 + 78 q^{2} + 47 q^{4} + 7 q^{6},
  \\
  & & \hspace{-2mm}
  - 15 q^{1} - 20 q^{3} - 5 q^{5},
  3 q^{2} + q^{4}
  \\
  & & \hspace{-17mm}
  LG_{9_{16}}
  =
  3 + 32 q^{2} + 98 q^{4} + 106 q^{6} + 28 q^{8},
  - 9 q^{1} - 57 q^{3} - 108 q^{5} - 65 q^{7} - 5 q^{9},
  16 q^{2} + 67 q^{4} + 81 q^{6} + 22 q^{8},
  \\
  & & \hspace{-2mm}
  - 19 q^{3} - 55 q^{5} - 39 q^{7} - 3 q^{9},
  16 q^{4} + 31 q^{6} + 10 q^{8},
  - 9 q^{5} - 10 q^{7} - q^{9},
  3 q^{6} + q^{8}
  \\
  & & \hspace{-17mm}
  LG_{9_{17}}
  =
  4 q^{-4} + 44 q^{-2} + 115 + 108 q^{2} + 24 q^{4},
  - 12 q^{-3} - 67 q^{-1} - 115 q^{1} - 64 q^{3} - 4 q^{5},
  16 q^{-2} + 66 + 83 q^{2} + 24 q^{4},
  \\
  & & \hspace{-2mm}
  - 15 q^{-1} - 50 q^{1} - 39 q^{3} - 4 q^{5},
  11 + 24 q^{2} + 8 q^{4},
  - 5 q^{1} - 5 q^{3},
  q^{2}
  \\
  & & \hspace{-17mm}
  LG_{9_{18}}
  =
  9 + 90 q^{2} + 184 q^{4} + 104 q^{6} + 14 q^{8},
  - 27 q^{1} - 128 q^{3} - 141 q^{5} - 42 q^{7} - 2 q^{9},
  37 q^{2} + 98 q^{4} + 60 q^{6} + 9 q^{8},
  \\
  & & \hspace{-2mm}
  - 27 q^{3} - 39 q^{5} - 13 q^{7} - q^{9},
  9 q^{4} + 6 q^{6} + q^{8}
  \\
  & & \hspace{-17mm}
  LG_{9_{19}}
  =
  8 q^{-4} + 96 q^{-2} + 221 + 144 q^{2} + 28 q^{4},
  - 25 q^{-3} - 133 q^{-1} - 161 q^{1} - 57 q^{3} - 4 q^{5},
  29 q^{-2} + 81 + 50 q^{2} + 8 q^{4},
  \\
  & & \hspace{-2mm}
  - 15 q^{-1} - 20 q^{1} - 5 q^{3},
  3 + q^{2}
  \\
  & & \hspace{-17mm}
  LG_{9_{20}}
  =
  4 q^{-2} + 43 + 124 q^{2} + 130 q^{4} + 34 q^{6},
  - 12 q^{-1} - 71 q^{1} - 131 q^{3} - 78 q^{5} - 6 q^{7},
  18 + 75 q^{2} + 91 q^{4} + 25 q^{6},
  \\
  & & \hspace{-2mm}
  - 18 q^{1} - 53 q^{3} - 38 q^{5} - 3 q^{7},
  12 q^{2} + 24 q^{4} + 7 q^{6},
  - 5 q^{3} - 5 q^{5},
  q^{4}
  \\
  & & \hspace{-17mm}
  LG_{9_{21}}
  =
  12 q^{-2} + 119 + 244 q^{2} + 142 q^{4} + 22 q^{6},
  - 34 q^{-1} - 156 q^{1} - 172 q^{3} - 53 q^{5} - 3 q^{7},
  36 + 93 q^{2} + 53 q^{4} + 7 q^{6},
  \\
  & & \hspace{-2mm}
  - 17 q^{1} - 22 q^{3} - 5 q^{5},
  3 q^{2} + q^{4}
  \\
  & & \hspace{-17mm}
  LG_{9_{22}}
  =
  4 q^{-4} + 50 q^{-2} + 143 + 142 q^{2} + 34 q^{4},
  - 13 q^{-3} - 81 q^{-1} - 146 q^{1} - 84 q^{3} - 6 q^{5},
  19 q^{-2} + 81 + 101 q^{2} + 29 q^{4},
  \\
  & & \hspace{-2mm}
  - 18 q^{-1} - 57 q^{1} - 43 q^{3} - 4 q^{5},
  12 + 25 q^{2} + 8 q^{4},
  - 5 q^{1} - 5 q^{3},
  q^{2}
  \\
  & & \hspace{-17mm}
  LG_{9_{23}}
  =
  9 + 102 q^{2} + 226 q^{4} + 140 q^{6} + 22 q^{8},
  - 29 q^{1} - 150 q^{3} - 177 q^{5} - 59 q^{7} - 3 q^{9},
  41 q^{2} + 114 q^{4} + 74 q^{6} + 12 q^{8},
  \\
  & & \hspace{-2mm}
  - 29 q^{3} - 43 q^{5} - 15 q^{7} - q^{9},
  9 q^{4} + 6 q^{6} + q^{8}
  \\
  & & \hspace{-17mm}
  LG_{9_{24}}
  =
  4 q^{-4} + 82 q^{-2} + 193 + 122 q^{2} + 20 q^{4},
  - 26 q^{-3} - 134 q^{-1} - 156 q^{1} - 51 q^{3} - 3 q^{5},
  \\
  & & \hspace{-2mm}
  4 q^{-4} + 53 q^{-2} + 117 + 67 q^{2} + 9 q^{4},
  - 10 q^{-3} - 51 q^{-1} - 53 q^{1} - 12 q^{3},
  10 q^{-2} + 25 + 10 q^{2},
  - 5 q^{-1} - 5 q^{1},
  1
  \\
  & & \hspace{-17mm}
  LG_{9_{25}}
  =
  10 q^{-2} + 119 + 268 q^{2} + 170 q^{4} + 28 q^{6},
  - 32 q^{-1} - 170 q^{1} - 204 q^{3} - 70 q^{5} - 4 q^{7},
  41 + 117 q^{2} + 76 q^{4} + 12 q^{6},
  \\
  & & \hspace{-2mm}
  - 25 q^{1} - 36 q^{3} - 11 q^{5},
  6 q^{2} + 3 q^{4}
  \\
  & & \hspace{-17mm}
  LG_{9_{26}}
  =
  28 q^{-2} + 143 + 206 q^{2} + 82 q^{4} + 4 q^{6},
  - 4 q^{-3} - 62 q^{-1} - 173 q^{1} - 141 q^{3} - 26 q^{5},
  \\
  & & \hspace{-2mm}
  10 q^{-2} + 74 + 128 q^{2} + 57 q^{4} + 4 q^{6},
  - 12 q^{-1} - 57 q^{1} - 56 q^{3} - 11 q^{5},
  10 + 26 q^{2} + 11 q^{4},
  - 5 q^{1} - 5 q^{3},
  q^{2}
  \\
  & & \hspace{-17mm}
  LG_{9_{27}}
  =
  4 q^{-4} + 94 q^{-2} + 235 + 158 q^{2} + 28 q^{4},
  - 28 q^{-3} - 157 q^{-1} - 193 q^{1} - 68 q^{3} - 4 q^{5},
  \\
  & & \hspace{-2mm}
  4 q^{-4} + 58 q^{-2} + 136 + 83 q^{2} + 12 q^{4},
  - 10 q^{-3} - 55 q^{-1} - 60 q^{1} - 15 q^{3},
  10 q^{-2} + 26 + 11 q^{2},
  - 5 q^{-1} - 5 q^{1},
  1
\end{eqnarray*}

\begin{eqnarray*}
  & & \hspace{-17mm}
  LG_{9_{28}}
  =
  30 q^{-2} + 171 + 254 q^{2} + 104 q^{4} + 6 q^{6},
  - 4 q^{-3} - 73 q^{-1} - 210 q^{1} - 172 q^{3} - 31 q^{5},
  \\
  & & \hspace{-2mm}
  12 q^{-2} + 91 + 150 q^{2} + 62 q^{4} + 3 q^{6},
  - 16 q^{-1} - 66 q^{1} - 59 q^{3} - 9 q^{5},
  12 + 27 q^{2} + 10 q^{4},
  - 5 q^{1} - 5 q^{3},
  q^{2}
  \\
  & & \hspace{-17mm}
  LG_{9_{29}}
  =
  6 q^{-4} + 76 q^{-2} + 217 + 214 q^{2} + 52 q^{4},
  - 20 q^{-3} - 123 q^{-1} - 216 q^{1} - 122 q^{3} - 9 q^{5},
  \\
  & & \hspace{-2mm}
  29 q^{-2} + 116 + 136 q^{2} + 37 q^{4},
  - 25 q^{-1} - 71 q^{1} - 50 q^{3} - 4 q^{5},
  14 + 27 q^{2} + 8 q^{4},
  - 5 q^{1} - 5 q^{3},
  q^{2}
  \\
  & & \hspace{-17mm}
  LG_{9_{30}}
  =
  6 q^{-4} + 120 q^{-2} + 287 + 186 q^{2} + 30 q^{4},
  - 35 q^{-3} - 192 q^{-1} - 230 q^{1} - 77 q^{3} - 4 q^{5},
  \\
  & & \hspace{-2mm}
  4 q^{-4} + 67 q^{-2} + 158 + 96 q^{2} + 13 q^{4},
  - 10 q^{-3} - 60 q^{-1} - 67 q^{1} - 17 q^{3},
  10 q^{-2} + 27 + 12 q^{2},
  - 5 q^{-1} - 5 q^{1},
  1
  \\
  & & \hspace{-17mm}
  LG_{9_{31}}
  =
  32 q^{-2} + 197 + 308 q^{2} + 134 q^{4} + 8 q^{6},
  - 4 q^{-3} - 80 q^{-1} - 246 q^{1} - 211 q^{3} - 41 q^{5},
  \\
  & & \hspace{-2mm}
  12 q^{-2} + 100 + 173 q^{2} + 76 q^{4} + 4 q^{6},
  - 16 q^{-1} - 71 q^{1} - 66 q^{3} - 11 q^{5},
  12 + 28 q^{2} + 11 q^{4},
  - 5 q^{1} - 5 q^{3},
  q^{2}
  \\
  & & \hspace{-17mm}
  LG_{9_{32}}
  =
  46 q^{-2} + 241 + 334 q^{2} + 128 q^{4} + 6 q^{6},
  - 6 q^{-3} - 106 q^{-1} - 284 q^{1} - 222 q^{3} - 38 q^{5},
  \\
  & & \hspace{-2mm}
  17 q^{-2} + 124 + 199 q^{2} + 83 q^{4} + 5 q^{6},
  - 21 q^{-1} - 87 q^{1} - 80 q^{3} - 14 q^{5},
  15 + 35 q^{2} + 14 q^{4},
  - 6 q^{1} - 6 q^{3},
  q^{2}
  \\
  & & \hspace{-17mm}
  LG_{9_{33}}
  =
  8 q^{-4} + 162 q^{-2} + 377 + 242 q^{2} + 38 q^{4},
  - 49 q^{-3} - 257 q^{-1} - 302 q^{1} - 99 q^{3} - 5 q^{5},
  \\
  & & \hspace{-2mm}
  6 q^{-4} + 94 q^{-2} + 211 + 125 q^{2} + 16 q^{4},
  - 15 q^{-3} - 82 q^{-1} - 88 q^{1} - 21 q^{3},
  14 q^{-2} + 35 + 15 q^{2},
  \\
  & & \hspace{-2mm}
  - 6 q^{-1} - 6 q^{1},
  1
  \\
  & & \hspace{-17mm}
  LG_{9_{34}}
  =
  12 q^{-4} + 216 q^{-2} + 503 + 330 q^{2} + 54 q^{4},
  - 63 q^{-3} - 334 q^{-1} - 400 q^{1} - 136 q^{3} - 7 q^{5},
  \\
  & & \hspace{-2mm}
  6 q^{-4} + 112 q^{-2} + 262 + 162 q^{2} + 22 q^{4},
  - 15 q^{-3} - 92 q^{-1} - 104 q^{1} - 27 q^{3},
  14 q^{-2} + 37 + 17 q^{2},
  \\
  & & \hspace{-2mm}
  - 6 q^{-1} - 6 q^{1},
  1
  \\
  & & \hspace{-17mm}
  LG_{9_{35}}
  =
  19 + 92 q^{2} + 100 q^{4} + 48 q^{6} + 8 q^{8},
  - 37 q^{1} - 75 q^{3} - 53 q^{5} - 16 q^{7} - q^{9},
  19 q^{2} + 20 q^{4} + 9 q^{6} + q^{8}
  \\
  & & \hspace{-17mm}
  LG_{9_{36}}
  =
  4 q^{-2} + 41 + 100 q^{2} + 94 q^{4} + 22 q^{6},
  - 12 q^{-1} - 62 q^{1} - 101 q^{3} - 55 q^{5} - 4 q^{7},
  17 + 62 q^{2} + 72 q^{4} + 19 q^{6},
  \\
  & & \hspace{-2mm}
  - 16 q^{1} - 46 q^{3} - 33 q^{5} - 3 q^{7},
  11 q^{2} + 23 q^{4} + 7 q^{6},
  - 5 q^{3} - 5 q^{5},
  q^{4}
  \\
  & & \hspace{-17mm}
  LG_{9_{37}}
  =
  10 q^{-4} + 120 q^{-2} + 273 + 176 q^{2} + 32 q^{4},
  - 31 q^{-3} - 164 q^{-1} - 196 q^{1} - 67 q^{3} - 4 q^{5},
  \\
  & & \hspace{-2mm}
  35 q^{-2} + 96 + 58 q^{2} + 8 q^{4},
  - 17 q^{-1} - 22 q^{1} - 5 q^{3},
  3 + q^{2}
  \\
  & & \hspace{-17mm}
  LG_{9_{38}}
  =
  17 + 178 q^{2} + 364 q^{4} + 214 q^{6} + 30 q^{8},
  - 53 q^{1} - 249 q^{3} - 279 q^{5} - 87 q^{7} - 4 q^{9},
  70 q^{2} + 183 q^{4} + 116 q^{6} + 17 q^{8},
  \\
  & & \hspace{-2mm}
  - 47 q^{3} - 69 q^{5} - 23 q^{7} - q^{9},
  14 q^{4} + 10 q^{6} + q^{8}
  \\
  & & \hspace{-17mm}
  LG_{9_{39}}
  =
  16 q^{-2} + 177 + 382 q^{2} + 238 q^{4} + 38 q^{6},
  - 49 q^{-1} - 242 q^{1} - 282 q^{3} - 94 q^{5} - 5 q^{7},
  57 + 154 q^{2} + 97 q^{4} + 14 q^{6},
  \\
  & & \hspace{-2mm}
  - 30 q^{1} - 42 q^{3} - 12 q^{5},
  6 q^{2} + 3 q^{4}
  \\
  & & \hspace{-17mm}
  LG_{9_{40}}
  =
  72 q^{-2} + 399 + 566 q^{2} + 228 q^{4} + 12 q^{6},
  - 9 q^{-3} - 172 q^{-1} - 471 q^{1} - 375 q^{3} - 67 q^{5},
  \\
  & & \hspace{-2mm}
  27 q^{-2} + 200 + 316 q^{2} + 132 q^{4} + 7 q^{6},
  - 33 q^{-1} - 130 q^{1} - 116 q^{3} - 19 q^{5},
  21 + 46 q^{2} + 18 q^{4},
  - 7 q^{1} - 7 q^{3},
  q^{2}
  \\
  & & \hspace{-17mm}
  LG_{9_{41}}
  =
  72 q^{-2} + 261 + 246 q^{2} + 82 q^{4} + 6 q^{6},
  - 12 q^{-3} - 138 q^{-1} - 232 q^{1} - 126 q^{3} - 20 q^{5},
  \\
  & & \hspace{-2mm}
  29 q^{-2} + 111 + 94 q^{2} + 24 q^{4},
  - 23 q^{-1} - 36 q^{1} - 13 q^{3},
  6 + 3 q^{2}
  \\
  & & \hspace{-17mm}
  LG_{9_{42}}
  =
  2 q^{-4} + 6 q^{-2} + 3 ,
  - 3 q^{-3} - 3 q^{-1} - q^{1} - q^{3},
  q^{-2} + 1 + 3 q^{2} + q^{4},
  - 2 q^{1} - 2 q^{3},
  q^{2}
  \\
  & & \hspace{-17mm}
  LG_{9_{43}}
  =
  2 q^{-2} + 11 + 12 q^{2} + 4 q^{4},
  - 4 q^{-1} - 10 q^{1} - 7 q^{3} - q^{5},
  3 + 6 q^{2} + 5 q^{4},
  - 2 q^{1} - 7 q^{3} - 5 q^{5},
  3 q^{2} + 8 q^{4} + 2 q^{6},
  \\
  & & \hspace{-2mm}
  - 3 q^{3} - 3 q^{5},
  q^{4}
  \\
  & & \hspace{-17mm}
  LG_{9_{44}}
  =
  8 q^{-2} + 35 + 30 q^{2} + 10 q^{4},
  - q^{-3} - 17 q^{-1} - 29 q^{1} - 15 q^{3} - 2 q^{5},
  3 q^{-2} + 14 + 10 q^{2} + 3 q^{4},
  \\
  & & \hspace{-2mm}
  - 3 q^{-1} - 4 q^{1} - q^{3},
  1
  \\
  & & \hspace{-17mm}
  LG_{9_{45}}
  =
  19 + 70 q^{2} + 54 q^{4} + 12 q^{6},
  - 2 q^{-1} - 36 q^{1} - 56 q^{3} - 24 q^{5} - 2 q^{7},
  5 + 26 q^{2} + 19 q^{4} + 4 q^{6},
  - 4 q^{1} - 6 q^{3} - 2 q^{5},
  q^{2}
  \\
  & & \hspace{-17mm}
  LG_{9_{46}}
  =
  2 q^{-2} + 11 + 10 q^{2} + 8 q^{4} + 2 q^{6},
  - 4 q^{-1} - 7 q^{1} - 6 q^{3} - 3 q^{5},
  2 + q^{2} + q^{4}
  \\
  & & \hspace{-17mm}
  LG_{9_{47}}
  =
  18 q^{-2} + 57 + 50 q^{2} + 6 q^{4},
  - 3 q^{-3} - 31 q^{-1} - 55 q^{1} - 27 q^{3},
  6 q^{-2} + 30 + 41 q^{2} + 11 q^{4},
  \\
  & & \hspace{-2mm}
  - 6 q^{-1} - 26 q^{1} - 23 q^{3} - 3 q^{5},
  6 + 16 q^{2} + 6 q^{4},
  - 4 q^{1} - 4 q^{3},
  q^{2}
  \\
  & & \hspace{-17mm}
  LG_{9_{48}}
  =
  8 q^{-2} + 63 + 106 q^{2} + 42 q^{4} + 2 q^{6},
  - 20 q^{-1} - 74 q^{1} - 64 q^{3} - 10 q^{5},
  18 + 39 q^{2} + 14 q^{4},
  - 7 q^{1} - 7 q^{3},
  q^{2}
  \\
  & & \hspace{-17mm}
  LG_{9_{49}}
  =
  7 + 44 q^{2} + 66 q^{4} + 22 q^{6},
  - 16 q^{1} - 54 q^{3} - 44 q^{5} - 6 q^{7},
  18 q^{2} + 41 q^{4} + 18 q^{6} + q^{8},
  - 14 q^{3} - 18 q^{5} - 4 q^{7},
  \\
  & & \hspace{-2mm}
  6 q^{4} + 3 q^{6}
\end{eqnarray*}

\pagebreak

\begin{eqnarray*}
  & & \hspace{-17mm}
  LG_{10_{1}}
  =
  8 q^{-2} + 37 + 34 q^{2} + 22 q^{4} + 10 q^{6} + 2 q^{8},
  - 15 q^{-1} - 26 q^{1} - 18 q^{3} - 10 q^{5} - 3 q^{7},
  7 + 5 q^{2} + 3 q^{4} + q^{6}
  \\
  & & \hspace{-17mm}
  LG_{10_{2}}
  =
  2 q^{-2} + 11 + 18 q^{2} + 18 q^{4} + 14 q^{6} + 2 q^{8},
  - 4 q^{-1} - 13 q^{1} - 18 q^{3} - 17 q^{5} - 8 q^{7},
  \\
  & & \hspace{-2mm}
  4 + 13 q^{2} + 18 q^{4} + 14 q^{6} + 2 q^{8},
  - 4 q^{1} - 13 q^{3} - 17 q^{5} - 8 q^{7},
  4 q^{2} + 13 q^{4} + 14 q^{6} + 2 q^{8},
  - 4 q^{3} - 12 q^{5} - 8 q^{7},
  \\
  & & \hspace{-2mm}
  4 q^{4} + 9 q^{6} + 2 q^{8},
  - 3 q^{5} - 3 q^{7},
  q^{6}
  \\
  & & \hspace{-17mm}
  LG_{10_{3}}
  =
  6 q^{-4} + 42 q^{-2} + 95 + 68 q^{2} + 26 q^{4} + 4 q^{6},
  - 11 q^{-3} - 49 q^{-1} - 60 q^{1} - 29 q^{3} - 7 q^{5},
  5 q^{-2} + 18 + 10 q^{2} + 3 q^{4}
  \\
  & & \hspace{-17mm}
  LG_{10_{4}}
  =
  6 q^{-4} + 46 q^{-2} + 59 + 38 q^{2} + 14 q^{4} + 2 q^{6},
  - 24 q^{-3} - 51 q^{-1} - 42 q^{1} - 19 q^{3} - 4 q^{5},
  \\
  & & \hspace{-2mm}
  6 q^{-4} + 33 q^{-2} + 32 + 16 q^{2} + 4 q^{4},
  - 11 q^{-3} - 18 q^{-1} - 10 q^{1} - 3 q^{3},
  5 q^{-2} + 3 + q^{2}
  \\
  & & \hspace{-17mm}
  LG_{10_{5}}
  =
  10 q^{-2} + 35 + 50 q^{2} + 40 q^{4} + 10 q^{6},
  - 2 q^{-3} - 18 q^{-1} - 41 q^{1} - 48 q^{3} - 25 q^{5} - 2 q^{7},
  \\
  & & \hspace{-2mm}
  4 q^{-2} + 20 + 41 q^{2} + 40 q^{4} + 10 q^{6},
  - 4 q^{-1} - 20 q^{1} - 39 q^{3} - 25 q^{5} - 2 q^{7},
  4 + 20 q^{2} + 31 q^{4} + 10 q^{6},
  \\
  & & \hspace{-2mm}
  - 4 q^{1} - 18 q^{3} - 16 q^{5} - 2 q^{7},
  4 q^{2} + 11 q^{4} + 4 q^{6},
  - 3 q^{3} - 3 q^{5},
  q^{4}
  \\
  & & \hspace{-17mm}
  LG_{10_{6}}
  =
  4 q^{-2} + 33 + 78 q^{2} + 82 q^{4} + 28 q^{6} + 2 q^{8},
  - 10 q^{-1} - 47 q^{1} - 83 q^{3} - 56 q^{5} - 10 q^{7},
  \\
  & & \hspace{-2mm}
  12 + 49 q^{2} + 70 q^{4} + 28 q^{6} + 2 q^{8},
  - 12 q^{1} - 46 q^{3} - 44 q^{5} - 10 q^{7},
  12 q^{2} + 33 q^{4} + 17 q^{6} + 2 q^{8},
  \\
  & & \hspace{-2mm}
  - 9 q^{3} - 12 q^{5} - 3 q^{7},
  3 q^{4} + q^{6}
  \\
  & & \hspace{-17mm}
  LG_{10_{7}}
  =
  10 q^{-2} + 95 + 202 q^{2} + 138 q^{4} + 38 q^{6} + 2 q^{8},
  - 28 q^{-1} - 129 q^{1} - 158 q^{3} - 69 q^{5} - 12 q^{7},
  \\
  & & \hspace{-2mm}
  32 + 90 q^{2} + 66 q^{4} + 21 q^{6} + 2 q^{8},
  - 19 q^{1} - 30 q^{3} - 14 q^{5} - 3 q^{7},
  5 q^{2} + 3 q^{4} + q^{6}
  \\
  & & \hspace{-17mm}
  LG_{10_{8}}
  =
  2 q^{-4} + 14 q^{-2} + 37 + 46 q^{2} + 30 q^{4} + 4 q^{6},
  - 4 q^{-3} - 20 q^{-1} - 40 q^{1} - 40 q^{3} - 16 q^{5},
  \\
  & & \hspace{-2mm}
  4 q^{-2} + 20 + 37 q^{2} + 30 q^{4} + 4 q^{6},
  - 4 q^{-1} - 19 q^{1} - 31 q^{3} - 16 q^{5},
  4 + 16 q^{2} + 21 q^{4} + 4 q^{6},
  \\
  & & \hspace{-2mm}
  - 3 q^{1} - 10 q^{3} - 7 q^{5},
  q^{2} + 3 q^{4}
  \\
  & & \hspace{-17mm}
  LG_{10_{9}}
  =
  2 q^{-4} + 28 q^{-2} + 75 + 82 q^{2} + 28 q^{4} + 2 q^{6},
  - 10 q^{-3} - 46 q^{-1} - 82 q^{1} - 56 q^{3} - 10 q^{5},
  \\
  & & \hspace{-2mm}
  2 q^{-4} + 18 q^{-2} + 52 + 69 q^{2} + 28 q^{4} + 2 q^{6},
  - 4 q^{-3} - 20 q^{-1} - 49 q^{1} - 43 q^{3} - 10 q^{5},
  \\
  & & \hspace{-2mm}
  4 q^{-2} + 37 q^{2} + 20 + 18 q^{4} + 2 q^{6},
  - 4 q^{-1} - 18 q^{1} - 18 q^{3} - 4 q^{5},
  4 + 11 q^{2} + 4 q^{4},
  - 3 q^{1} - 3 q^{3},
  q^{2}
  \\
  & & \hspace{-17mm}
  LG_{10_{10}}
  =
  58 q^{-2} + 201 + 194 q^{2} + 82 q^{4} + 14 q^{6},
  - 10 q^{-3} - 109 q^{-1} - 180 q^{1} - 109 q^{3} - 30 q^{5} - 2 q^{7},
  \\
  & & \hspace{-2mm}
  24 q^{-2} + 88 + 78 q^{2} + 29 q^{4} + 4 q^{6},
  - 19 q^{-1} - 30 q^{1} - 14 q^{3} - 3 q^{5},
  5 + 3 q^{2} + q^{4}
  \\
  & & \hspace{-17mm}
  LG_{10_{11}}
  =
  4 q^{-4} + 42 q^{-2} + 143 + 182 q^{2} + 68 q^{4} + 4 q^{6},
  - 10 q^{-3} - 69 q^{-1} - 155 q^{1} - 118 q^{3} - 22 q^{5},
  \\
  & & \hspace{-2mm}
  12 q^{-2} + 64 + 102 q^{2} + 43 q^{4} + 4 q^{6},
  - 9 q^{-1} - 37 q^{1} - 35 q^{3} - 7 q^{5},
  3 + 10 q^{2} + 3 q^{4}
  \\
  & & \hspace{-17mm}
  LG_{10_{12}}
  =
  22 q^{-2} + 107 + 172 q^{2} + 88 q^{4} + 12 q^{6},
  - 4 q^{-3} - 48 q^{-1} - 139 q^{1} - 132 q^{3} - 39 q^{5} - 2 q^{7},
  \\
  & & \hspace{-2mm}
  10 q^{-2} + 62 + 119 q^{2} + 69 q^{4} + 12 q^{6},
  - 12 q^{-1} - 58 q^{1} - 68 q^{3} - 24 q^{5} - 2 q^{7},
  12 + 37 q^{2} + 23 q^{4} + 4 q^{6},
  \\
  & & \hspace{-2mm}
  - 9 q^{1} - 12 q^{3} - 3 q^{5},
  3 q^{2} + q^{4}
  \\
  & & \hspace{-17mm}
  LG_{10_{13}}
  =
  12 q^{-4} + 148 q^{-2} + 361 + 272 q^{2} + 76 q^{4} + 6 q^{6},
  - 37 q^{-3} - 206 q^{-1} - 271 q^{1} - 119 q^{3} - 17 q^{5},
  \\
  & & \hspace{-2mm}
  41 q^{-2} + 120 + 83 q^{2} + 17 q^{4},
  - 19 q^{-1} - 26 q^{1} - 7 q^{3},
  3 + q^{2}
  \\
  & & \hspace{-17mm}
  LG_{10_{14}}
  =
  6 q^{-2} + 67 + 202 q^{2} + 232 q^{4} + 82 q^{6} + 4 q^{8},
  - 18 q^{-1} - 111 q^{1} - 224 q^{3} - 157 q^{5} - 26 q^{7},
  \\
  & & \hspace{-2mm}
  26 + 123 q^{2} + 176 q^{4} + 71 q^{6} + 4 q^{8},
  - 28 q^{1} - 103 q^{3} - 94 q^{5} - 19 q^{7},
  24 q^{2} + 58 q^{4} + 28 q^{6} + 2 q^{8},
  \\
  & & \hspace{-2mm}
  - 13 q^{3} - 16 q^{5} - 3 q^{7},
  3 q^{4} + q^{6}
  \\
  & & \hspace{-17mm}
  LG_{10_{15}}
  =
  12 q^{-4} + 64 q^{-2} + 123 + 102 q^{2} + 22 q^{4},
  - 2 q^{-5} - 26 q^{-3} - 86 q^{-1} - 117 q^{1} - 59 q^{3} - 4 q^{5},
  \\
  & & \hspace{-2mm}
  4 q^{-4} + 32 q^{-2} + 82 + 85 q^{2} + 22 q^{4},
  - 4 q^{-3} - 30 q^{-1} - 64 q^{1} - 42 q^{3} - 4 q^{5},
  4 q^{-2} + 23 + 35 q^{2} + 10 q^{4},
  \\
  & & \hspace{-2mm}
  - 3 q^{-1} - 12 q^{1} - 9 q^{3},
  1 + 3 q^{2}
  \\
  & & \hspace{-17mm}
  LG_{10_{16}}
  =
  6 q^{-4} + 60 q^{-2} + 181 + 216 q^{2} + 78 q^{4} + 4 q^{6},
  - 16 q^{-3} - 92 q^{-1} - 188 q^{1} - 136 q^{3} - 24 q^{5},
  \\
  & & \hspace{-2mm}
  18 q^{-2} + 79 + 116 q^{2} + 47 q^{4} + 4 q^{6},
  - 11 q^{-1} - 41 q^{1} - 37 q^{3} - 7 q^{5},
  3 + 10 q^{2} + 3 q^{4}
  \\
  & & \hspace{-17mm}
  LG_{10_{17}}
  =
  10 q^{-4} + 60 q^{-2} + 109 + 60 q^{2} + 10 q^{4},
  - 2 q^{-5} - 28 q^{-3} - 86 q^{-1} - 86 q^{1} - 28 q^{3} - 2 q^{5},
  \\
  & & \hspace{-2mm}
  10 q^{-4} + 46 q^{-2} + 79 + 46 q^{2} + 10 q^{4},
  - 2 q^{-5} - 18 q^{-3} - 49 q^{-1} - 49 q^{1} - 18 q^{3} - 2 q^{5},
  \\
  & & \hspace{-2mm}
  4 q^{-4} + 20 q^{-2} + 37 + 20 q^{2} + 4 q^{4},
  - 4 q^{-3} - 18 q^{-1} - 18 q^{1} - 4 q^{3},
  4 q^{-2} + 11 + 4 q^{2},
  - 3 q^{-1} - 3 q^{1},
  1
  \\
  & & \hspace{-17mm}
  LG_{10_{18}}
  =
  6 q^{-4} + 72 q^{-2} + 247 + 322 q^{2} + 130 q^{4} + 8 q^{6},
  - 18 q^{-3} - 118 q^{-1} - 265 q^{1} - 204 q^{3} - 39 q^{5},
  \\
  & & \hspace{-2mm}
  22 q^{-2} + 103 + 155 q^{2} + 64 q^{4} + 4 q^{6},
  - 13 q^{-1} - 49 q^{1} - 43 q^{3} - 7 q^{5},
  3 + 10 q^{2} + 3 q^{4}
  \\
  & & \hspace{-17mm}
  LG_{10_{19}}
  =
  14 q^{-4} + 86 q^{-2} + 181 + 154 q^{2} + 34 q^{4},
  - 2 q^{-5} - 32 q^{-3} - 120 q^{-1} - 174 q^{1} - 90 q^{3} - 6 q^{5},
  \\
  & & \hspace{-2mm}
  4 q^{-4} + 40 q^{-2} + 116 + 125 q^{2} + 34 q^{4},
  - 4 q^{-3} - 38 q^{-1} - 89 q^{1} - 61 q^{3} - 6 q^{5},
  4 q^{-2} + 29 + 46 q^{2} + 14 q^{4},
  \\
  & & \hspace{-2mm}
  - 3 q^{-1} - 14 q^{1} - 11 q^{3},
  1 + 3 q^{2}
  \\
  & & \hspace{-17mm}
  LG_{10_{20}}
  =
  6 q^{-2} + 55 + 122 q^{2} + 88 q^{4} + 28 q^{6} + 2 q^{8},
  - 16 q^{-1} - 78 q^{1} - 100 q^{3} - 48 q^{5} - 10 q^{7},
  \\
  & & \hspace{-2mm}
  20 + 61 q^{2} + 47 q^{4} + 17 q^{6} + 2 q^{8},
  - 15 q^{1} - 24 q^{3} - 12 q^{5} - 3 q^{7},
  5 q^{2} + 3 q^{4} + q^{6}
  \\
  & & \hspace{-17mm}
  LG_{10_{21}}
  =
  6 q^{-2} + 53 + 126 q^{2} + 124 q^{4} + 38 q^{6} + 2 q^{8},
  - 16 q^{-1} - 77 q^{1} - 130 q^{3} - 81 q^{5} - 12 q^{7},
  \\
  & & \hspace{-2mm}
  20 + 79 q^{2} + 104 q^{4} + 38 q^{6} + 2 q^{8},
  - 20 q^{1} - 69 q^{3} - 61 q^{5} - 12 q^{7},
  18 q^{2} + 44 q^{4} + 21 q^{6} + 2 q^{8},
  \\
  & & \hspace{-2mm}
  - 11 q^{3} - 14 q^{5} - 3 q^{7},
  3 q^{4} + q^{6}
\end{eqnarray*}

\begin{eqnarray*}
  & & \hspace{-17mm}
  LG_{10_{22}}
  =
  4 q^{-4} + 68 q^{-2} + 185 + 148 q^{2} + 42 q^{4} + 2 q^{6},
  - 22 q^{-3} - 121 q^{-1} - 168 q^{1} - 81 q^{3} - 12 q^{5},
  \\
  & & \hspace{-2mm}
  4 q^{-4} + 48 q^{-2} + 123 + 93 q^{2} + 26 q^{4} + 2 q^{6},
  - 10 q^{-3} - 56 q^{-1} - 72 q^{1} - 30 q^{3} - 4 q^{5},
  \\
  & & \hspace{-2mm}
  12 q^{-2} + 37 + 23 q^{2} + 4 q^{4},
  - 9 q^{-1} - 12 q^{1} - 3 q^{3},
  3 + q^{2}
  \\
  & & \hspace{-17mm}
  LG_{10_{23}}
  =
  38 q^{-2} + 193 + 290 q^{2} + 134 q^{4} + 14 q^{6},
  - 6 q^{-3} - 86 q^{-1} - 241 q^{1} - 212 q^{3} - 53 q^{5} - 2 q^{7},
  \\
  & & \hspace{-2mm}
  16 q^{-2} + 108 + 192 q^{2} + 101 q^{4} + 14 q^{6},
  - 20 q^{-1} - 91 q^{1} - 99 q^{3} - 30 q^{5} - 2 q^{7},
  18 + 50 q^{2} + 29 q^{4} + 4 q^{6},
  \\
  & & \hspace{-2mm}
  - 11 q^{1} - 14 q^{3} - 3 q^{5},
  3 q^{2} + q^{4}
  \\
  & & \hspace{-17mm}
  LG_{10_{24}}
  =
  12 q^{-2} + 141 + 334 q^{2} + 240 q^{4} + 56 q^{6} + 2 q^{8},
  - 38 q^{-1} - 207 q^{1} - 270 q^{3} - 115 q^{5} - 14 q^{7},
  \\
  & & \hspace{-2mm}
  50 + 152 q^{2} + 115 q^{4} + 29 q^{6} + 2 q^{8},
  - 33 q^{1} - 53 q^{3} - 23 q^{5} - 3 q^{7},
  9 q^{2} + 6 q^{4} + q^{6}
  \\
  & & \hspace{-17mm}
  LG_{10_{25}}
  =
  6 q^{-2} + 79 + 274 q^{2} + 332 q^{4} + 120 q^{6} + 6 q^{8},
  - 20 q^{-1} - 145 q^{1} - 311 q^{3} - 221 q^{5} - 35 q^{7},
  \\
  & & \hspace{-2mm}
  34 + 167 q^{2} + 231 q^{4} + 88 q^{6} + 4 q^{8},
  - 38 q^{1} - 127 q^{3} - 108 q^{5} - 19 q^{7},
  28 q^{2} + 62 q^{4} + 28 q^{6} + 2 q^{8},
  \\
  & & \hspace{-2mm}
  - 13 q^{3} - 16 q^{5} - 3 q^{7},
  3 q^{4} + q^{6}
  \\
  & & \hspace{-17mm}
  LG_{10_{26}}
  =
  6 q^{-4} + 122 q^{-2} + 313 + 234 q^{2} + 56 q^{4} + 2 q^{6},
  - 38 q^{-3} - 209 q^{-1} - 274 q^{1} - 117 q^{3} - 14 q^{5},
  \\
  & & \hspace{-2mm}
  6 q^{-4} + 82 q^{-2} + 198 + 139 q^{2} + 32 q^{4} + 2 q^{6},
  - 16 q^{-3} - 87 q^{-1} - 105 q^{1} - 38 q^{3} - 4 q^{5},
  \\
  & & \hspace{-2mm}
  18 q^{-2} + 50 + 29 q^{2} + 4 q^{4},
  - 11 q^{-1} - 14 q^{1} - 3 q^{3},
  3 + q^{2}
  \\
  & & \hspace{-17mm}
  LG_{10_{27}}
  =
  42 q^{-2} + 259 + 442 q^{2} + 236 q^{4} + 30 q^{6},
  - 6 q^{-3} - 106 q^{-1} - 345 q^{1} - 338 q^{3} - 97 q^{5} - 4 q^{7},
  \\
  & & \hspace{-2mm}
  18 q^{-2} + 146 + 280 q^{2} + 158 q^{4} + 22 q^{6},
  - 26 q^{-1} - 125 q^{1} - 138 q^{3} - 41 q^{5} - 2 q^{7},
  24 + 64 q^{2} + 36 q^{4} + 4 q^{6},
  \\
  & & \hspace{-2mm}
  - 13 q^{1} - 16 q^{3} - 3 q^{5},
  3 q^{2} + q^{4}
  \\
  & & \hspace{-17mm}
  LG_{10_{28}}
  =
  70 q^{-2} + 263 + 270 q^{2} + 112 q^{4} + 16 q^{6},
  - 12 q^{-3} - 141 q^{-1} - 249 q^{1} - 156 q^{3} - 38 q^{5} - 2 q^{7},
  \\
  & & \hspace{-2mm}
  32 q^{-2} + 126 + 118 q^{2} + 41 q^{4} + 4 q^{6},
  - 29 q^{-1} - 49 q^{1} - 23 q^{3} - 3 q^{5},
  9 + 6 q^{2} + q^{4}
  \\
  & & \hspace{-17mm}
  LG_{10_{29}}
  =
  6 q^{-4} + 84 q^{-2} + 283 + 338 q^{2} + 122 q^{4} + 6 q^{6},
  - 20 q^{-3} - 145 q^{-1} - 309 q^{1} - 222 q^{3} - 38 q^{5},
  \\
  & & \hspace{-2mm}
  30 q^{-2} + 150 + 218 q^{2} + 89 q^{4} + 6 q^{6},
  - 29 q^{-1} - 105 q^{1} - 93 q^{3} - 17 q^{5},
  19 + 43 q^{2} + 17 q^{4},
  - 7 q^{1} - 7 q^{3},
  q^{2}
  \\
  & & \hspace{-17mm}
  LG_{10_{30}}
  =
  18 q^{-2} + 217 + 516 q^{2} + 384 q^{4} + 96 q^{6} + 4 q^{8},
  - 58 q^{-1} - 312 q^{1} - 412 q^{3} - 181 q^{5} - 23 q^{7},
  \\
  & & \hspace{-2mm}
  72 + 213 q^{2} + 162 q^{4} + 40 q^{6} + 2 q^{8},
  - 41 q^{1} - 65 q^{3} - 27 q^{5} - 3 q^{7},
  9 q^{2} + 6 q^{4} + q^{6}
  \\
  & & \hspace{-17mm}
  LG_{10_{31}}
  =
  30 q^{-4} + 200 q^{-2} + 373 + 222 q^{2} + 40 q^{4},
  - 4 q^{-5} - 74 q^{-3} - 257 q^{-1} - 270 q^{1} - 89 q^{3} - 6 q^{5},
  \\
  & & \hspace{-2mm}
  10 q^{-4} + 83 q^{-2} + 165 + 92 q^{2} + 14 q^{4},
  - 9 q^{-3} - 45 q^{-1} - 47 q^{1} - 11 q^{3},
  3 q^{-2} + 10 + 3 q^{2}
  \\
  & & \hspace{-17mm}
  LG_{10_{32}}
  =
  6 q^{-4} + 142 q^{-2} + 391 + 318 q^{2} + 86 q^{4} + 4 q^{6},
  - 42 q^{-3} - 255 q^{-1} - 357 q^{1} - 166 q^{3} - 22 q^{5},
  \\
  & & \hspace{-2mm}
  6 q^{-4} + 98 q^{-2} + 252 + 187 q^{2} + 44 q^{4} + 2 q^{6},
  - 18 q^{-3} - 109 q^{-1} - 136 q^{1} - 49 q^{3} - 4 q^{5},
  \\
  & & \hspace{-2mm}
  22 q^{-2} + 62 + 36 q^{2} + 4 q^{4},
  - 13 q^{-1} - 16 q^{1} - 3 q^{3},
  3 + q^{2}
  \\
  & & \hspace{-17mm}
  LG_{10_{33}}
  =
  46 q^{-4} + 286 q^{-2} + 505 + 286 q^{2} + 46 q^{4},
  - 6 q^{-5} - 107 q^{-3} - 351 q^{-1} - 351 q^{1} - 107 q^{3} - 6 q^{5},
  \\
  & & \hspace{-2mm}
  14 q^{-4} + 110 q^{-2} + 208 + 110 q^{2} + 14 q^{4},
  - 11 q^{-3} - 53 q^{-1} - 53 q^{1} - 11 q^{3},
  3 q^{-2} + 10 + 3 q^{2}
  \\
  & & \hspace{-17mm}
  LG_{10_{34}}
  =
  34 q^{-2} + 123 + 122 q^{2} + 58 q^{4} + 12 q^{6},
  - 6 q^{-3} - 68 q^{-1} - 114 q^{1} - 74 q^{3} - 24 q^{5} - 2 q^{7},
  \\
  & & \hspace{-2mm}
  16 q^{-2} + 61 + 55 q^{2} + 23 q^{4} + 4 q^{6},
  - 15 q^{-1} - 24 q^{1} - 12 q^{3} - 3 q^{5},
  5 + 3 q^{2} + q^{4}
  \\
  & & \hspace{-17mm}
  LG_{10_{35}}
  =
  10 q^{-4} + 122 q^{-2} + 299 + 230 q^{2} + 70 q^{4} + 6 q^{6},
  - 31 q^{-3} - 171 q^{-1} - 228 q^{1} - 105 q^{3} - 17 q^{5},
  \\
  & & \hspace{-2mm}
  35 q^{-2} + 103 + 73 q^{2} + 17 q^{4},
  - 17 q^{-1} - 24 q^{1} - 7 q^{3},
  3 + q^{2}
  \\
  & & \hspace{-17mm}
  LG_{10_{36}}
  =
  10 q^{-2} + 121 + 292 q^{2} + 224 q^{4} + 66 q^{6} + 4 q^{8},
  - 32 q^{-1} - 175 q^{1} - 235 q^{3} - 111 q^{5} - 19 q^{7},
  \\
  & & \hspace{-2mm}
  40 + 120 q^{2} + 93 q^{4} + 28 q^{6} + 2 q^{8},
  - 23 q^{1} - 36 q^{3} - 16 q^{5} - 3 q^{7},
  5 q^{2} + 3 q^{4} + q^{6}
  \\
  & & \hspace{-17mm}
  LG_{10_{37}}
  =
  28 q^{-4} + 178 q^{-2} + 319 + 178 q^{2} + 28 q^{4},
  - 4 q^{-5} - 68 q^{-3} - 227 q^{-1} - 227 q^{1} - 68 q^{3} - 4 q^{5},
  \\
  & & \hspace{-2mm}
  10 q^{-4} + 77 q^{-2} + 147 + 77 q^{2} + 10 q^{4},
  - 9 q^{-3} - 43 q^{-1} - 43 q^{1} - 9 q^{3},
  3 q^{-2} + 10 + 3 q^{2}
  \\
  & & \hspace{-17mm}
  LG_{10_{38}}
  =
  12 q^{-2} + 153 + 380 q^{2} + 294 q^{4} + 80 q^{6} + 4 q^{8},
  - 40 q^{-1} - 229 q^{1} - 314 q^{3} - 146 q^{5} - 21 q^{7},
  \\
  & & \hspace{-2mm}
  54 + 168 q^{2} + 133 q^{4} + 36 q^{6} + 2 q^{8},
  - 35 q^{1} - 57 q^{3} - 25 q^{5} - 3 q^{7},
  9 q^{2} + 6 q^{4} + q^{6}
  \\
  & & \hspace{-17mm}
  LG_{10_{39}}
  =
  6 q^{-2} + 77 + 242 q^{2} + 276 q^{4} + 94 q^{6} + 4 q^{8},
  - 20 q^{-1} - 133 q^{1} - 267 q^{3} - 182 q^{5} - 28 q^{7},
  \\
  & & \hspace{-2mm}
  32 + 147 q^{2} + 201 q^{4} + 77 q^{6} + 4 q^{8},
  - 34 q^{1} - 115 q^{3} - 100 q^{5} - 19 q^{7},
  26 q^{2} + 60 q^{4} + 28 q^{6} + 2 q^{8},
  \\
  & & \hspace{-2mm}
  - 13 q^{3} - 16 q^{5} - 3 q^{7},
  3 q^{4} + q^{6}
  \\
  & & \hspace{-17mm}
  LG_{10_{40}}
  =
  46 q^{-2} + 301 + 510 q^{2} + 266 q^{4} + 32 q^{6},
  - 6 q^{-3} - 122 q^{-1} - 398 q^{1} - 383 q^{3} - 105 q^{5} - 4 q^{7},
  \\
  & & \hspace{-2mm}
  20 q^{-2} + 168 + 313 q^{2} + 170 q^{4} + 22 q^{6},
  - 30 q^{-1} - 137 q^{1} - 146 q^{3} - 41 q^{5} - 2 q^{7},
  26 + 66 q^{2} + 36 q^{4} + 4 q^{6},
  \\
  & & \hspace{-2mm}
  - 13 q^{1} - 16 q^{3} - 3 q^{5},
  3 q^{2} + q^{4}
  \\
  & & \hspace{-17mm}
  LG_{10_{41}}
  =
  8 q^{-4} + 387 + 116 q^{-2} + 448 q^{2} + 154 q^{4} + 6 q^{6},
  - 28 q^{-3} - 200 q^{-1} - 411 q^{1} - 283 q^{3} - 44 q^{5},
  \\
  & & \hspace{-2mm}
  42 q^{-2} + 197 + 269 q^{2} + 103 q^{4} + 6 q^{6},
  - 37 q^{-1} - 123 q^{1} - 103 q^{3} - 17 q^{5},
  21 + 45 q^{2} + 17 q^{4},
  - 7 q^{1} - 7 q^{3},
  q^{2}
  \\
  & & \hspace{-17mm}
  LG_{10_{42}}
  =
  60 q^{-4} + 394 q^{-2} + 679 + 368 q^{2} + 50 q^{4},
  - 8 q^{-5} - 157 q^{-3} - 507 q^{-1} - 490 q^{1} - 138 q^{3} - 6 q^{5},
  \\
  & & \hspace{-2mm}
  26 q^{-4} + 198 q^{-2} + 352 + 181 q^{2} + 20 q^{4},
  - 33 q^{-3} - 133 q^{-1} - 127 q^{1} - 27 q^{3},
  21 q^{-2} + 47 + 19 q^{2},
  \\
  & & \hspace{-2mm}
  - 7 q^{-1} - 7 q^{1},
  1
\end{eqnarray*}

\begin{eqnarray*}
  & & \hspace{-17mm}
  LG_{10_{43}}
  =
  44 q^{-4} + 296 q^{-2} + 527 + 296 q^{2} + 44 q^{4},
  - 6 q^{-5} - 118 q^{-3} - 393 q^{-1} - 393 q^{1} - 118 q^{3} - 6 q^{5},
  \\
  & & \hspace{-2mm}
  20 q^{-4} + 157 q^{-2} + 291 + 157 q^{2} + 20 q^{4},
  - 27 q^{-3} - 115 q^{-1} - 115 q^{1} - 27 q^{3},
  19 q^{-2} + 45 + 19 q^{2},
  \\
  & & \hspace{-2mm}
  - 7 q^{-1} - 7 q^{1},
  1
  \\
  & & \hspace{-17mm}
  LG_{10_{44}}
  =
  8 q^{-4} + 130 q^{-2} + 477 + 594 q^{2} + 226 q^{4} + 12 q^{6},
  - 30 q^{-3} - 236 q^{-1} - 520 q^{1} - 379 q^{3} - 65 q^{5},
  \\
  & & \hspace{-2mm}
  48 q^{-2} + 237 + 330 q^{2} + 128 q^{4} + 6 q^{6},
  - 43 q^{-1} - 141 q^{1} - 115 q^{3} - 17 q^{5},
  23 + 47 q^{2} + 17 q^{4},
  - 7 q^{1} - 7 q^{3},
  q^{2}
  \\
  & & \hspace{-17mm}
  LG_{10_{45}}
  =
  68 q^{-4} + 478 q^{-2} + 851 + 478 q^{2} + 68 q^{4},
  - 8 q^{-5} - 181 q^{-3} - 616 q^{-1} - 616 q^{1} - 181 q^{3} - 8 q^{5},
  \\
  & & \hspace{-2mm}
  26 q^{-4} + 224 q^{-2} + 417 + 224 q^{2} + 26 q^{4},
  - 33 q^{-3} - 145 q^{-1} - 145 q^{1} - 33 q^{3},
  21 q^{-2} + 49 + 21 q^{2},
  \\
  & & \hspace{-2mm}
  - 7 q^{-1} - 7 q^{1},
  1
  \\
  & & \hspace{-17mm}
  LG_{10_{46}}
  =
  2 q^{-2} + 15 + 34 q^{2} + 42 q^{4} + 30 q^{6} + 4 q^{8},
  - 5 q^{-1} - 22 q^{1} - 39 q^{3} - 38 q^{5} - 16 q^{7},
  \\
  & & \hspace{-2mm}
  7 + 26 q^{2} + 38 q^{4} + 28 q^{6} + 4 q^{8},
  - 8 q^{1} - 26 q^{3} - 32 q^{5} - 14 q^{7},
  8 q^{2} + 23 q^{4} + 22 q^{6} + 3 q^{8},
  \\
  & & \hspace{-2mm}
  - 7 q^{3} - 17 q^{5} - 10 q^{7},
  5 q^{4} + 10 q^{6} + 2 q^{8},
  - 3 q^{5} - 3 q^{7},
  q^{6}
  \\
  & & \hspace{-17mm}
  LG_{10_{47}}
  =
  12 q^{-2} + 53 + 88 q^{2} + 70 q^{4} + 16 q^{6},
  - 2 q^{-3} - 25 q^{-1} - 68 q^{1} - 84 q^{3} - 42 q^{5} - 3 q^{7},
  \\
  & & \hspace{-2mm}
  5 q^{-2} + 31 + 69 q^{2} + 66 q^{4} + 16 q^{6},
  - 6 q^{-1} - 32 q^{1} - 61 q^{3} - 38 q^{5} - 3 q^{7},
  6 + 30 q^{2} + 43 q^{4} + 13 q^{6},
  \\
  & & \hspace{-2mm}
  - 6 q^{1} - 23 q^{3} - 19 q^{5} - 2 q^{7},
  5 q^{2} + 12 q^{4} + 4 q^{6},
  - 3 q^{3} - 3 q^{5},
  q^{4}
  \\
  & & \hspace{-17mm}
  LG_{10_{48}}
  =
  16 q^{-4} + 94 q^{-2} + 165 + 88 q^{2} + 12 q^{4},
  - 3 q^{-5} - 43 q^{-3} - 132 q^{-1} - 129 q^{1} - 39 q^{3} - 2 q^{5},
  \\
  & & \hspace{-2mm}
  13 q^{-4} + 67 q^{-2} + 119 + 68 q^{2} + 12 q^{4},
  - 2 q^{-5} - 22 q^{-3} - 69 q^{-1} - 72 q^{1} - 25 q^{3} - 2 q^{5},
  \\
  & & \hspace{-2mm}
  4 q^{-4} + 24 q^{-2} + 50 + 29 q^{2} + 5 q^{4},
  - 4 q^{-3} - 21 q^{-1} - 23 q^{1} - 6 q^{3},
  4 q^{-2} + 12 + 5 q^{2},
  - 3 q^{-1} - 3 q^{1},
  1
  \\
  & & \hspace{-17mm}
  LG_{10_{49}}
  =
  7 + 72 q^{2} + 206 q^{4} + 232 q^{6} + 82 q^{8} + 4 q^{10},
  - 21 q^{1} - 119 q^{3} - 229 q^{5} - 157 q^{7} - 26 q^{9},
  \\
  & & \hspace{-2mm}
  32 q^{2} + 135 q^{4} + 182 q^{6} + 71 q^{8} + 4 q^{10},
  - 35 q^{3} - 116 q^{5} - 100 q^{7} - 19 q^{9},
  31 q^{4} + 70 q^{6} + 33 q^{8} + 2 q^{10},
  \\
  & & \hspace{-2mm}
  - 19 q^{5} - 24 q^{7} - 5 q^{9},
  6 q^{6} + 3 q^{8}
  \\
  & & \hspace{-17mm}
  LG_{10_{50}}
  =
  6 q^{-2} + 67 + 182 q^{2} + 194 q^{4} + 64 q^{6} + 4 q^{8},
  - 19 q^{-1} - 108 q^{1} - 195 q^{3} - 126 q^{5} - 20 q^{7},
  \\
  & & \hspace{-2mm}
  29 + 115 q^{2} + 147 q^{4} + 53 q^{6} + 3 q^{8},
  - 29 q^{1} - 89 q^{3} - 74 q^{5} - 14 q^{7},
  21 q^{2} + 48 q^{4} + 22 q^{6} + 2 q^{8},
  \\
  & & \hspace{-2mm}
  - 11 q^{3} - 14 q^{5} - 3 q^{7},
  3 q^{4} + q^{6}
  \\
  & & \hspace{-17mm}
  LG_{10_{51}}
  =
  44 q^{-2} + 255 + 402 q^{2} + 194 q^{4} + 22 q^{6},
  - 6 q^{-3} - 109 q^{-1} - 323 q^{1} - 292 q^{3} - 75 q^{5} - 3 q^{7},
  \\
  & & \hspace{-2mm}
  19 q^{-2} + 140 + 247 q^{2} + 128 q^{4} + 17 q^{6},
  - 26 q^{-1} - 110 q^{1} - 115 q^{3} - 33 q^{5} - 2 q^{7},
  21 + 54 q^{2} + 30 q^{4} + 4 q^{6},
  \\
  & & \hspace{-2mm}
  - 11 q^{1} - 14 q^{3} - 3 q^{5},
  3 q^{2} + q^{4}
  \\
  & & \hspace{-17mm}
  LG_{10_{52}}
  =
  16 q^{-4} + 114 q^{-2} + 259 + 228 q^{2} + 52 q^{4},
  - 2 q^{-5} - 41 q^{-3} - 168 q^{-1} - 250 q^{1} - 130 q^{3} - 9 q^{5},
  \\
  & & \hspace{-2mm}
  5 q^{-4} + 56 q^{-2} + 161 + 166 q^{2} + 43 q^{4},
  - 6 q^{-3} - 51 q^{-1} - 109 q^{1} - 70 q^{3} - 6 q^{5},
  5 q^{-2} + 33 + 49 q^{2} + 14 q^{4},
  \\
  & & \hspace{-2mm}
  - 3 q^{-1} - 14 q^{1} - 11 q^{3},
  1 + 3 q^{2}
  \\
  & & \hspace{-17mm}
  LG_{10_{53}}
  =
  21 + 246 q^{2} + 568 q^{4} + 410 q^{6} + 96 q^{8} + 4 q^{10},
  - 69 q^{1} - 361 q^{3} - 466 q^{5} - 197 q^{7} - 23 q^{9},
  \\
  & & \hspace{-2mm}
  94 q^{2} + 272 q^{4} + 207 q^{6} + 49 q^{8} + 2 q^{10},
  - 63 q^{3} - 103 q^{5} - 45 q^{7} - 5 q^{9},
  18 q^{4} + 15 q^{6} + 3 q^{8}
  \\
  & & \hspace{-17mm}
  LG_{10_{54}}
  =
  14 q^{-4} + 80 q^{-2} + 153 + 126 q^{2} + 28 q^{4},
  - 2 q^{-5} - 33 q^{-3} - 108 q^{-1} - 144 q^{1} - 72 q^{3} - 5 q^{5},
  \\
  & & \hspace{-2mm}
  5 q^{-4} + 42 q^{-2} + 101 + 99 q^{2} + 25 q^{4},
  - 6 q^{-3} - 37 q^{-1} - 72 q^{1} - 45 q^{3} - 4 q^{5},
  5 q^{-2} + 25 + 36 q^{2} + 10 q^{4},
  \\
  & & \hspace{-2mm}
  - 3 q^{-1} - 12 q^{1} - 9 q^{3},
  1 + 3 q^{2}
  \\
  & & \hspace{-17mm}
  LG_{10_{55}}
  =
  13 + 158 q^{2} + 388 q^{4} + 298 q^{6} + 80 q^{8} + 4 q^{10},
  - 43 q^{1} - 239 q^{3} - 326 q^{5} - 151 q^{7} - 21 q^{9},
  \\
  & & \hspace{-2mm}
  60 q^{2} + 184 q^{4} + 148 q^{6} + 41 q^{8} + 2 q^{10},
  - 41 q^{3} - 70 q^{5} - 34 q^{7} - 5 q^{9},
  12 q^{4} + 10 q^{6} + 3 q^{8}
  \\
  & & \hspace{-17mm}
  LG_{10_{56}}
  =
  6 q^{-2} + 83 + 276 q^{2} + 328 q^{4} + 118 q^{6} + 6 q^{8},
  - 21 q^{-1} - 148 q^{1} - 309 q^{3} - 218 q^{5} - 36 q^{7},
  \\
  & & \hspace{-2mm}
  35 + 165 q^{2} + 229 q^{4} + 90 q^{6} + 5 q^{8},
  - 37 q^{1} - 125 q^{3} - 109 q^{5} - 21 q^{7},
  27 q^{2} + 62 q^{4} + 29 q^{6} + 2 q^{8},
  \\
  & & \hspace{-2mm}
  - 13 q^{3} - 16 q^{5} - 3 q^{7},
  3 q^{4} + q^{6}
  \\
  & & \hspace{-17mm}
  LG_{10_{57}}
  =
  48 q^{-2} + 331 + 580 q^{2} + 314 q^{4} + 40 q^{6},
  - 6 q^{-3} - 131 q^{-1} - 443 q^{1} - 438 q^{3} - 125 q^{5} - 5 q^{7},
  \\
  & & \hspace{-2mm}
  21 q^{-2} + 182 + 345 q^{2} + 191 q^{4} + 25 q^{6},
  - 32 q^{-1} - 146 q^{1} - 156 q^{3} - 44 q^{5} - 2 q^{7},
  27 + 68 q^{2} + 37 q^{4} + 4 q^{6},
  \\
  & & \hspace{-2mm}
  - 13 q^{1} - 16 q^{3} - 3 q^{5},
  3 q^{2} + q^{4}
  \\
  & & \hspace{-17mm}
  LG_{10_{58}}
  =
  14 q^{-4} + 192 q^{-2} + 505 + 416 q^{2} + 124 q^{4} + 8 q^{6},
  - 47 q^{-3} - 283 q^{-1} - 405 q^{1} - 197 q^{3} - 28 q^{5},
  \\
  & & \hspace{-2mm}
  58 q^{-2} + 183 + 143 q^{2} + 34 q^{4},
  - 31 q^{-1} - 48 q^{1} - 17 q^{3},
  6 + 3 q^{2}
  \\
  & & \hspace{-17mm}
  LG_{10_{59}}
  =
  8 q^{-4} + 128 q^{-2} + 439 + 514 q^{2} + 180 q^{4} + 8 q^{6},
  - 30 q^{-3} - 225 q^{-1} - 466 q^{1} - 322 q^{3} - 51 q^{5},
  \\
  & & \hspace{-2mm}
  47 q^{-2} + 220 + 297 q^{2} + 112 q^{4} + 6 q^{6},
  - 41 q^{-1} - 132 q^{1} - 108 q^{3} - 17 q^{5},
  22 + 46 q^{2} + 17 q^{4},
  - 7 q^{1} - 7 q^{3},
  q^{2}
  \\
  & & \hspace{-17mm}
  LG_{10_{60}}
  =
  12 q^{-4} + 274 q^{-2} + 725 + 574 q^{2} + 148 q^{4} + 8 q^{6},
  - 76 q^{-3} - 453 q^{-1} - 619 q^{1} - 274 q^{3} - 32 q^{5},
  \\
  & & \hspace{-2mm}
  8 q^{-4} + 147 q^{-2} + 373 + 266 q^{2} + 52 q^{4},
  - 20 q^{-3} - 124 q^{-1} - 149 q^{1} - 45 q^{3},
  18 q^{-2} + 48 + 23 q^{2},
  \\
  & & \hspace{-2mm}
  - 7 q^{-1} - 7 q^{1},
  1
  \\
  & & \hspace{-17mm}
  LG_{10_{61}}
  =
  2 q^{-4} + 18 q^{-2} + 49 + 62 q^{2} + 42 q^{4} + 6 q^{6},
  - 5 q^{-3} - 28 q^{-1} - 55 q^{1} - 54 q^{3} - 22 q^{5},
  \\
  & & \hspace{-2mm}
  7 q^{-2} + 30 + 50 q^{2} + 38 q^{4} + 5 q^{6},
  - 7 q^{-1} - 26 q^{1} - 37 q^{3} - 18 q^{5},
  5 + 18 q^{2} + 22 q^{4} + 4 q^{6},
  \\
  & & \hspace{-2mm}
  - 3 q^{1} - 10 q^{3} - 7 q^{5},
  q^{2} + 3 q^{4}
\end{eqnarray*}

\begin{eqnarray*}
  & & \hspace{-17mm}
  LG_{10_{62}}
  =
  12 q^{-2} + 61 + 112 q^{2} + 94 q^{4} + 22 q^{6},
  - 2 q^{-3} - 27 q^{-1} - 84 q^{1} - 110 q^{3} - 55 q^{5} - 4 q^{7},
  \\
  & & \hspace{-2mm}
  5 q^{-2} + 38 + 88 q^{2} + 82 q^{4} + 19 q^{6},
  - 7 q^{-1} - 41 q^{1} - 73 q^{3} - 42 q^{5} - 3 q^{7},
  8 + 35 q^{2} + 46 q^{4} + 13 q^{6},
  \\
  & & \hspace{-2mm}
  - 7 q^{1} - 24 q^{3} - 19 q^{5} - 2 q^{7},
  5 q^{2} + 12 q^{4} + 4 q^{6},
  - 3 q^{3} - 3 q^{5},
  q^{4}
  \\
  & & \hspace{-17mm}
  LG_{10_{63}}
  =
  13 + 144 q^{2} + 332 q^{4} + 244 q^{6} + 66 q^{8} + 4 q^{10},
  - 41 q^{1} - 212 q^{3} - 276 q^{5} - 124 q^{7} - 19 q^{9},
  \\
  & & \hspace{-2mm}
  56 q^{2} + 164 q^{4} + 128 q^{6} + 36 q^{8} + 2 q^{10},
  - 39 q^{3} - 65 q^{5} - 31 q^{7} - 5 q^{9},
  12 q^{4} + 10 q^{6} + 3 q^{8}
  \\
  & & \hspace{-17mm}
  LG_{10_{64}}
  =
  2 q^{-4} + 42 q^{-2} + 141 + 166 q^{2} + 58 q^{4} + 4 q^{6},
  - 12 q^{-3} - 80 q^{-1} - 160 q^{1} - 111 q^{3} - 19 q^{5},
  \\
  & & \hspace{-2mm}
  2 q^{-4} + 27 q^{-2} + 97 + 127 q^{2} + 48 q^{4} + 3 q^{6},
  - 5 q^{-3} - 36 q^{-1} - 85 q^{1} - 67 q^{3} - 13 q^{5},
  \\
  & & \hspace{-2mm}
  7 q^{-2} + 34 + 53 q^{2} + 22 q^{4} + 2 q^{6},
  - 7 q^{-1} - 24 q^{1} - 21 q^{3} - 4 q^{5},
  5 + 12 q^{2} + 4 q^{4},
  - 3 q^{1} - 3 q^{3},
  q^{2}
  \\
  & & \hspace{-17mm}
  LG_{10_{65}}
  =
  40 q^{-2} + 217 + 342 q^{2} + 168 q^{4} + 20 q^{6},
  - 6 q^{-3} - 94 q^{-1} - 276 q^{1} - 253 q^{3} - 68 q^{5} - 3 q^{7},
  \\
  & & \hspace{-2mm}
  17 q^{-2} + 120 + 218 q^{2} + 118 q^{4} + 17 q^{6},
  - 22 q^{-1} - 99 q^{1} - 108 q^{3} - 33 q^{5} - 2 q^{7},
  19 + 52 q^{2} + 30 q^{4} + 4 q^{6},
  \\
  & & \hspace{-2mm}
  - 11 q^{1} - 14 q^{3} - 3 q^{5},
  3 q^{2} + q^{4}
  \\
  & & \hspace{-17mm}
  LG_{10_{66}}
  =
  7 + 98 q^{2} + 348 q^{4} + 430 q^{6} + 162 q^{8} + 8 q^{10},
  - 25 q^{1} - 185 q^{3} - 401 q^{5} - 290 q^{7} - 49 q^{9},
  \\
  & & \hspace{-2mm}
  46 q^{2} + 220 q^{4} + 303 q^{6} + 119 q^{8} + 6 q^{10},
  - 54 q^{3} - 175 q^{5} - 147 q^{7} - 26 q^{9},
  43 q^{4} + 91 q^{6} + 41 q^{8} + 2 q^{10},
  \\
  & & \hspace{-2mm}
  - 22 q^{5} - 27 q^{7} - 5 q^{9},
  6 q^{6} + 3 q^{8}
  \\
  & & \hspace{-17mm}
  LG_{10_{67}}
  =
  14 q^{-2} + 177 + 442 q^{2} + 344 q^{4} + 92 q^{6} + 4 q^{8},
  - 46 q^{-1} - 263 q^{1} - 363 q^{3} - 169 q^{5} - 23 q^{7},
  \\
  & & \hspace{-2mm}
  60 + 188 q^{2} + 150 q^{4} + 40 q^{6} + 2 q^{8},
  - 37 q^{1} - 61 q^{3} - 27 q^{5} - 3 q^{7},
  9 q^{2} + 6 q^{4} + q^{6}
  \\
  & & \hspace{-17mm}
  LG_{10_{68}}
  =
  82 q^{-2} + 309 + 322 q^{2} + 134 q^{4} + 18 q^{6},
  - 14 q^{-3} - 163 q^{-1} - 292 q^{1} - 185 q^{3} - 44 q^{5} - 2 q^{7},
  \\
  & & \hspace{-2mm}
  36 q^{-2} + 142 + 135 q^{2} + 47 q^{4} + 4 q^{6},
  - 31 q^{-1} - 53 q^{1} - 25 q^{3} - 3 q^{5},
  9 + 6 q^{2} + q^{4}
  \\
  & & \hspace{-17mm}
  LG_{10_{69}}
  =
  72 q^{-2} + 471 + 798 q^{2} + 426 q^{4} + 56 q^{6},
  - 8 q^{-3} - 182 q^{-1} - 594 q^{1} - 574 q^{3} - 161 q^{5} - 7 q^{7},
  \\
  & & \hspace{-2mm}
  24 q^{-2} + 218 + 406 q^{2} + 216 q^{4} + 25 q^{6},
  - 30 q^{-1} - 142 q^{1} - 146 q^{3} - 34 q^{5},
  20 + 49 q^{2} + 22 q^{4},
  - 7 q^{1} - 7 q^{3},
  q^{2}
  \\
  & & \hspace{-17mm}
  LG_{10_{70}}
  =
  6 q^{-4} + 92 q^{-2} + 329 + 398 q^{2} + 142 q^{4} + 6 q^{6},
  - 21 q^{-3} - 165 q^{-1} - 360 q^{1} - 258 q^{3} - 42 q^{5},
  \\
  & & \hspace{-2mm}
  33 q^{-2} + 171 + 246 q^{2} + 98 q^{4} + 6 q^{6},
  - 32 q^{-1} - 114 q^{1} - 99 q^{3} - 17 q^{5},
  20 + 44 q^{2} + 17 q^{4},
  - 7 q^{1} - 7 q^{3},
  q^{2}
  \\
  & & \hspace{-17mm}
  LG_{10_{71}}
  =
  48 q^{-4} + 340 q^{-2} + 605 + 334 q^{2} + 46 q^{4},
  - 6 q^{-5} - 131 q^{-3} - 448 q^{-1} - 446 q^{1} - 129 q^{3} - 6 q^{5},
  \\
  & & \hspace{-2mm}
  20 q^{-4} + 172 q^{-2} + 323 + 174 q^{2} + 21 q^{4},
  - 27 q^{-3} - 122 q^{-1} - 124 q^{1} - 29 q^{3},
  19 q^{-2} + 46 + 20 q^{2},
  \\
  & & \hspace{-2mm}
  - 7 q^{-1} - 7 q^{1},
  1
  \\
  & & \hspace{-17mm}
  LG_{10_{72}}
  =
  6 q^{-2} + 93 + 344 q^{2} + 430 q^{4} + 162 q^{6} + 8 q^{8},
  - 22 q^{-1} - 176 q^{1} - 395 q^{3} - 290 q^{5} - 49 q^{7},
  \\
  & & \hspace{-2mm}
  39 + 204 q^{2} + 294 q^{4} + 119 q^{6} + 6 q^{8},
  - 44 q^{1} - 156 q^{3} - 138 q^{5} - 26 q^{7},
  33 q^{2} + 75 q^{4} + 35 q^{6} + 2 q^{8},
  \\
  & & \hspace{-2mm}
  - 15 q^{3} - 18 q^{5} - 3 q^{7},
  3 q^{4} + q^{6}
  \\
  & & \hspace{-17mm}
  LG_{10_{73}}
  =
  68 q^{-2} + 429 + 712 q^{2} + 372 q^{4} + 48 q^{6},
  - 8 q^{-3} - 170 q^{-1} - 539 q^{1} - 511 q^{3} - 140 q^{5} - 6 q^{7},
  \\
  & & \hspace{-2mm}
  24 q^{-2} + 205 + 373 q^{2} + 194 q^{4} + 22 q^{6},
  - 30 q^{-1} - 136 q^{1} - 137 q^{3} - 31 q^{5},
  20 + 48 q^{2} + 21 q^{4},
  - 7 q^{1} - 7 q^{3},
  q^{2}
  \\
  & & \hspace{-17mm}
  LG_{10_{74}}
  =
  18 q^{-2} + 203 + 460 q^{2} + 320 q^{4} + 70 q^{6} + 2 q^{8},
  - 56 q^{-1} - 286 q^{1} - 360 q^{3} - 146 q^{5} - 16 q^{7},
  \\
  & & \hspace{-2mm}
  68 + 195 q^{2} + 142 q^{4} + 33 q^{6} + 2 q^{8},
  - 39 q^{1} - 61 q^{3} - 25 q^{5} - 3 q^{7},
  9 q^{2} + 6 q^{4} + q^{6}
  \\
  & & \hspace{-17mm}
  LG_{10_{75}}
  =
  8 q^{-4} + 232 q^{-2} + 641 + 522 q^{2} + 140 q^{4} + 8 q^{6},
  - 64 q^{-3} - 399 q^{-1} - 558 q^{1} - 254 q^{3} - 31 q^{5},
  \\
  & & \hspace{-2mm}
  8 q^{-4} + 134 q^{-2} + 341 + 245 q^{2} + 49 q^{4},
  - 20 q^{-3} - 118 q^{-1} - 140 q^{1} - 42 q^{3},
  18 q^{-2} + 47 + 22 q^{2},
  \\
  & & \hspace{-2mm}
  - 7 q^{-1} - 7 q^{1},
  1
  \\
  & & \hspace{-17mm}
  LG_{10_{76}}
  =
  4 q^{-2} + 53 + 198 q^{2} + 258 q^{4} + 100 q^{6} + 6 q^{8},
  - 13 q^{-1} - 101 q^{1} - 234 q^{3} - 177 q^{5} - 31 q^{7},
  \\
  & & \hspace{-2mm}
  23 + 122 q^{2} + 180 q^{4} + 73 q^{6} + 4 q^{8},
  - 27 q^{1} - 97 q^{3} - 87 q^{5} - 17 q^{7},
  21 q^{2} + 50 q^{4} + 24 q^{6} + 2 q^{8},
  \\
  & & \hspace{-2mm}
  - 11 q^{3} - 14 q^{5} - 3 q^{7},
  3 q^{4} + q^{6}
  \\
  & & \hspace{-17mm}
  LG_{10_{77}}
  =
  28 q^{-2} + 189 + 346 q^{2} + 196 q^{4} + 28 q^{6},
  - 4 q^{-3} - 75 q^{-1} - 262 q^{1} - 271 q^{3} - 84 q^{5} - 4 q^{7},
  \\
  & & \hspace{-2mm}
  13 q^{-2} + 109 + 219 q^{2} + 129 q^{4} + 20 q^{6},
  - 20 q^{-1} - 97 q^{1} - 110 q^{3} - 35 q^{5} - 2 q^{7},
  19 + 52 q^{2} + 30 q^{4} + 4 q^{6},
  \\
  & & \hspace{-2mm}
  - 11 q^{1} - 14 q^{3} - 3 q^{5},
  3 q^{2} + q^{4}
  \\
  & & \hspace{-17mm}
  LG_{10_{78}}
  =
  8 q^{-2} + 109 + 356 q^{2} + 420 q^{4} + 152 q^{6} + 8 q^{8},
  - 28 q^{-1} - 189 q^{1} - 384 q^{3} - 266 q^{5} - 43 q^{7},
  \\
  & & \hspace{-2mm}
  44 + 192 q^{2} + 251 q^{4} + 91 q^{6} + 4 q^{8},
  - 40 q^{1} - 120 q^{3} - 93 q^{5} - 13 q^{7},
  22 q^{2} + 44 q^{4} + 15 q^{6},
  - 7 q^{3} - 7 q^{5},
  q^{4}
  \\
  & & \hspace{-17mm}
  LG_{10_{79}}
  =
  20 q^{-4} + 154 q^{-2} + 283 + 154 q^{2} + 20 q^{4},
  - 3 q^{-5} - 64 q^{-3} - 221 q^{-1} - 221 q^{1} - 64 q^{3} - 3 q^{5},
  \\
  & & \hspace{-2mm}
  16 q^{-4} + 107 q^{-2} + 194 + 107 q^{2} + 16 q^{4},
  - 2 q^{-5} - 32 q^{-3} - 107 q^{-1} - 107 q^{1} - 32 q^{3} - 2 q^{5},
  \\
  & & \hspace{-2mm}
  5 q^{-4} + 36 q^{-2} + 69 + 36 q^{2} + 5 q^{4},
  - 6 q^{-3} - 27 q^{-1} - 27 q^{1} - 6 q^{3},
  5 q^{-2} + 13 + 5 q^{2},
  - 3 q^{-1} - 3 q^{1},
  1
  \\
  & & \hspace{-17mm}
  LG_{10_{80}}
  =
  7 + 96 q^{2} + 316 q^{4} + 366 q^{6} + 128 q^{8} + 6 q^{10},
  - 25 q^{1} - 175 q^{3} - 354 q^{5} - 242 q^{7} - 38 q^{9},
  \\
  & & \hspace{-2mm}
  45 q^{2} + 203 q^{4} + 268 q^{6} + 100 q^{8} + 5 q^{10},
  - 52 q^{3} - 163 q^{5} - 134 q^{7} - 23 q^{9},
  42 q^{4} + 88 q^{6} + 39 q^{8} + 2 q^{10},
  \\
  & & \hspace{-2mm}
  - 22 q^{5} - 27 q^{7} - 5 q^{9},
  6 q^{6} + 3 q^{8}
  \\
  & & \hspace{-17mm}
  LG_{10_{81}}
  =
  56 q^{-4} + 408 q^{-2} + 731 + 408 q^{2} + 56 q^{4},
  - 7 q^{-5} - 158 q^{-3} - 543 q^{-1} - 543 q^{1} - 158 q^{3} - 7 q^{5},
  \\
  & & \hspace{-2mm}
  25 q^{-4} + 213 q^{-2} + 396 + 213 q^{2} + 25 q^{4},
  - 35 q^{-3} - 152 q^{-1} - 152 q^{1} - 35 q^{3},
  24 q^{-2} + 56 + 24 q^{2},
  \\
  & & \hspace{-2mm}
  - 8 q^{-1} - 8 q^{1},
  1
\end{eqnarray*}

\begin{eqnarray*}
  & & \hspace{-17mm}
  LG_{10_{82}}
  =
  4 q^{-4} + 72 q^{-2} + 217 + 240 q^{2} + 82 q^{4} + 4 q^{6},
  - 22 q^{-3} - 127 q^{-1} - 240 q^{1} - 161 q^{3} - 26 q^{5},
  \\
  & & \hspace{-2mm}
  3 q^{-4} + 43 q^{-2} + 148 + 191 q^{2} + 75 q^{4} + 4 q^{6},
  - 7 q^{-3} - 51 q^{-1} - 131 q^{1} - 109 q^{3} - 22 q^{5},
  \\
  & & \hspace{-2mm}
  8 q^{-2} + 49 + 86 q^{2} + 40 q^{4} + 3 q^{6},
  - 8 q^{-1} - 37 q^{1} - 36 q^{3} - 7 q^{5},
  7 + 18 q^{2} + 7 q^{4},
  - 4 q^{1} - 4 q^{3},
  q^{2}
  \\
  & & \hspace{-17mm}
  LG_{10_{83}}
  =
  12 q^{-4} + 250 q^{-2} + 641 + 492 q^{2} + 118 q^{4} + 4 q^{6},
  - 74 q^{-3} - 420 q^{-1} - 562 q^{1} - 244 q^{3} - 28 q^{5},
  \\
  & & \hspace{-2mm}
  9 q^{-4} + 153 q^{-2} + 381 + 278 q^{2} + 63 q^{4} + 3 q^{6},
  - 25 q^{-3} - 150 q^{-1} - 189 q^{1} - 71 q^{3} - 7 q^{5},
  \\
  & & \hspace{-2mm}
  27 q^{-2} + 76 + 47 q^{2} + 7 q^{4},
  - 14 q^{-1} - 18 q^{1} - 4 q^{3},
  3 + q^{2}
  \\
  & & \hspace{-17mm}
  LG_{10_{84}}
  =
  56 q^{-2} + 391 + 700 q^{2} + 394 q^{4} + 54 q^{6},
  - 7 q^{-3} - 153 q^{-1} - 528 q^{1} - 538 q^{3} - 163 q^{5} - 7 q^{7},
  \\
  & & \hspace{-2mm}
  24 q^{-2} + 213 + 416 q^{2} + 242 q^{4} + 35 q^{6},
  - 36 q^{-1} - 171 q^{1} - 191 q^{3} - 59 q^{5} - 3 q^{7},
  30 + 79 q^{2} + 46 q^{4} + 6 q^{6},
  \\
  & & \hspace{-2mm}
  - 14 q^{1} - 18 q^{3} - 4 q^{5},
  3 q^{2} + q^{4}
  \\
  & & \hspace{-17mm}
  LG_{10_{85}}
  =
  20 q^{-2} + 101 + 182 q^{2} + 146 q^{4} + 34 q^{6},
  - 3 q^{-3} - 45 q^{-1} - 137 q^{1} - 175 q^{3} - 86 q^{5} - 6 q^{7},
  \\
  & & \hspace{-2mm}
  8 q^{-2} + 61 + 142 q^{2} + 130 q^{4} + 31 q^{6},
  - 11 q^{-1} - 65 q^{1} - 118 q^{3} - 69 q^{5} - 5 q^{7},
  12 + 57 q^{2} + 75 q^{4} + 22 q^{6},
  \\
  & & \hspace{-2mm}
  - 11 q^{1} - 39 q^{3} - 31 q^{5} - 3 q^{7},
  8 q^{2} + 18 q^{4} + 6 q^{6},
  - 4 q^{3} - 4 q^{5},
  q^{4}
  \\
  & & \hspace{-17mm}
  LG_{10_{86}}
  =
  68 q^{-2} + 395 + 618 q^{2} + 304 q^{4} + 36 q^{6},
  - 9 q^{-3} - 168 q^{-1} - 497 q^{1} - 453 q^{3} - 120 q^{5} - 5 q^{7},
  \\
  & & \hspace{-2mm}
  28 q^{-2} + 214 + 379 q^{2} + 202 q^{4} + 28 q^{6},
  - 38 q^{-1} - 165 q^{1} - 176 q^{3} - 52 q^{5} - 3 q^{7},
  30 + 77 q^{2} + 44 q^{4} + 6 q^{6},
  \\
  & & \hspace{-2mm}
  - 14 q^{1} - 18 q^{3} - 4 q^{5},
  3 q^{2} + q^{4}
  \\
  & & \hspace{-17mm}
  LG_{10_{87}}
  =
  8 q^{-4} + 188 q^{-2} + 541 + 468 q^{2} + 136 q^{4} + 6 q^{6},
  - 54 q^{-3} - 341 q^{-1} - 505 q^{1} - 253 q^{3} - 35 q^{5},
  \\
  & & \hspace{-2mm}
  7 q^{-4} + 124 q^{-2} + 337 + 269 q^{2} + 70 q^{4} + 3 q^{6},
  - 21 q^{-3} - 135 q^{-1} - 180 q^{1} - 73 q^{3} - 7 q^{5},
  \\
  & & \hspace{-2mm}
  25 q^{-2} + 74 + 47 q^{2} + 7 q^{4},
  - 14 q^{-1} - 18 q^{1} - 4 q^{3},
  3 + q^{2}
  \\
  & & \hspace{-17mm}
  LG_{10_{88}}
  =
  86 q^{-4} + 618 q^{-2} + 1099 + 618 q^{2} + 86 q^{4},
  - 10 q^{-5} - 233 q^{-3} - 797 q^{-1} - 797 q^{1} - 233 q^{3} - 10 q^{5},
  \\
  & & \hspace{-2mm}
  33 q^{-4} + 290 q^{-2} + 538 + 290 q^{2} + 33 q^{4},
  - 42 q^{-3} - 185 q^{-1} - 185 q^{1} - 42 q^{3},
  26 q^{-2} + 60 + 26 q^{2},
  \\
  & & \hspace{-2mm}
  - 8 q^{-1} - 8 q^{1},
  1
  \\
  & & \hspace{-17mm}
  LG_{10_{89}}
  =
  104 q^{-2} + 641 + 1032 q^{2} + 528 q^{4} + 66 q^{6},
  - 12 q^{-3} - 257 q^{-1} - 786 q^{1} - 727 q^{3} - 194 q^{5} - 8 q^{7},
  \\
  & & \hspace{-2mm}
  36 q^{-2} + 298 + 523 q^{2} + 266 q^{4} + 29 q^{6},
  - 43 q^{-1} - 185 q^{1} - 182 q^{3} - 40 q^{5},
  26 + 60 q^{2} + 26 q^{4},
  - 8 q^{1} - 8 q^{3},
  q^{2}
  \\
  & & \hspace{-17mm}
  LG_{10_{90}}
  =
  10 q^{-4} + 198 q^{-2} + 519 + 408 q^{2} + 104 q^{4} + 4 q^{6},
  - 58 q^{-3} - 335 q^{-1} - 459 q^{1} - 208 q^{3} - 26 q^{5},
  \\
  & & \hspace{-2mm}
  7 q^{-4} + 119 q^{-2} + 307 + 232 q^{2} + 57 q^{4} + 3 q^{6},
  - 19 q^{-3} - 119 q^{-1} - 156 q^{1} - 63 q^{3} - 7 q^{5},
  \\
  & & \hspace{-2mm}
  21 q^{-2} + 63 + 41 q^{2} + 7 q^{4},
  - 12 q^{-1} - 16 q^{1} - 4 q^{3},
  3 + q^{2}
  \\
  & & \hspace{-17mm}
  LG_{10_{91}}
  =
  32 q^{-4} + 218 q^{-2} + 387 + 212 q^{2} + 28 q^{4},
  - 5 q^{-5} - 95 q^{-3} - 308 q^{-1} - 305 q^{1} - 91 q^{3} - 4 q^{5},
  \\
  & & \hspace{-2mm}
  25 q^{-4} + 154 q^{-2} + 274 + 155 q^{2} + 24 q^{4},
  - 3 q^{-5} - 47 q^{-3} - 155 q^{-1} - 158 q^{1} - 50 q^{3} - 3 q^{5},
  \\
  & & \hspace{-2mm}
  7 q^{-4} + 52 q^{-2} + 103 + 57 q^{2} + 8 q^{4},
  - 8 q^{-3} - 40 q^{-1} - 42 q^{1} - 10 q^{3},
  7 q^{-2} + 19 + 8 q^{2},
  - 4 q^{-1} - 4 q^{1},
  1
  \\
  & & \hspace{-17mm}
  LG_{10_{92}}
  =
  10 q^{-2} + 157 + 552 q^{2} + 660 q^{4} + 242 q^{6} + 12 q^{8},
  - 38 q^{-1} - 289 q^{1} - 611 q^{3} - 431 q^{5} - 71 q^{7},
  \\
  & & \hspace{-2mm}
  66 + 315 q^{2} + 426 q^{4} + 165 q^{6} + 8 q^{8},
  - 68 q^{1} - 218 q^{3} - 182 q^{5} - 32 q^{7},
  44 q^{2} + 93 q^{4} + 41 q^{6} + 2 q^{8},
  \\
  & & \hspace{-2mm}
  - 17 q^{3} - 20 q^{5} - 3 q^{7},
  3 q^{4} + q^{6}
  \\
  & & \hspace{-17mm}
  LG_{10_{93}}
  =
  24 q^{-4} + 160 q^{-2} + 341 + 286 q^{2} + 64 q^{4},
  - 3 q^{-5} - 60 q^{-3} - 228 q^{-1} - 320 q^{1} - 160 q^{3} - 11 q^{5},
  \\
  & & \hspace{-2mm}
  8 q^{-4} + 80 q^{-2} + 211 + 206 q^{2} + 52 q^{4},
  - 10 q^{-3} - 70 q^{-1} - 137 q^{1} - 84 q^{3} - 7 q^{5},
  8 q^{-2} + 42 + 58 q^{2} + 16 q^{4},
  \\
  & & \hspace{-2mm}
  - 4 q^{-1} - 16 q^{1} - 12 q^{3},
  1 + 3 q^{2}
  \\
  & & \hspace{-17mm}
  LG_{10_{94}}
  =
  4 q^{-4} + 86 q^{-2} + 283 + 322 q^{2} + 112 q^{4} + 6 q^{6},
  - 24 q^{-3} - 161 q^{-1} - 317 q^{1} - 215 q^{3} - 35 q^{5},
  \\
  & & \hspace{-2mm}
  3 q^{-4} + 52 q^{-2} + 193 + 248 q^{2} + 95 q^{4} + 5 q^{6},
  - 8 q^{-3} - 67 q^{-1} - 167 q^{1} - 133 q^{3} - 25 q^{5},
  \\
  & & \hspace{-2mm}
  11 q^{-2} + 63 + 102 q^{2} + 44 q^{4} + 3 q^{6},
  - 11 q^{-1} - 43 q^{1} - 39 q^{3} - 7 q^{5},
  8 + 19 q^{2} + 7 q^{4},
  - 4 q^{1} - 4 q^{3},
  q^{2}
  \\
  & & \hspace{-17mm}
  LG_{10_{95}}
  =
  78 q^{-2} + 479 + 778 q^{2} + 398 q^{4} + 48 q^{6},
  - 10 q^{-3} - 198 q^{-1} - 610 q^{1} - 570 q^{3} - 154 q^{5} - 6 q^{7},
  \\
  & & \hspace{-2mm}
  32 q^{-2} + 254 + 454 q^{2} + 241 q^{4} + 30 q^{6},
  - 44 q^{-1} - 190 q^{1} - 197 q^{3} - 53 q^{5} - 2 q^{7},
  34 + 83 q^{2} + 44 q^{4} + 4 q^{6},
  \\
  & & \hspace{-2mm}
  - 15 q^{1} - 18 q^{3} - 3 q^{5},
  3 q^{2} + q^{4}
  \\
  & & \hspace{-17mm}
  LG_{10_{96}}
  =
  18 q^{-4} + 346 q^{-2} + 895 + 706 q^{2} + 182 q^{4} + 10 q^{6},
  - 95 q^{-3} - 553 q^{-1} - 753 q^{1} - 334 q^{3} - 39 q^{5},
  \\
  & & \hspace{-2mm}
  8 q^{-4} + 169 q^{-2} + 436 + 315 q^{2} + 62 q^{4},
  - 20 q^{-3} - 135 q^{-1} - 167 q^{1} - 52 q^{3},
  18 q^{-2} + 50 + 25 q^{2},
  \\
  & & \hspace{-2mm}
  - 7 q^{-1} - 7 q^{1},
  1
  \\
  & & \hspace{-17mm}
  LG_{10_{97}}
  =
  26 q^{-2} + 343 + 864 q^{2} + 686 q^{4} + 180 q^{6} + 8 q^{8},
  - 88 q^{-1} - 507 q^{1} - 707 q^{3} - 329 q^{5} - 41 q^{7},
  \\
  & & \hspace{-2mm}
  113 + 350 q^{2} + 281 q^{4} + 68 q^{6} + 2 q^{8},
  - 65 q^{1} - 107 q^{3} - 45 q^{5} - 3 q^{7},
  14 q^{2} + 10 q^{4} + q^{6}
  \\
  & & \hspace{-17mm}
  LG_{10_{98}}
  =
  10 q^{-2} + 141 + 462 q^{2} + 534 q^{4} + 190 q^{6} + 10 q^{8},
  - 36 q^{-1} - 250 q^{1} - 503 q^{3} - 344 q^{5} - 55 q^{7},
  \\
  & & \hspace{-2mm}
  60 + 267 q^{2} + 348 q^{4} + 129 q^{6} + 6 q^{8},
  - 60 q^{1} - 183 q^{3} - 148 q^{5} - 25 q^{7},
  38 q^{2} + 79 q^{4} + 34 q^{6} + 2 q^{8},
  \\
  & & \hspace{-2mm}
  - 15 q^{3} - 18 q^{5} - 3 q^{7},
  3 q^{4} + q^{6}
  \\
  & & \hspace{-17mm}
  LG_{10_{99}}
  =
  36 q^{-4} + 272 q^{-2} + 491 + 272 q^{2} + 36 q^{4},
  - 5 q^{-5} - 115 q^{-3} - 388 q^{-1} - 388 q^{1} - 115 q^{3} - 5 q^{5},
  \\
  & & \hspace{-2mm}
  28 q^{-4} + 192 q^{-2} + 344 + 192 q^{2} + 28 q^{4},
  - 3 q^{-5} - 57 q^{-3} - 192 q^{-1} - 192 q^{1} - 57 q^{3} - 3 q^{5},
  \\
  & & \hspace{-2mm}
  8 q^{-4} + 64 q^{-2} + 122 + 64 q^{2} + 8 q^{4},
  - 10 q^{-3} - 46 q^{-1} - 46 q^{1} - 10 q^{3},
  8 q^{-2} + 20 + 8 q^{2},
  - 4 q^{-1} - 4 q^{1},
  1
\end{eqnarray*}

\begin{eqnarray*}
  & & \hspace{-17mm}
  LG_{10_{100}}
  =
  22 q^{-2} + 129 + 252 q^{2} + 204 q^{4} + 46 q^{6},
  - 3 q^{-3} - 54 q^{-1} - 184 q^{1} - 242 q^{3} - 117 q^{5} - 8 q^{7},
  \\
  & & \hspace{-2mm}
  9 q^{-2} + 79 + 191 q^{2} + 173 q^{4} + 40 q^{6},
  - 14 q^{-1} - 86 q^{1} - 152 q^{3} - 86 q^{5} - 6 q^{7},
  16 + 72 q^{2} + 90 q^{4} + 25 q^{6},
  \\
  & & \hspace{-2mm}
  - 14 q^{1} - 45 q^{3} - 34 q^{5} - 3 q^{7},
  9 q^{2} + 19 q^{4} + 6 q^{6},
  - 4 q^{3} - 4 q^{5},
  q^{4}
  \\
  & & \hspace{-17mm}
  LG_{10_{101}}
  =
  25 + 308 q^{2} + 752 q^{4} + 582 q^{6} + 148 q^{8} + 6 q^{10},
  - 84 q^{1} - 464 q^{3} - 636 q^{5} - 292 q^{7} - 36 q^{9},
  \\
  & & \hspace{-2mm}
  117 q^{2} + 358 q^{4} + 293 q^{6} + 76 q^{8} + 3 q^{10},
  - 80 q^{3} - 139 q^{5} - 67 q^{7} - 8 q^{9},
  23 q^{4} + 21 q^{6} + 5 q^{8}
  \\
  & & \hspace{-17mm}
  LG_{10_{102}}
  =
  8 q^{-4} + 166 q^{-2} + 447 + 364 q^{2} + 100 q^{4} + 4 q^{6},
  - 50 q^{-3} - 288 q^{-1} - 404 q^{1} - 192 q^{3} - 26 q^{5},
  \\
  & & \hspace{-2mm}
  7 q^{-4} + 107 q^{-2} + 274 + 210 q^{2} + 55 q^{4} + 3 q^{6},
  - 19 q^{-3} - 111 q^{-1} - 144 q^{1} - 59 q^{3} - 7 q^{5},
  \\
  & & \hspace{-2mm}
  21 q^{-2} + 61 + 39 q^{2} + 7 q^{4},
  - 12 q^{-1} - 16 q^{1} - 4 q^{3},
  3 + q^{2}
  \\
  & & \hspace{-17mm}
  LG_{10_{103}}
  =
  52 q^{-2} + 309 + 502 q^{2} + 258 q^{4} + 34 q^{6},
  - 7 q^{-3} - 130 q^{-1} - 396 q^{1} - 374 q^{3} - 106 q^{5} - 5 q^{7},
  \\
  & & \hspace{-2mm}
  22 q^{-2} + 168 + 307 q^{2} + 170 q^{4} + 26 q^{6},
  - 30 q^{-1} - 132 q^{1} - 145 q^{3} - 46 q^{5} - 3 q^{7},
  \\
  & & \hspace{-2mm}
  24 + 64 q^{2} + 38 q^{4} + 6 q^{6},
  - 12 q^{1} - 16 q^{3} - 4 q^{5},
  3 q^{2} + q^{4}
  \\
  & & \hspace{-17mm}
  LG_{10_{104}}
  =
  30 q^{-4} + 244 q^{-2} + 449 + 250 q^{2} + 34 q^{4},
  - 4 q^{-5} - 101 q^{-3} - 351 q^{-1} - 354 q^{1} - 105 q^{3} - 5 q^{5},
  \\
  & & \hspace{-2mm}
  25 q^{-4} + 173 q^{-2} + 309 + 172 q^{2} + 26 q^{4},
  - 3 q^{-5} - 53 q^{-3} - 173 q^{-1} - 170 q^{1} - 50 q^{3} - 3 q^{5},
  \\
  & & \hspace{-2mm}
  8 q^{-4} + 60 q^{-2} + 109 + 55 q^{2} + 7 q^{4},
  - 10 q^{-3} - 43 q^{-1} - 41 q^{1} - 8 q^{3},
  8 q^{-2} + 19 + 7 q^{2},
  - 4 q^{-1} - 4 q^{1},
  1
  \\
  & & \hspace{-17mm}
  LG_{10_{105}}
  =
  10 q^{-4} + 174 q^{-2} + 645 + 794 q^{2} + 300 q^{4} + 16 q^{6},
  - 39 q^{-3} - 319 q^{-1} - 697 q^{1} - 503 q^{3} - 86 q^{5},
  \\
  & & \hspace{-2mm}
  64 q^{-2} + 317 + 435 q^{2} + 168 q^{4} + 8 q^{6},
  - 57 q^{-1} - 183 q^{1} - 148 q^{3} - 22 q^{5},
  29 + 58 q^{2} + 21 q^{4},
  - 8 q^{1} - 8 q^{3},
  q^{2}
  \\
  & & \hspace{-17mm}
  LG_{10_{106}}
  =
  4 q^{-4} + 94 q^{-2} + 327 + 376 q^{2} + 128 q^{4} + 6 q^{6},
  - 25 q^{-3} - 182 q^{-1} - 367 q^{1} - 248 q^{3} - 38 q^{5},
  \\
  & & \hspace{-2mm}
  3 q^{-4} + 56 q^{-2} + 219 + 281 q^{2} + 105 q^{4} + 5 q^{6},
  - 8 q^{-3} - 73 q^{-1} - 184 q^{1} - 145 q^{3} - 26 q^{5},
  \\
  & & \hspace{-2mm}
  11 q^{-2} + 67 + 108 q^{2} + 46 q^{4} + 3 q^{6},
  - 11 q^{-1} - 44 q^{1} - 40 q^{3} - 7 q^{5},
  8 + 19 q^{2} + 7 q^{4},
  - 4 q^{1} - 4 q^{3},
  q^{2}
  \\
  & & \hspace{-17mm}
  LG_{10_{107}}
  =
  78 q^{-4} + 524 q^{-2} + 901 + 490 q^{2} + 66 q^{4},
  - 10 q^{-5} - 207 q^{-3} - 672 q^{-1} - 651 q^{1} - 184 q^{3} - 8 q^{5},
  \\
  & & \hspace{-2mm}
  33 q^{-4} + 260 q^{-2} + 463 + 241 q^{2} + 27 q^{4},
  - 42 q^{-3} - 171 q^{-1} - 165 q^{1} - 36 q^{3},
  26 q^{-2} + 58 + 24 q^{2},
  \\
  & & \hspace{-2mm}
  - 8 q^{-1} - 8 q^{1},
  1
  \\
  & & \hspace{-17mm}
  LG_{10_{108}}
  =
  24 q^{-4} + 148 q^{-2} + 293 + 236 q^{2} + 52 q^{4},
  - 3 q^{-5} - 58 q^{-3} - 202 q^{-1} - 271 q^{1} - 133 q^{3} - 9 q^{5},
  \\
  & & \hspace{-2mm}
  8 q^{-4} + 74 q^{-2} + 185 + 179 q^{2} + 46 q^{4},
  - 10 q^{-3} - 64 q^{-1} - 125 q^{1} - 78 q^{3} - 7 q^{5},
  8 q^{-2} + 40 + 56 q^{2} + 16 q^{4},
  \\
  & & \hspace{-2mm}
  - 4 q^{-1} - 16 q^{1} - 12 q^{3},
  1 + 3 q^{2}
  \\
  & & \hspace{-17mm}
  LG_{10_{109}}
  =
  38 q^{-4} + 308 q^{-2} + 561 + 308 q^{2} + 38 q^{4},
  - 5 q^{-5} - 126 q^{-3} - 440 q^{-1} - 440 q^{1} - 126 q^{3} - 5 q^{5},
  \\
  & & \hspace{-2mm}
  29 q^{-4} + 212 q^{-2} + 383 + 212 q^{2} + 29 q^{4},
  - 3 q^{-5} - 60 q^{-3} - 208 q^{-1} - 208 q^{1} - 60 q^{3} - 3 q^{5},
  \\
  & & \hspace{-2mm}
  8 q^{-4} + 67 q^{-2} + 128 + 67 q^{2} + 8 q^{4},
  - 10 q^{-3} - 47 q^{-1} - 47 q^{1} - 10 q^{3},
  8 q^{-2} + 20 + 8 q^{2},
  - 4 q^{-1} - 4 q^{1},
  1
  \\
  & & \hspace{-17mm}
  LG_{10_{110}}
  =
  10 q^{-4} + 158 q^{-2} + 537 + 622 q^{2} + 218 q^{4} + 10 q^{6},
  - 37 q^{-3} - 277 q^{-1} - 568 q^{1} - 391 q^{3} - 63 q^{5},
  \\
  & & \hspace{-2mm}
  58 q^{-2} + 271 + 364 q^{2} + 139 q^{4} + 8 q^{6},
  - 51 q^{-1} - 163 q^{1} - 134 q^{3} - 22 q^{5},
  27 + 56 q^{2} + 21 q^{4},
  - 8 q^{1} - 8 q^{3},
  q^{2}
  \\
  & & \hspace{-17mm}
  LG_{10_{111}}
  =
  10 q^{-2} + 137 + 416 q^{2} + 460 q^{4} + 158 q^{6} + 8 q^{8},
  - 36 q^{-1} - 233 q^{1} - 445 q^{3} - 295 q^{5} - 47 q^{7},
  \\
  & & \hspace{-2mm}
  58 + 243 q^{2} + 312 q^{4} + 116 q^{6} + 6 q^{8},
  - 56 q^{1} - 170 q^{3} - 139 q^{5} - 25 q^{7},
  36 q^{2} + 77 q^{4} + 34 q^{6} + 2 q^{8},
  \\
  & & \hspace{-2mm}
  - 15 q^{3} - 18 q^{5} - 3 q^{7},
  3 q^{4} + q^{6}
  \\
  & & \hspace{-17mm}
  LG_{10_{112}}
  =
  6 q^{-4} + 130 q^{-2} + 431 + 488 q^{2} + 170 q^{4} + 8 q^{6},
  - 36 q^{-3} - 245 q^{-1} - 482 q^{1} - 325 q^{3} - 52 q^{5},
  \\
  & & \hspace{-2mm}
  4 q^{-4} + 78 q^{-2} + 294 + 373 q^{2} + 143 q^{4} + 7 q^{6},
  - 11 q^{-3} - 100 q^{-1} - 251 q^{1} - 199 q^{3} - 37 q^{5},
  \\
  & & \hspace{-2mm}
  15 q^{-2} + 93 + 151 q^{2} + 66 q^{4} + 4 q^{6},
  - 15 q^{-1} - 62 q^{1} - 57 q^{3} - 10 q^{5},
  11 + 26 q^{2} + 10 q^{4},
  - 5 q^{1} - 5 q^{3},
  q^{2}
  \\
  & & \hspace{-17mm}
  LG_{10_{113}}
  =
  94 q^{-2} + 663 + 1178 q^{2} + 666 q^{4} + 90 q^{6},
  - 11 q^{-3} - 257 q^{-1} - 884 q^{1} - 898 q^{3} - 271 q^{5} - 11 q^{7},
  \\
  & & \hspace{-2mm}
  38 q^{-2} + 349 + 675 q^{2} + 391 q^{4} + 53 q^{6},
  - 56 q^{-1} - 266 q^{1} - 293 q^{3} - 86 q^{5} - 3 q^{7},
  44 + 112 q^{2} + 63 q^{4} + 6 q^{6},
  \\
  & & \hspace{-2mm}
  - 18 q^{1} - 22 q^{3} - 4 q^{5},
  3 q^{2} + q^{4}
  \\
  & & \hspace{-17mm}
  LG_{10_{114}}
  =
  14 q^{-4} + 288 q^{-2} + 753 + 602 q^{2} + 156 q^{4} + 6 q^{6},
  - 85 q^{-3} - 488 q^{-1} - 673 q^{1} - 309 q^{3} - 39 q^{5},
  \\
  & & \hspace{-2mm}
  10 q^{-4} + 176 q^{-2} + 449 + 342 q^{2} + 84 q^{4} + 4 q^{6},
  - 28 q^{-3} - 173 q^{-1} - 226 q^{1} - 91 q^{3} - 10 q^{5},
  \\
  & & \hspace{-2mm}
  30 q^{-2} + 87 + 57 q^{2} + 10 q^{4},
  - 15 q^{-1} - 20 q^{1} - 5 q^{3},
  3 + q^{2}
  \\
  & & \hspace{-17mm}
  LG_{10_{115}}
  =
  100 q^{-4} + 712 q^{-2} + 1261 + 712 q^{2} + 100 q^{4},
  - 12 q^{-5} - 272 q^{-3} - 921 q^{-1} - 921 q^{1} - 272 q^{3} - 12 q^{5},
  \\
  & & \hspace{-2mm}
  40 q^{-4} + 342 q^{-2} + 630 + 342 q^{2} + 40 q^{4},
  - 51 q^{-3} - 220 q^{-1} - 220 q^{1} - 51 q^{3},
  31 q^{-2} + 71 + 31 q^{2},
  \\
  & & \hspace{-2mm}
  - 9 q^{-1} - 9 q^{1},
  1
  \\
  & & \hspace{-17mm}
  LG_{10_{116}}
  =
  6 q^{-4} + 150 q^{-2} + 525 + 602 q^{2} + 210 q^{4} + 10 q^{6},
  - 39 q^{-3} - 293 q^{-1} - 590 q^{1} - 399 q^{3} - 63 q^{5},
  \\
  & & \hspace{-2mm}
  4 q^{-4} + 90 q^{-2} + 356 + 451 q^{2} + 170 q^{4} + 8 q^{6},
  - 12 q^{-3} - 120 q^{-1} - 299 q^{1} - 232 q^{3} - 41 q^{5},
  \\
  & & \hspace{-2mm}
  18 q^{-2} + 110 + 172 q^{2} + 72 q^{4} + 4 q^{6},
  - 18 q^{-1} - 69 q^{1} - 61 q^{3} - 10 q^{5},
  12 + 27 q^{2} + 10 q^{4},
  - 5 q^{1} - 5 q^{3},
  q^{2}
  \\
  & & \hspace{-17mm}
  LG_{10_{117}}
  =
  90 q^{-2} + 593 + 1014 q^{2} + 552 q^{4} + 72 q^{6},
  - 11 q^{-3} - 236 q^{-1} - 772 q^{1} - 759 q^{3} - 221 q^{5} - 9 q^{7},
  \\
  & & \hspace{-2mm}
  36 q^{-2} + 308 + 578 q^{2} + 326 q^{4} + 44 q^{6},
  - 50 q^{-1} - 229 q^{1} - 249 q^{3} - 73 q^{5} - 3 q^{7},
  38 + 97 q^{2} + 55 q^{4} + 6 q^{6},
  \\
  & & \hspace{-2mm}
  - 16 q^{1} - 20 q^{3} - 4 q^{5},
  3 q^{2} + q^{4}
\end{eqnarray*}

\begin{eqnarray*}
  & & \hspace{-17mm}
  LG_{10_{118}}
  =
  52 q^{-4} + 392 q^{-2} + 703 + 392 q^{2} + 52 q^{4},
  - 7 q^{-5} - 166 q^{-3} - 557 q^{-1} - 557 q^{1} - 166 q^{3} - 7 q^{5},
  \\
  & & \hspace{-2mm}
  40 q^{-4} + 277 q^{-2} + 493 + 277 q^{2} + 40 q^{4},
  - 4 q^{-5} - 82 q^{-3} - 276 q^{-1} - 276 q^{1} - 82 q^{3} - 4 q^{5},
  \\
  & & \hspace{-2mm}
  11 q^{-4} + 92 q^{-2} + 174 + 92 q^{2} + 11 q^{4},
  - 14 q^{-3} - 65 q^{-1} - 65 q^{1} - 14 q^{3},
  11 q^{-2} + 27 + 11 q^{2},
  \\
  & & \hspace{-2mm}
  - 5 q^{-1} - 5 q^{1},
  1
  \\
  & & \hspace{-17mm}
  LG_{10_{119}}
  =
  16 q^{-4} + 342 q^{-2} + 917 + 748 q^{2} + 196 q^{4} + 8 q^{6},
  - 98 q^{-3} - 583 q^{-1} - 819 q^{1} - 380 q^{3} - 46 q^{5},
  \\
  & & \hspace{-2mm}
  11 q^{-4} + 203 q^{-2} + 528 + 403 q^{2} + 93 q^{4} + 3 q^{6},
  - 31 q^{-3} - 196 q^{-1} - 254 q^{1} - 96 q^{3} - 7 q^{5},
  \\
  & & \hspace{-2mm}
  33 q^{-2} + 94 + 58 q^{2} + 7 q^{4},
  - 16 q^{-1} - 20 q^{1} - 4 q^{3},
  3 + q^{2}
  \\
  & & \hspace{-17mm}
  LG_{10_{120}}
  =
  41 + 508 q^{2} + 1192 q^{4} + 890 q^{6} + 210 q^{8} + 8 q^{10},
  - 140 q^{1} - 744 q^{3} - 983 q^{5} - 426 q^{7} - 47 q^{9},
  \\
  & & \hspace{-2mm}
  189 q^{2} + 548 q^{4} + 429 q^{6} + 99 q^{8} + 3 q^{10},
  - 120 q^{3} - 200 q^{5} - 88 q^{7} - 8 q^{9},
  31 q^{4} + 28 q^{6} + 5 q^{8}
  \\
  & & \hspace{-17mm}
  LG_{10_{121}}
  =
  124 q^{-2} + 783 + 1278 q^{2} + 666 q^{4} + 82 q^{6},
  - 15 q^{-3} - 319 q^{-1} - 991 q^{1} - 935 q^{3} - 258 q^{5} - 10 q^{7},
  \\
  & & \hspace{-2mm}
  49 q^{-2} + 402 + 719 q^{2} + 388 q^{4} + 49 q^{6},
  - 66 q^{-1} - 285 q^{1} - 298 q^{3} - 82 q^{5} - 3 q^{7},
  47 + 114 q^{2} + 62 q^{4} + 6 q^{6},
  \\
  & & \hspace{-2mm}
  - 18 q^{1} - 22 q^{3} - 4 q^{5},
  3 q^{2} + q^{4}
  \\
  & & \hspace{-17mm}
  LG_{10_{122}}
  =
  14 q^{-4} + 336 q^{-2} + 951 + 820 q^{2} + 232 q^{4} + 10 q^{6},
  - 96 q^{-3} - 596 q^{-1} - 877 q^{1} - 434 q^{3} - 57 q^{5},
  \\
  & & \hspace{-2mm}
  12 q^{-4} + 211 q^{-2} + 564 + 450 q^{2} + 113 q^{4} + 4 q^{6},
  - 34 q^{-3} - 212 q^{-1} - 282 q^{1} - 114 q^{3} - 10 q^{5},
  \\
  & & \hspace{-2mm}
  36 q^{-2} + 104 + 67 q^{2} + 10 q^{4},
  - 17 q^{-1} - 22 q^{1} - 5 q^{3},
  3 + q^{2}
  \\
  & & \hspace{-17mm}
  LG_{10_{123}}
  =
  80 q^{-4} + 622 q^{-2} + 1113 + 622 q^{2} + 80 q^{4},
  - 10 q^{-5} - 260 q^{-3} - 882 q^{-1} - 882 q^{1} - 260 q^{3} - 10 q^{5},
  \\
  & & \hspace{-2mm}
  59 q^{-4} + 435 q^{-2} + 776 + 435 q^{2} + 59 q^{4},
  - 5 q^{-5} - 124 q^{-3} - 429 q^{-1} - 429 q^{1} - 124 q^{3} - 5 q^{5},
  \\
  & & \hspace{-2mm}
  15 q^{-4} + 139 q^{-2} + 263 + 139 q^{2} + 15 q^{4},
  - 20 q^{-3} - 94 q^{-1} - 94 q^{1} - 20 q^{3},
  15 q^{-2} + 36 + 15 q^{2},
  \\
  & & \hspace{-2mm}
  - 6 q^{-1} - 6 q^{1},
  1
  \\
  & & \hspace{-17mm}
  LG_{10_{124}}
  =
  1 + 2 q^{2} + 2 q^{4} + 2 q^{6},
  - q^{1} - q^{3} - q^{5} - q^{7},
  - q^{4},
  q^{3} + 2 q^{5} + q^{7},
  - q^{4} - 2 q^{6} - q^{8},
  q^{5} + q^{7},
  q^{8},
  - q^{7} - q^{9},
  q^{8}
  \\
  & & \hspace{-17mm}
  LG_{10_{125}}
  =
  4 q^{-4} + 8 q^{-2} + 5 + 2 q^{2},
  - q^{-5} - 5 q^{-3} - 6 q^{-1} - 3 q^{1} - q^{3},
  q^{-4} + 3 q^{-2} + 4 + 4 q^{2},
  - 2 q^{-1} - 5 q^{1} - 3 q^{3},
  \\
  & & \hspace{-2mm}
  2 + 5 q^{2} + q^{4},
  - 2 q^{1} - 2 q^{3},
  q^{2}
  \\
  & & \hspace{-17mm}
  LG_{10_{126}}
  =
  11 + 32 q^{2} + 20 q^{4} + 4 q^{6},
  - 3 q^{-1} - 20 q^{1} - 26 q^{3} - 10 q^{5} - q^{7},
  q^{-2} + 7 + 19 q^{2} + 13 q^{4} + 4 q^{6},
  \\
  & & \hspace{-2mm}
  - 2 q^{-1} - 8 q^{1} - 10 q^{3} - 5 q^{5} - q^{7},
  2 + 6 q^{2} + 3 q^{4} + q^{6},
  - 2 q^{1} - 2 q^{3},
  q^{2}
  \\
  & & \hspace{-17mm}
  LG_{10_{127}}
  =
  9 + 46 q^{2} + 70 q^{4} + 28 q^{6} + 2 q^{8},
  - q^{-1} - 20 q^{1} - 59 q^{3} - 50 q^{5} - 10 q^{7},
  3 + 26 q^{2} + 50 q^{4} + 23 q^{6} + 2 q^{8},
  \\
  & & \hspace{-2mm}
  - 4 q^{1} - 24 q^{3} - 27 q^{5} - 7 q^{7},
  4 q^{2} + 15 q^{4} + 8 q^{6} + q^{8},
  - 3 q^{3} - 4 q^{5} - q^{7},
  q^{4}
  \\
  & & \hspace{-17mm}
  LG_{10_{128}}
  =
  3 + 12 q^{2} + 12 q^{4} + 2 q^{6},
  - 5 q^{1} - 9 q^{3} - 3 q^{5} + q^{7},
  2 q^{2} - 2 q^{6} - q^{8},
  q^{3} - q^{5} - 2 q^{7},
  2 q^{4} + 8 q^{6} + 3 q^{8},
  \\
  & & \hspace{-2mm}
  - 5 q^{5} - 6 q^{7} - q^{9},
  3 q^{6} + q^{8}
  \\
  & & \hspace{-17mm}
  LG_{10_{129}}
  =
  20 q^{-2} + 69 + 56 q^{2} + 16 q^{4},
  - 5 q^{-3} - 39 q^{-1} - 58 q^{1} - 27 q^{3} - 3 q^{5},
  q^{-4} + 12 q^{-2} + 33 + 21 q^{2} + 5 q^{4},
  \\
  & & \hspace{-2mm}
  - 2 q^{-3} - 10 q^{-1} - 10 q^{1} - 2 q^{3},
  q^{-2} + 3 
  \\
  & & \hspace{-17mm}
  LG_{10_{130}}
  =
  13 + 24 q^{2} + 22 q^{4} + 6 q^{6},
  - 7 q^{-1} - 18 q^{1} - 20 q^{3} - 10 q^{5} - q^{7},
  3 q^{-2} + 14 + 12 q^{2} + 6 q^{4} + q^{6},
  \\
  & & \hspace{-2mm}
  - 6 q^{-1} - 8 q^{1} - 2 q^{3},
  3 + q^{2}
  \\
  & & \hspace{-17mm}
  LG_{10_{131}}
  =
  25 + 102 q^{2} + 96 q^{4} + 32 q^{6} + 2 q^{8},
  - 3 q^{-1} - 52 q^{1} - 92 q^{3} - 52 q^{5} - 9 q^{7},
  9 + 45 q^{2} + 40 q^{4} + 13 q^{6} + q^{8},
  \\
  & & \hspace{-2mm}
  - 9 q^{1} - 15 q^{3} - 7 q^{5} - q^{7},
  3 q^{2} + q^{4}
  \\
  & & \hspace{-17mm}
  LG_{10_{132}}
  =
  - 1 - 2 q^{2} + 4 q^{4} + 4 q^{6},
  - q^{-1} + 2 q^{1} - 4 q^{5} - q^{7},
  q^{-2} + 2 - q^{2} + q^{6},
  - 2 q^{-1} - q^{1} + q^{3},
  1
  \\
  & & \hspace{-17mm}
  LG_{10_{133}}
  =
  9 + 36 q^{2} + 36 q^{4} + 18 q^{6} + 2 q^{8},
  - q^{-1} - 18 q^{1} - 32 q^{3} - 22 q^{5} - 7 q^{7},
  3 + 15 q^{2} + 13 q^{4} + 7 q^{6} + q^{8},
  \\
  & & \hspace{-2mm}
  - 3 q^{1} - 4 q^{3} - 2 q^{5} - q^{7},
  q^{2}
  \\
  & & \hspace{-17mm}
  LG_{10_{134}}
  =
  3 + 20 q^{2} + 34 q^{4} + 22 q^{6} + 2 q^{8},
  - 7 q^{1} - 25 q^{3} - 29 q^{5} - 11 q^{7},
  8 q^{2} + 22 q^{4} + 22 q^{6} + 4 q^{8},
  \\
  & & \hspace{-2mm}
  - 7 q^{3} - 21 q^{5} - 15 q^{7} - q^{9},
  8 q^{4} + 18 q^{6} + 6 q^{8},
  - 7 q^{5} - 8 q^{7} - q^{9},
  3 q^{6} + q^{8}
  \\
  & & \hspace{-17mm}
  LG_{10_{135}}
  =
  2 q^{-4} + 60 q^{-2} + 159 + 108 q^{2} + 20 q^{4},
  - 15 q^{-3} - 99 q^{-1} - 126 q^{1} - 45 q^{3} - 3 q^{5},
  \\
  & & \hspace{-2mm}
  q^{-4} + 28 q^{-2} + 76 + 47 q^{2} + 7 q^{4},
  - 2 q^{-3} - 21 q^{-1} - 25 q^{1} - 6 q^{3},
  q^{-2} + 6 + 2 q^{2}
  \\
  & & \hspace{-17mm}
  LG_{10_{136}}
  =
  4 q^{-4} + 18 q^{-2} + 19 + 14 q^{2} + 4 q^{4},
  - 8 q^{-3} - 15 q^{-1} - 15 q^{1} - 9 q^{3} - q^{5},
  5 q^{-2} + 8 + 10 q^{2} + 3 q^{4},
  \\
  & & \hspace{-2mm}
  - q^{-1} - 4 q^{1} - 3 q^{3},
  q^{2}
  \\
  & & \hspace{-17mm}
  LG_{10_{137}}
  =
  16 q^{-2} + 69 + 72 q^{2} + 34 q^{4} + 4 q^{6},
  - 2 q^{-3} - 32 q^{-1} - 60 q^{1} - 40 q^{3} - 10 q^{5},
  5 q^{-2} + 24 + 21 q^{2} + 8 q^{4},
  \\
  & & \hspace{-2mm}
  - 4 q^{-1} - 6 q^{1} - 2 q^{3},
  1
  \\
  & & \hspace{-17mm}
  LG_{10_{138}}
  =
  4 q^{-4} + 44 q^{-2} + 99 + 74 q^{2} + 8 q^{4},
  12 q^{-3} + 62 q^{-1} + 89 q^{1} + 39 q^{3},
  15 q^{-2} + 55 + 64 q^{2} + 16 q^{4},
  \\
  & & \hspace{-2mm}
  13 q^{-1} + 43 q^{1} + 34 q^{3} + 4 q^{5},
  10 + 23 q^{2} + 8 q^{4},
  5 q^{1} + 5 q^{3},
  q^{2}
  \\
  & & \hspace{-17mm}
  LG_{10_{139}}
  =
  1 + 2 q^{2} + 4 q^{4} + 10 q^{6} + 4 q^{8},
  - q^{1} - q^{3} - 5 q^{5} - 6 q^{7} - q^{9},
  - 2 q^{4} + q^{6} + q^{8},
  2 q^{3} + 5 q^{5} + 3 q^{7},
  \\
  & & \hspace{-2mm}
  - 3 q^{4} - 5 q^{6} - 2 q^{8},
  2 q^{5} + 2 q^{7},
  q^{8},
  - q^{7} - q^{9},
  q^{8}
  \\
  & & \hspace{-17mm}
  LG_{10_{140}}
  =
  3 + 4 q^{2} + 8 q^{4} + 4 q^{6},
  - 2 q^{-1} - 3 q^{1} - 5 q^{3} - 5 q^{5} - q^{7},
  q^{-2} + 4 + 2 q^{2} + 2 q^{4} + q^{6},
  - 2 q^{-1} - 2 q^{1},
  1
  \\
  & & \hspace{-17mm}
  LG_{10_{141}}
  =
  2 q^{-2} + 23 + 32 q^{2} + 18 q^{4} + 2 q^{6},
  - 12 q^{-1} - 28 q^{1} - 23 q^{3} - 7 q^{5},
  5 q^{-2} + 21 + 19 q^{2} + 8 q^{4} + q^{6},
  \\
  & & \hspace{-2mm}
  - q^{-3} - 12 q^{-1} - 15 q^{1} - 5 q^{3} - q^{5},
  3 q^{-2} + 10 + 4 q^{2},
  - 3 q^{-1} - 3 q^{1},
  1
\end{eqnarray*}

\begin{eqnarray*}
  & & \hspace{-17mm}
  LG_{10_{142}}
  =
  3 + 12 q^{2} + 14 q^{4} + 6 q^{6},
  - 5 q^{1} - 12 q^{3} - 9 q^{5} - 2 q^{7},
  4 q^{2} + 8 q^{4} + 6 q^{6},
  - 3 q^{3} - 8 q^{5} - 5 q^{7},
  4 q^{4} + 10 q^{6} + 3 q^{8},
  \\
  & & \hspace{-2mm}
  - 5 q^{5} - 6 q^{7} - q^{9},
  3 q^{6} + q^{8}
  \\
  & & \hspace{-17mm}
  LG_{10_{143}}
  =
  21 + 64 q^{2} + 46 q^{4} + 10 q^{6},
  - 5 q^{-1} - 40 q^{1} - 56 q^{3} - 23 q^{5} - 2 q^{7},
  q^{-2} + 14 + 40 q^{2} + 28 q^{4} + 7 q^{6},
  \\
  & & \hspace{-2mm}
  - 3 q^{-1} - 17 q^{1} - 21 q^{3} - 8 q^{5} - q^{7},
  4 + 11 q^{2} + 5 q^{4} + q^{6},
  - 3 q^{1} - 3 q^{3},
  q^{2}
  \\
  & & \hspace{-17mm}
  LG_{10_{144}}
  =
  14 q^{-2} + 103 + 180 q^{2} + 84 q^{4} + 6 q^{6},
  - q^{-3} - 36 q^{-1} - 132 q^{1} - 124 q^{3} - 27 q^{5},
  \\
  & & \hspace{-2mm}
  3 q^{-2} + 42 + 88 q^{2} + 42 q^{4} + 3 q^{6},
  - 3 q^{-1} - 25 q^{1} - 27 q^{3} - 5 q^{5},
  1 + 6 q^{2} + 2 q^{4}
  \\
  & & \hspace{-17mm}
  LG_{10_{145}}
  =
  3 + 4 q^{2} - 2 q^{4} + 4 q^{6} + 4 q^{8},
  - 3 q^{1} + 2 q^{3} + 2 q^{5} - 4 q^{7} - q^{9},
  - q^{2} - 4 q^{4} - q^{6} + q^{8},
  q^{3} + q^{5},
  q^{4}
  \\
  & & \hspace{-17mm}
  LG_{10_{146}}
  =
  2 q^{-4} + 46 q^{-2} + 133 + 100 q^{2} + 24 q^{4},
  - 10 q^{-3} - 75 q^{-1} - 106 q^{1} - 45 q^{3} - 4 q^{5},
  16 q^{-2} + 54 + 38 q^{2} + 8 q^{4},
  \\
  & & \hspace{-2mm}
  - 11 q^{-1} - 16 q^{1} - 5 q^{3},
  3 + q^{2}
  \\
  & & \hspace{-17mm}
  LG_{10_{147}}
  =
  4 q^{-4} + 32 q^{-2} + 71 + 66 q^{2} + 14 q^{4},
  - 10 q^{-3} - 42 q^{-1} - 65 q^{1} - 35 q^{3} - 2 q^{5},
  10 q^{-2} + 31 + 36 q^{2} + 8 q^{4},
  \\
  & & \hspace{-2mm}
  - 5 q^{-1} - 14 q^{1} - 9 q^{3},
  1 + 3 q^{2}
  \\
  & & \hspace{-17mm}
  LG_{10_{148}}
  =
  29 + 92 q^{2} + 66 q^{4} + 12 q^{6},
  - 6 q^{-1} - 56 q^{1} - 79 q^{3} - 31 q^{5} - 2 q^{7},
  q^{-2} + 17 + 53 q^{2} + 38 q^{4} + 8 q^{6},
  \\
  & & \hspace{-2mm}
  - 3 q^{-1} - 20 q^{1} - 26 q^{3} - 10 q^{5} - q^{7},
  4 + 12 q^{2} + 6 q^{4} + q^{6},
  - 3 q^{1} - 3 q^{3},
  q^{2}
  \\
  & & \hspace{-17mm}
  LG_{10_{149}}
  =
  13 + 94 q^{2} + 156 q^{4} + 68 q^{6} + 4 q^{8},
  - q^{-1} - 36 q^{1} - 126 q^{3} - 113 q^{5} - 22 q^{7},
  4 + 52 q^{2} + 102 q^{4} + 48 q^{6} + 3 q^{8},
  \\
  & & \hspace{-2mm}
  - 7 q^{1} - 45 q^{3} - 49 q^{5} - 11 q^{7},
  7 q^{2} + 23 q^{4} + 12 q^{6} + q^{8},
  - 4 q^{3} - 5 q^{5} - q^{7},
  q^{4}
  \\
  & & \hspace{-17mm}
  LG_{10_{150}}
  =
  4 q^{-2} + 33 + 62 q^{2} + 48 q^{4} + 8 q^{6},
  - 11 q^{-1} - 44 q^{1} - 58 q^{3} - 26 q^{5} - q^{7},
  14 + 41 q^{2} + 41 q^{4} + 8 q^{6},
  \\
  & & \hspace{-2mm}
  - 12 q^{1} - 30 q^{3} - 19 q^{5} - q^{7},
  8 q^{2} + 16 q^{4} + 4 q^{6},
  - 4 q^{3} - 4 q^{5},
  q^{4}
  \\
  & & \hspace{-17mm}
  LG_{10_{151}}
  =
  28 q^{-2} + 135 + 178 q^{2} + 60 q^{4} + 2 q^{6},
  - 4 q^{-3} - 62 q^{-1} - 155 q^{1} - 112 q^{3} - 15 q^{5},
  \\
  & & \hspace{-2mm}
  11 q^{-2} + 70 + 105 q^{2} + 37 q^{4} + q^{6},
  - 13 q^{-1} - 48 q^{1} - 40 q^{3} - 5 q^{5},
  9 + 20 q^{2} + 7 q^{4},
  - 4 q^{1} - 4 q^{3},
  q^{2}
  \\
  & & \hspace{-17mm}
  LG_{10_{152}}
  =
  1 + 2 q^{2} + 10 q^{4} + 30 q^{6} + 18 q^{8} + 2 q^{10},
  - q^{1} - q^{3} - 16 q^{5} - 23 q^{7} - 7 q^{9},
  - q^{2} - 2 q^{4} + 10 q^{6} + 8 q^{8} + q^{10},
  \\
  & & \hspace{-2mm}
  4 q^{3} + 6 q^{5} + q^{7} - q^{9},
  - 5 q^{4} - 8 q^{6} - 4 q^{8},
  4 q^{5} + 5 q^{7} + q^{9},
  - q^{6},
  - q^{7} - q^{9},
  q^{8}
  \\
  & & \hspace{-17mm}
  LG_{10_{153}}
  =
  4 q^{-4} + 4 q^{-2} + 1 + 2 q^{2} + 4 q^{4},
  - q^{-5} - 4 q^{-3} + 2 q^{1} - 2 q^{3} - q^{5},
  q^{-4} - q^{-2} - 4 - 3 q^{2},
  \\
  & & \hspace{-2mm}
  q^{-3} + 4 q^{-1} + 3 q^{1},
  - 2 q^{-2} - 2 + 2 q^{2} + q^{4},
  q^{-1} - q^{1} - 2 q^{3},
  q^{2}
  \\
  & & \hspace{-17mm}
  LG_{10_{154}}
  =
  1 + 6 q^{2} + 26 q^{4} + 30 q^{6} + 18 q^{8} + 2 q^{10},
  - q^{1} - 9 q^{3} - 20 q^{5} - 19 q^{7} - 7 q^{9},
  - q^{2} + 2 q^{6} + 6 q^{8} + q^{10},
  \\
  & & \hspace{-2mm}
  4 q^{3} + 8 q^{5} + 3 q^{7} - q^{9},
  - 3 q^{4} - 4 q^{6} - q^{8},
  0,
  q^{6}
  \\
  & & \hspace{-17mm}
  LG_{10_{155}}
  =
  4 q^{-2} + 39 + 50 q^{2} + 24 q^{4} + 2 q^{6},
  - 19 q^{-1} - 45 q^{1} - 34 q^{3} - 8 q^{5},
  6 q^{-2} + 30 + 29 q^{2} + 11 q^{4} + q^{6},
  \\
  & & \hspace{-2mm}
  - q^{-3} - 14 q^{-1} - 20 q^{1} - 8 q^{3} - q^{5},
  3 q^{-2} + 11 + 5 q^{2},
  - 3 q^{-1} - 3 q^{1},
  1
  \\
  & & \hspace{-17mm}
  LG_{10_{156}}
  =
  24 q^{-2} + 91 + 102 q^{2} + 26 q^{4},
  - 4 q^{-3} - 46 q^{-1} - 98 q^{1} - 62 q^{3} - 6 q^{5},
  9 q^{-2} + 48 + 70 q^{2} + 24 q^{4} + q^{6},
  \\
  & & \hspace{-2mm}
  - 9 q^{-1} - 36 q^{1} - 32 q^{3} - 5 q^{5},
  7 + 18 q^{2} + 7 q^{4},
  - 4 q^{1} - 4 q^{3},
  q^{2}
  \\
  & & \hspace{-17mm}
  LG_{10_{157}}
  =
  17 + 134 q^{2} + 226 q^{4} + 102 q^{6} + 6 q^{8},
  - q^{-1} - 50 q^{1} - 181 q^{3} - 165 q^{5} - 33 q^{7},
  \\
  & & \hspace{-2mm}
  5 + 74 q^{2} + 146 q^{4} + 70 q^{6} + 4 q^{8},
  - 10 q^{1} - 64 q^{3} - 69 q^{5} - 15 q^{7},
  10 q^{2} + 31 q^{4} + 16 q^{6} + q^{8},
  \\
  & & \hspace{-2mm}
  - 5 q^{3} - 6 q^{5} - q^{7},
  q^{4}
  \\
  & & \hspace{-17mm}
  LG_{10_{158}}
  =
  6 q^{-4} + 100 q^{-2} + 217 + 122 q^{2} + 14 q^{4},
  - 31 q^{-3} - 150 q^{-1} - 162 q^{1} - 44 q^{3} - q^{5},
  \\
  & & \hspace{-2mm}
  4 q^{-4} + 54 q^{-2} + 115 + 61 q^{2} + 6 q^{4},
  - 9 q^{-3} - 45 q^{-1} - 46 q^{1} - 10 q^{3},
  8 q^{-2} + 20 + 8 q^{2},
  - 4 q^{-1} - 4 q^{1},
  1
  \\
  & & \hspace{-17mm}
  LG_{10_{159}}
  =
  45 + 144 q^{2} + 108 q^{4} + 20 q^{6},
  - 9 q^{-1} - 88 q^{1} - 127 q^{3} - 51 q^{5} - 3 q^{7},
  q^{-2} + 27 + 85 q^{2} + 62 q^{4} + 12 q^{6},
  \\
  & & \hspace{-2mm}
  - 4 q^{-1} - 32 q^{1} - 42 q^{3} - 15 q^{5} - q^{7},
  6 + 18 q^{2} + 9 q^{4} + q^{6},
  - 4 q^{1} - 4 q^{3},
  q^{2}
  \\
  & & \hspace{-17mm}
  LG_{10_{160}}
  =
  4 q^{-2} + 25 + 32 q^{2} + 14 q^{4},
  - 9 q^{-1} - 26 q^{1} - 23 q^{3} - 6 q^{5},
  8 + 19 q^{2} + 18 q^{4} + 3 q^{6},
  - 6 q^{1} - 18 q^{3} - 13 q^{5} - q^{7},
  \\
  & & \hspace{-2mm}
  6 q^{2} + 14 q^{4} + 4 q^{6},
  - 4 q^{3} - 4 q^{5},
  q^{4}
  \\
  & & \hspace{-17mm}
  LG_{10_{161}}
  =
  1 + 6 q^{4} + 8 q^{6} + 4 q^{8},
  - q^{3} - 5 q^{5} - 5 q^{7} - q^{9},
  - q^{2} + q^{8},
  q^{3} + 3 q^{5} + 2 q^{7},
  - q^{4} - 2 q^{6} - q^{8},
  0,
  q^{6}
  \\
  & & \hspace{-17mm}
  LG_{10_{163}}
  =
  14 q^{-2} + 87 + 138 q^{2} + 58 q^{4} + 4 q^{6},
  - q^{-3} - 33 q^{-1} - 106 q^{1} - 93 q^{3} - 19 q^{5},
  \\
  & & \hspace{-2mm}
  3 q^{-2} + 36 + 72 q^{2} + 33 q^{4} + 3 q^{6},
  - 3 q^{-1} - 22 q^{1} - 24 q^{3} - 5 q^{5},
  1 + 6 q^{2} + 2 q^{4}
  \\
  & & \hspace{-17mm}
  LG_{10_{164}}
  =
  42 q^{-2} + 195 + 246 q^{2} + 80 q^{4} + 2 q^{6},
  - 6 q^{-3} - 91 q^{-1} - 219 q^{1} - 154 q^{3} - 20 q^{5},
  \\
  & & \hspace{-2mm}
  16 q^{-2} + 100 + 147 q^{2} + 53 q^{4} + 2 q^{6},
  - 18 q^{-1} - 67 q^{1} - 57 q^{3} - 8 q^{5},
  12 + 27 q^{2} + 10 q^{4},
  - 5 q^{1} - 5 q^{3},
  q^{2}
  \\
  & & \hspace{-17mm}
  LG_{10_{165}}
  =
  4 q^{-4} + 94 q^{-2} + 245 + 172 q^{2} + 34 q^{4},
  - 23 q^{-3} - 148 q^{-1} - 192 q^{1} - 72 q^{3} - 5 q^{5},
  \\
  & & \hspace{-2mm}
  q^{-4} + 38 q^{-2} + 105 + 68 q^{2} + 11 q^{4},
  - 2 q^{-3} - 25 q^{-1} - 31 q^{1} - 8 q^{3},
  q^{-2} + 6 + 2 q^{2}
  \\
  & & \hspace{-17mm}
  LG_{10_{166}}
  =
  41 + 166 q^{2} + 164 q^{4} + 58 q^{6} + 4 q^{8},
  - 5 q^{-1} - 82 q^{1} - 149 q^{3} - 88 q^{5} - 16 q^{7},
  13 + 65 q^{2} + 61 q^{4} + 20 q^{6} + q^{8},
  \\
  & & \hspace{-2mm}
  - 11 q^{1} - 19 q^{3} - 9 q^{5} - q^{7},
  3 q^{2} + q^{4}
\end{eqnarray*}
\normalsize

\pagebreak

\bibliographystyle{plain}
\bibliography{DeWit99a}

\end{document}